\pdfoutput=1
\RequirePackage{ifpdf}
\ifpdf 
\documentclass[pdftex]{sigma}
\else
\documentclass{sigma}
\fi

\numberwithin{equation}{section}

\newtheorem{Theorem}{Theorem}[section]
\newtheorem{Corollary}[Theorem]{Corollary}
\newtheorem{Lemma}[Theorem]{Lemma}
\newtheorem{Proposition}[Theorem]{Proposition}
{ \theoremstyle{definition}
\newtheorem{Definition}[Theorem]{Definition}
\newtheorem{Remark}[Theorem]{Remark}}

\usepackage{mathrsfs}
\usepackage{tikz}
 \usetikzlibrary{cd}
\usepackage[all]{xy}
\usepackage{float}
\usepackage{enumitem}

\DeclareMathOperator{\tr}{tr}
\DeclareMathOperator{\Sym}{Sym}

\newcommand{\Z}{\mathbb{Z}}

\newcommand{\C}{\mathbb{C}}
\renewcommand{\P}{\mathbb{P}}

\renewcommand{\S}{\mathcal{S}}
\DeclareMathOperator{\Con}{Con}
\DeclareMathOperator{\Bun}{Bun}
\DeclareMathOperator{\Higgs}{Higgs}
\DeclareMathOperator{\Rep}{Rep}
\DeclareMathOperator{\Res}{Res}

\DeclareMathOperator{\Fix}{Fix}
\newcommand{\SP}{\Sym^2\big(\P^2\big)}
\newcommand{\bfw}{\mathbf{w}}
\newcommand{\bfW}{\mathbf{W}}
\newcommand{\bft}{\mathbf{t}}
\newcommand{\bfD}{\mathbf{D}}
\newcommand{\bfp}{\mathbf{p}}
\newcommand{\bfPsi}{\boldsymbol{\Psi}}
\newcommand{\bfnabla}{\boldsymbol{\nabla}}
\newcommand{\bfTheta}{\boldsymbol{\Theta}}
\newcommand{\bfmu}{\bar{\boldsymbol{\mu}}}
\newcommand{\bfnu}{\bar{\boldsymbol{\nu}}}

\begin{document}
\allowdisplaybreaks

\newcommand{\arXivNumber}{1910.13535}

\renewcommand{\PaperNumber}{125}

\FirstPageHeading

\ShortArticleName{A Map Between Moduli Spaces of Connections}

\ArticleName{A Map Between Moduli Spaces of Connections}

\Author{Frank LORAY~$^\dag$ and Valente RAM\'{I}REZ~$^\ddag$}
\AuthorNameForHeading{F.~Loray and V.~Ram\'{i}rez}

\Address{$^\dag$~Univ Rennes, CNRS, IRMAR - UMR 6625, F-35000 Rennes, France}
\EmailD{\href{mailto:frank.loray@univ-rennes1.fr}{frank.loray@univ-rennes1.fr}}

\Address{$^\ddag$~University of Twente, Department of Applied Mathematics, 7500 AE Enschede,\\
\hphantom{$^\ddag$}~The Netherlands}
\EmailD{\href{mailto:v.ramirez@utwente.nl}{v.ramirez@utwente.nl}}

\ArticleDates{Received December 17, 2019, in final form November 24, 2020; Published online December 02, 2020}

\Abstract{We are interested in studying moduli spaces of rank~2 logarithmic connections on elliptic curves having two poles. To do so, we investigate certain logarithmic rank 2 connections defined on the Riemann sphere and a transformation rule to lift such connections to an elliptic curve. The transformation is as follows: given an elliptic curve~$C$ with elliptic quotient $\pi\colon C\to\P^1$, and the logarithmic connection $(E,\nabla)$ on $\P^1$, we may pullback the connection to the elliptic curve to obtain a new connection $(\pi^*E, \pi^*\nabla)$ on $C$. After suitable birational modifications we bring the connection to a particular normal form. The whole transformation is equivariant with respect to bundle automorphisms and therefore defines a~map between the corresponding moduli spaces of connections. The aim of this paper is to describe the moduli spaces involved and compute explicit expressions for the above map in the case where the target space is the moduli space of rank~2 logarithmic connections on an elliptic curve~$C$ with two simple poles and trivial determinant.}


\Classification{14D20; 32G34; 34M55; 14H52; 53D30}

\section{Introduction} \label{sec:intro}

Let $X$ be a compact complex curve, $E$ a rank~2 holomorphic vector bundle, and $\nabla\colon E\to E\otimes\Omega^1_X(D)$ a connection having simple poles at the (reduced) divisor $D=t_1+\dots+t_n$. At each pole $t_i$, consider the residue matrix $\Res_{t_i}(\nabla)$ and denote by $\nu_i^+$, $\nu_i^-$ its eigenvalues.
Fixing the base curve $(X, D)$, the \textit{spectral data} $\bar\nu = (\nu_1^\pm, \dots, \nu_n^\pm)$, the \textit{trace connection} $(\det E, \tr\nabla)$, and introducing weights $\bar\mu$ for stability, we may construct the moduli space
$\Con^{\bar\mu}_{\bar\nu}(X,D)$
of $\bar\mu$-semistable $\bar\nu$-parabolic connections $(E,\nabla,\bar\ell)$ using geometric invariant theory (GIT)
\cite{InabaIwasakiSaito2006I, Nitsure1993}.
This moduli space is a separated irreducible quasi-projective variety of dimension $2N$, where $N=3g-3+n$ is the dimension of deformation of the base curve, and $g$ is the genus of $X$. This variety is moreover endowed with a holomorphic symplectic structure (which is in fact algebraic)
\cite{Boalch2001, Inaba2013, InabaIwasakiSaito2006I, Iwasaki1992}.

Moduli spaces of connections over the Riemann sphere have been extensively studied, in particular as these correspond to spaces of initial conditions for Garnier systems. The elliptic case with one and two poles have been studied in \cite{Loray2016} and \cite{FassarellaLoray2018}, respectively.

Closely related, we have moduli spaces $\Bun^{\bar\mu}(X,D)$ of $\bar\mu$-semistable \textit{parabolic bundles}, and a~natural map (which we denote $\operatorname{Bun}$) that assigns to a parabolic connection $(E,\nabla,\bar\ell)$ its underlying parabolic bundle $(E,\bar\ell)$. This correspondence is a Lagrangian fibration~\cite{LoraySaito2015}, and over the set of \textit{simple bundles} it defines an affine $\C^N$-bundle which is an affine extension of the cotangent bundle of $\Bun^{\bar\mu}(X,D)$
\cite{ArinkinLysenko1997A, ArinkinLysenko1997B}.

Let $(C,T)$ be an elliptic curve with two marked points, and let $\iota$ be the unique elliptic involution that permutes the marked points. Taking the quotient by this involution defines an elliptic covering $\pi \colon C \to \P^1$. Via this ramified covering we can pull bundles and connections from~$\P^1$ back to the elliptic curve~$C$. This correspondence, subject to some normalizations, defines a~map between the corresponding moduli spaces.
In this paper we aim to study a~particular map
\[
 \Phi \colon \ \Con^{\bar\mu}_{\bar\nu}\big(\P^1,D\big) \longrightarrow \Con^{\bfmu}_{\bfnu}(C,T),
\]
obtained in this way.
The divisor $D$ above contains 5 points: the four branch points of $\pi \colon C \to \P^1$ and the unique point $t\in\P^1$ which satisfies $\pi^{-1}(t) = T$ (see Section~\ref{subsec:Phi-connections} for the explicit construction).
This transformation was originally introduced in \cite{DiarraLoray2015}, using the associated monodromy representations. There, it was shown to be dominant and generically~$2:1$.
The same transformation rule induces also a map between moduli spaces of parabolic bundles (which we denote by a lowercase $\phi$), making the following diagram commute:
\begin{equation}\label{eq:Phi-square}
 \begin{tikzcd}
 \Con^{\bar\mu}_{\bar\nu}\big(\P^1,D\big) \arrow[r, "\Phi"] \arrow[d, swap, "\operatorname{Bun}"]
 &
 \Con^{\bfmu}_{\bfnu}(C,T) \arrow[d, "\operatorname{Bun}"]
 \\
 \Bun^{\bar\mu}\big(\P^1,D\big) \arrow[r, swap, "\phi"]
 &
 \Bun^{\bfmu}(C,T) .
 \end{tikzcd}
\end{equation}

The moduli spaces $\Con^{\bar\mu}_{\bar\nu}\big(\P^1,D\big)$ and $\Bun^{\bar\mu}\big(\P^1,D\big)$ have been explicitly described in~\cite{LoraySaito2015}, as well as the fibration $\operatorname{Bun}$ between them.
The moduli space of parabolic bundles $\Bun^{\bfmu}(C,T)$ was later studied in~\cite{FernandezVargas2016}. Moreover, the latter paper also describes geometrically and in coordinates the map $\phi \colon \Bun^{\bar\mu}\big(\P^1,D\big) \to \Bun^{\bfmu}(C,T)$.

The objective of this paper is to complete the explicit description of the commutative diagram~(\ref{eq:Phi-square}) by describing the space $\Con^{\bfmu}_{\bfnu}(C,T)$, endowing it with a coordinate system, and computing the map $\Phi \colon \Con^{\bar\mu}_{\bar\nu}\big(\P^1,D\big) \to \Con^{\bfmu}_{\bfnu}(C,T)$ in such coordinates.
In order to do so, we first study the associated map $\Phi^{{\rm top}} \colon \Rep_{\bar\nu}\big(\P^1,D\big)\to\Rep_{\bfnu}(C,T)$ between monodromy representations.

\subsection{Structure of the manuscript}

Later in the present section we will make some clarifications about notation and discuss related works.

In Section~\ref{sec:moduli-spaces} we recall general facts and definitions about parabolic bundles and connections, and about their moduli spaces.

In Section~\ref{sec:pullback-map} we define explicitly the transformation that takes a connection $\nabla$ on $\big(\P^1,D\big)$ and returns a connection $\Phi(\nabla)$ on $(C,T)$, thus defining the main object of study of the present paper: the map $\Phi \colon \Con^{\bar\mu}_{\bar\nu}\big(\P^1,D\big) \to \Con^{\bfmu}_{\bfnu}(C,T)$ between moduli spaces.
We also describe analogous transformations for parabolic bundles, and monodromy representations.
In this section we define the weights $\bar\mu$ and spectral data $\bar\nu$ that we will use throughout the present work.

The map $\Phi$ was originally defined for monodromy representations in \cite{DiarraLoray2015}.
In Section~\ref{sec:monodromy} we further discuss several properties of this map. Through the Riemann--Hilbert correspondence, we conclude that the map $\Phi$ between moduli spaces of connections enjoys analogous properties.

In Section~\ref{sec:genericity} we define and explain the genericity assumptions assumed in the statement of the main results.

The present work relies heavily on the constructions, results and ideas that appear in \cite{FernandezVargas2016, LoraySaito2015}. For the sake of the reader's convenience, we provide in Section~\ref{sec:recap} a brief survey of the results needed from these papers.

The main results of the present work are stated in Section~\ref{sec:main-body}. In Section~\ref{sec:geom-pic} we discuss some geometric properties of the map $\Phi$, before proving the main results in Section~\ref{sec:comp-coords}.
Finally, Section~\ref{sec:app} contains additional results, with the corresponding proofs, about the \textit{apparent map} in the elliptic case.

Below we present a short summary of our results.

Let $C\subset\P^2$ be an elliptic curve given by the affine equation $y^2 = x(x-1)(x-\lambda)$, and $\pi \colon C \to \P^1$ the elliptic quotient $(x,y) \mapsto x$. Let $t\in\P^1$ be a point different from $0,1,\lambda,\infty$. We fix the divisors $D=0+1+\lambda+\infty+t$ on $\P^1$, and $T=\pi^*(t)$ on $C$.

We begin with the map ${\Phi}^{{\rm top}} \colon \Rep_{\bar\nu}\big(\P^1,D\big)\to\Rep_{\bfnu}(C,T)$, which is the topological counterpart of our map $\Phi \colon \Con^{\bar\mu}_{\bar\nu}\big(\P^1,D\big) \to \Con^{\bfmu}_{\bfnu}(C,T)$ via the Riemann--Hilbert correspondence. It is proved in \cite[Theorem~1]{DiarraLoray2015} that this map is dominant and generically $2:1$. In Section~\ref{sec:monodromy} this result is extended by proving, in Theorem~\ref{thm:PhiMonodromySide},
that ${\Phi}^{{\rm top}}$ is surjective and a ramified cover between GIT spaces of representations; the ramification and branch loci of ${\Phi}^{{\rm top}}$ are described.
Moreover, we show in Theorem~\ref{thm:PhiMonodromySide} that ${\Phi}^{{\rm top}}$ is symplectic, up to a scalar factor of $2$ (the symplectic form on the codomain is pulled back to twice the symplectic form on the domain). As a consequence we obtain analogous results for the map $\Phi$ (cf.~Corollaries~\ref{cor:ConIsSingular} and~\ref{cor:PhiIsSymplectic}).

\begin{Remark}
 The fact that $\Phi$ is only symplectic up to a scalar factor is analogous to the case of the quadratic transformation of the Painlev\'{e} VI equation. This transformation is also induced by a ramified cover $\P^1\to\P^1$ of degree two, and the induced map between moduli spaces of connections is symplectic up to a factor of~2 \cite[Remark~3.1]{TsudaOkamotoSakai2005}.
\end{Remark}

In order to describe the map $\Phi$ in coordinates, we construct in Section~\ref{subsec:universal-family} a family of~connections over~$C$, denoted~$\mathcal{U}_C$, birationally parametrized by $\Con^{\bar\mu}_{\bar\nu}\big(\P^1,D\big) \stackrel{\sim}{\dasharrow} \Bun^{\bar\mu}\big(\P^1,D\big) \times \C^2$.
This family is the image under~$\Phi$ of the \textit{universal family} for $\Con^{\bar\mu}_{\bar\nu}\big(\P^1,D\big)$ constructed in \cite[Section~5]{LoraySaito2015} (which will be discussed in Section~\ref{subsec:universal-family-P1}).
We can choose a suitable set of generators $\bfnabla_0$, $\bfTheta_z$, $\bfTheta_w$ in such a way that any element $\bfnabla\in\mathcal{U}_C$ is given by a unique combination
\[
 \bfnabla = \bfnabla_0(u) + \kappa_1\bfTheta_z(u) + \kappa_2\bfTheta_w(u),
 \qquad
 u\in\Bun^{\bar\mu}\big(\P^1,D\big), \qquad (\kappa_1,\kappa_2)\in\C^2.
\]
The natural map into the moduli space $\mathcal{U}_C \dashrightarrow \Con^{\bfmu}_{\bfnu}(C,T)$ is a rational dominant map, generically $2:1$. Using this family we are able to give an explicit birational equivalence
\[
 \Con^{\bfmu}_{\bfnu}(C,T) \, \stackrel{\sim}{\dashrightarrow}\, \Bun^{\bfmu}(C,T) \times \C^2.
\]
This gives a trivialization of the affine $\C^2$-bundle $\Con^{\bfmu}_{\bfnu}(C,T) \to \Bun^{\bfmu}(C,T)$ over some open and dense subset of $\Bun^{\bfmu}(C,T)$.
Furthermore, over this dense set, it identifies the moduli space $\Con^{\bfmu}_{\bfnu}(C,T)$ to the moduli space of \emph{parabolic Higgs bundles} $\Higgs^{\bfmu}(C,T)$ (see~\cite{BiswasMukherjee2003}).
The latter is further identified with the cotangent bundle $T^* \Bun^{\bfmu}(C,T)$ in a natural way.
We can check {\it a posteriori} that these identifications are symplectic. The authors are not aware of a reference for this fact, that might be true in a more general setting.
This is why we make a large detour towards the Betti side of moduli spaces (i.e., representations) in order to deduce the symplectic structure from $\Con^{\bar\mu}_{\bar\nu}\big(\P^1,D\big)$ to $\Con^{\bfmu}_{\bfnu}(C,T)$.

Using the isomorphism $\Bun^{\bfmu}(C,T) \cong \P^1_z\times\P^1_w$ constructed in \cite[Section~4.3]{FernandezVargas2016}, we obtain a~coordinate system for the moduli space of connections
\[
 \Con^{\bfmu}_{\bfnu}(C,T) \stackrel{\sim}{\dashrightarrow} \P^1_z\times\P^1_w \times \C^2_{(\kappa_1, \kappa_2)}.
\]
We have explicitly computed the map $\Phi$ in these coordinates. Computations in coordinates appear throughout Section~\ref{sec:comp-coords}, and those corresponding to the map $\Phi$ are given in Section~\ref{subsec:base-change}.

Using the fact that the map $\Phi$ is symplectic up to a constant factor of~$2$, that is, $\Phi^*\omega_C = 2\omega_{\P^1}$, we show that the 2-form $\omega_C$ defining the symplectic structure of $\Con^{\bfmu}_{\bfnu}(C,T)$ is given, in the above coordinates, by
\[
 \omega_C = {\rm d}\kappa_1 \wedge {\rm d}z + {\rm d}\kappa_2 \wedge {\rm d}w,
\]
which coincides, under our identification, with the canonical 2-form defining the symplectic structure on $T^* \Bun^{\bfmu}(C,T)$ and $\Higgs^{\bfmu}(C,T)$.
This proves that the reduction from $\Con^{\bfmu}_{\bfnu}(C,T)$ to $\Higgs^{\bfmu}(C,T)$ is symplectic.

Unlike $\Con^{\bar\mu}_{\bar\nu}\big(\P^1,D\big)$, the moduli space $\Con^{\bfmu}_{\bfnu}(C,T)$ is singular. We describe in Section~\ref{subsec:singular-locus} the singular locus and describe the local analytic type of such singularities, together with its symplectic structure.

In Section~\ref{sec:app} we define an \emph{apparent map} for connections from the family $\mathcal{U}_C$. This map is defined as the set of tangencies of the connection with respect to \textit{two fixed subbundles}. The image of this map belongs to $\P^2\times\P^2$.
This map is not well defined on the moduli space, but after symmetrization, i.e., after passing to the quotient $\P^2\times\P^2\to\SP$, we obtain a well defined map which we denote $\operatorname{App}_C$.
Note that this is a map between spaces of the same dimension, thus not a Lagrangian fibration. The map is rational, dominant, and the generic fiber consists of exactly~12 points (cf.~Theorem~\ref{thm:App-fibers}).

Finally, inspired by the results of \cite{LoraySaito2015}, we combine the maps $\operatorname{App}$ and $\operatorname{Bun}$ to obtain a~generically injective map $\operatorname{App}_C\times\operatorname{Bun} \colon \Con^{\bfmu}_{\bfnu}(C,T) \dashrightarrow \SP\times\P^1_z\times\P^1_w$, showing that a~generic connection is completely determined by its underlying parabolic bundle together with its image under the apparent map (cf.~Theorem~\ref{thm:App-embedding}).

\subsection{Code repository}

All the computations mentioned in the present work have been carried out using the computer algebra system SageMath \cite{sage}. The code is available at the following repository \cite{github}.

\subsection{Related work}\label{sec:related-work}

It is well-known that compact Riemann surfaces of genus $g$ with $n$ punctures are hyperelliptic for
\[
 (g,n)=(2,0),\ (1,2),\ (1,1),\ \ \text{and}\ \ (1,0).
\]
It has been observed by W.~Goldman in \cite[Theorem~10.2]{Goldman97} that, $\mathrm{SL}_2(\mathbb C)$-representations of the fundamental group of these surfaces, with parabolic representation around each puncture, are invariant under the hyperelliptic involution; moreover, they come from the orbifold quotient representations.
From the Riemann--Hilbert correspondence, this means that a similar result should hold true for logarithmic connections, providing a dominant map between the correspon\-ding moduli spaces of connections.
This has been studied in details in the genus~2 case in \cite{HeuLoray2019}. The genus~1 case has been considered much earlier in \cite{Hitchin95} (see also \cite{LorayUlmerVanderPut2008}).
For the genus~1 case with one puncture, the same results also revealed to be true with arbitrary local monodromy at the puncture, which has been studied in \cite{Loray2016}.

The case studied here, 2~punctures on genus~1 curves, was first considered in \cite{DiarraLoray2015} for representations.
There it was proved that the result of Goldman recalled above, \cite[Theorem~10.2]{Goldman97}, extends as follows.
Consider the unique elliptic involution permuting the two punctures; then any $\mathrm{SL}_2(\mathbb C)$-representation whose image is Zariski dense, and whose boundary components have image into the same conjugacy class, is invariant under the involution and comes from a representation of the orbifold quotient.
The goal of the present paper was to provide the similar property for logarithmic connections, and therefore complete the whole picture for hyperelliptic curves. We note that similar constructions also hold within the class of connections on the 4-punctured sphere (see~\cite{MazzoccoVidunas13}).

The present work relies strongly on several results from \cite{FernandezVargas2016, LoraySaito2015}, which we discuss in Section~\ref{sec:recap}.

Finally, we remark the following for the 2-punctured elliptic curve case. Let $E$ be a rank 2 vector bundle over the elliptic curve $C$ of degree~$d$. By tensoring $E$ with a line bundle $L$, we can change the degree to any desired value as long as it has the same parity as~$d$. Therefore, the study of moduli spaces of rank~2 connections falls into two cases: odd degree and even degree. Usually the determinant of the bundle is fixed to be either~$\mathcal{O}_C$ in the even case (as in the present paper), or $\mathcal{O}_C(\bfw_\infty)$, where $\bfw_\infty\in\C$ is the identity element for the group structure of~$C$.
The moduli space of connections on $C$ with two poles and fixed determinant $\mathcal{O}_C(\bfw_\infty)$ has already been described in detail in \cite{FassarellaLoray2018}, together with its symplectic structure and apparent map. As pointed out in \cite{FernandezVargas2016}, it is possible to pass from the moduli space in the even degree case to that in the odd degree case. This is done by one elementary transformation followed by a~twist by a~rank~1 connection of degree zero. However, the transformation is not canonical, and this passage makes explicit computations hard to obtain.

\subsection{A note about notation}

In this text a curve will mean a nonsingular complex projective algebraic curve, which is identified to its associated analytic object, namely, a compact Riemann surface.
We will also identify a vector bundle with its associated locally free sheaf.

We are going to deal with a lot of objects that are defined over the elliptic curve $C$, and analogous objects defined over $\P^1$. In order to avoid confusion, we will try to use bold typography for objects in $C$ that have a counterpart in $\P^1$ (e.g., $\bfnabla$ and~$\nabla$).

Throughout this work we will use $\Phi$ to denote the transformation described in Section~\ref{subsec:Phi-connections}, which takes a connection defined over $\big(\P^1,D\big)$, and returns a connection over $(C,T)$. We use the same symbol for the analogous transformation acting on parabolic Higgs bundles.
The corresponding map between parabolic bunldes is denoted by $\phi$, and that between monodromy representations $\Phi^{{\rm top}}$.
All these maps are generically $2:1$, and we denote the corresponding Galois involutions, permuting the two sheets of the cover, by $\Psi$, $\psi$ and $\Psi^{{\rm top}}$, respectively.

Finally, we remark that we write $\P^1_z$ whenever we want to make explicit the fact that the space $\P^1$ is endowed with an affine coordinate $z\in\C$. This will allow us to distinguish different occurrences of $\P^1$. Similar for affine spaces such as $\C^2_{(c_1,c_2)}$.

\section{General aspects about parabolic bundles and connections} \label{sec:moduli-spaces}

Let $X$ be a smooth projective complex curve and $D = t_1 + \dots + t_n$ a reduced divisor.
A \emph{quasi-parabolic bundle} of rank 2 on $(X, D)$ is a pair $(E, \bar\ell)$, where $E$ is a holomorphic vector bundle of rank 2 over $X$, and $\bar\ell = \{\ell_1, \dots, \ell_n\}$ a collection of rank 1 subspaces $\ell_i \subset E|_{t_i}$.
A \emph{parabolic bundle} is a quasi-parabolic bundle endowed with a vector of weights $\bar\mu = (\mu_1,\dots,\mu_n)$, where $\mu_i\in[0,1]$. We will usually omit the vector $\bar\mu$ in the notation and denote a parabolic bundle simply by $(E,\bar\ell)$.

A \emph{logarithmic connection} on $X$ with poles at $D$ is a pair $(E,\nabla)$, where $E$ is a holomorphic vector bundle over $X$, and $\nabla \colon E \to E \otimes \Omega^1_X(D)$ is a $\C$-linear map satisfying Leibniz' rule.
The eigenvalues of the residue $\Res_{t_i}(\nabla)$, $\nu_i^+$, $\nu_i^-$ are called the \emph{local exponents}, and the collection $\bar\nu = (\nu_1^\pm, \dots, \nu_n^\pm)$ is the \textit{spectral data}.
We have the following equality known as \emph{Fuchs' relation}:
\begin{gather}\label{eq:Fuchs}
 \sum_{i=1}^n (\nu_i^+ + \nu_i^-) + \deg E = 0.
\end{gather}

\begin{Definition} Let $\bar\nu$ be fixed spectral data and $\bar\mu$ a fixed vector of weights. A $\bar\nu$-\emph{parabolic connection} of rank~2 on $(X,D)$ is a triple $(E, \nabla, \bar\ell)$ where $(E, \bar\ell)$ is a rank~2 parabolic bundle and $(E, \nabla)$ is a logarithmic connection with poles on~$D$, such that at each subspace~$\ell_i$ the residue $\Res_{t_i}(\nabla)$ acts by multiplication by $\nu_i^+$.
\end{Definition}

We remark that the difference between two connections is an $\mathcal{O}_X$-linear operator (known as a~\emph{Higgs field}), and the space of connections with a fixed parabolic structure is a finite dimensional affine space.

\begin{Definition}
 A \emph{parabolic Higgs bundle} of rank 2 on $(X,D)$ is a triple $(E, \Theta, \bar\ell)$, where $(E, \bar\ell)$ is a parabolic vector bundle, $\Theta \colon E \to E \otimes \Omega^1_X(D)$ is a $\mathcal{O}_X$-linear map, and such that for each $t_i\in D$ the residue $\Res_{t_i}(\Theta)$ is nilpotent with null space given by $\ell_i$.
\end{Definition}

We now introduce the notion of $\bar\mu$-semistability.

\begin{Definition}\label{def:stability}
 Let $(E, \bar\ell)$ be a rank 2 parabolic bundle and $\bar\mu = (\mu_1, \dots, \mu_n) \in [0,1]^n$ its weight vector. We define the $\bar\mu$-\emph{parabolic degree} of a line subbundle $L \subset E$ as
 \[
 \deg E - 2\deg L + \sum_{\ell_i \not\subset L} \mu_i - \sum_{\ell_i \subset L} \mu_i .
 \]
 The parabolic bundle $(E, \bar\ell)$ is said to be $\bar\mu$-\emph{semistable} ($\bar\mu$-\emph{stable}) if the parabolic degree is non-negative (resp.~positive) for every subline bundle $L$.

 A parabolic connection $(\nabla, E, \bar\ell)$ is said to be $\bar\mu$-\emph{semistable} ($\bar\mu$-\emph{stable}) if the parabolic degree is non-negative (resp.~positive) for every subline bundle invariant by the connection.
\end{Definition}

In order to define moduli spaces it is convenient to fix the determinant bundle $\det(E)$, and the trace $\tr(\nabla)$ in the case of connections. These choices will not appear explicitly in the notation, but we always assume these objects have been defined and fixed.
The moduli space does not depend on the choice of the prescribed determinant, as we can
freely change it by twisting by a~line bundle and performing elementary transformations. This is further explained in Section~\ref{subsec:elem-transf}.

\begin{Definition}\label{def:Con}
 We denote by $\Con^{\bar\mu}_{\bar\nu}(X,D)$ the moduli space of $\bar\mu$-semistable $\bar\nu$-parabolic connections on $(X,D)$, where the determinant bundle and trace connection equal some fixed pair~$(L, \eta)$.
 Similarly, we denote by $\Higgs^{\bar\mu}(X,D)$ the moduli space of $\bar\mu$-semistable parabolic Higgs bundles with given trace and determinant bundle, and by $\Bun^{\bar\mu}(X, D)$ the moduli space of $\bar\mu$-semistable parabolic bundles with a given determinant bundle.
\end{Definition}

From now on, unless otherwise specified, connections, Higgs fields and bundles are assumed to have trivial determinant and zero trace.

\begin{Remark}
 The precise spectral data and weight vectors to be used throughout this text are given in Definition~\ref{def:mu-nu}.
 In Lemma~\ref{lemma:definition-Phi} it is proved that, for such spectral data, any $\bar\nu$-parabolic connection is $\bar\mu$-semistable.
 Therefore, the moduli space of $\bar\mu$-semistable $\bar\nu$-parabolic connections defined above coincides
 with the moduli space $\Con_{\bar\nu}(X,D)$ of $\bar\nu$-parabolic connections.
\end{Remark}

\begin{Remark} \label{rmk:unstable-bundles}
 As can be seen in Definition~\ref{def:stability}, a parabolic connection is semistable if and only if every subbundle \textit{invariant by the connection} has non-negative parabolic degree.
 Thus it is possible for a connection to be semistable while the underlying parabolic bundle is not.
 These are exceptional cases and will be excluded from our definition of a \textit{generic connection}.
\end{Remark}

\subsection{Twists and elementary transformations}\label{subsec:elem-transf}

In the previous section we have defined the moduli spaces of parabolic bundles and connections by restricting ourselves to bundles with a prescribed determinant. By imposing this condition we do not loose any generality. Indeed, the determinant bundle can be arbitrarily modified by a suitable sequence of \textit{elementary transformations} and \textit{twists} by rank~1 connections. In this section we briefly describe these transformations for the special case of rank~2.

A twist by a line bundle $L$ is simply the transformation $(E,\bar\ell) \mapsto (E\otimes L, \bar\ell)$. Twisting by a~rank~1 connection is defined as follows.

\begin{Definition}\label{def:twist} Let $(E,\nabla,\bar\ell)$ be a parabolic connection on $(X,D)$ with spectral data $\bar\nu$.
 The \textit{twist} by a rank~1 connection $(L,\xi)$ on $(X,D)$ with spectral data $\bar\vartheta$ is defined to be the parabolic connection
 \[
 (E,\nabla,\bar\ell) \otimes (L, \xi) = (E\otimes L, \nabla\otimes\xi, \bar\ell).
 \]
\end{Definition}
The spectral data over a pole $p\in D$ transforms as $(\nu^+, \nu^-) \mapsto (\nu^+ + \vartheta, \nu^- + \vartheta)$.

Let us now define the elementary transformations.
We refer the reader to \cite[Section~6]{Machu2007} for a detailed discussion and equivalent ways of defining these transformations.

\begin{Definition}\label{def:elem-transf}
 Let $(E, \bar\ell)$ be a parabolic bundle over $(X,D)$, and consider a point $p\in D$.
 The parabolic structure $\bar\ell$ provides a rank~1 subspace of $\ell_p \subset E\vert_p$. Let us regard $\ell_p$ as a sky-scrapper sheaf supported at $p$.
 We define a new parabolic bundle $(E^-, \bar\ell') = \operatorname{elem}^-_p(E, \bar\ell)$ as follows. The underlying vector bundle is characterized by the exact sequence of sheaves
 \[
 0 \longrightarrow E^- \longrightarrow E \longrightarrow E\vert_p/\ell_p \longrightarrow 0.
 \]
 The parabolic structure is unchanged outside of $p$. The new parabolic direction $\ell^-_p \subset E^-\vert_p$ is defined by the kernel of $E\vert_p^- \to E\vert_p$.
 We remark that $\det(E^-) = \det(E)\otimes\mathcal{O}_X(-p)$.
 Thus this is called a \textit{negative elementary transformation}.

 \noindent We define $(E^+, \bar\ell') = \operatorname{elem}^+_p(E, \bar\ell)$ as a twist by the line bundle $\mathcal{O}_X(p)$ followed by a negative elementary transformation at $p$. Thus $E^+$ fits into the exact sequence
 \[
 0 \rightarrow E^+ \rightarrow E(p) \rightarrow E\vert_p/\ell_p \rightarrow 0.
 \]
 The parabolic direction $\ell^+_p \subset E^+\vert_p$ is once again defined by the kernel of $E\vert_p^+ \to E(p)\vert_p$.
 This time $\det(E^+) = \det(E)\otimes\mathcal{O}_X(p)$, and we call this a \textit{positive elementary transformation}.
\end{Definition}

From the perspective of the ruled surface $\P(E)$, the parabolic structure $\ell_p$ is nothing but a point on the projectivized fiber $F=\P(E\vert_p)$. The elementary transformation $\operatorname{elem}^+_p(E, \bar\ell)$ corresponds to the birational transformation of the total space $\P(E)$ given by blowing-up the point $\ell_p \in F$, and then contracting the strict transform of the fiber $F$. The point resulting from this contraction gives the new parabolic direction.
Since the two transformations $\operatorname{elem}^-_p(E, \bar\ell)$ and $\operatorname{elem}^+_p(E, \bar\ell)$ differ only by twisting by a line bundle, they coincide projectively. Namely, they define the same birational transformation on the ruled surface $\P(E)$.

Finally, we remark that a connection $\nabla$ on $E$ induces a connection $\nabla'$ on the subsheaf $E^- \subset E$. Over~$p$ the residual eigenvalues are changed by the rule $(\nu^+, \nu^-) \mapsto (\nu^- + 1, \nu^+)$.
This means that the eigenvalue associated to the parabolic direction $\ell'_p$ is now $\nu^- + 1$.
In a positive elementary transformation we must twist by a rank~1 connection on $\mathcal{O}_X(p)$ having a unique pole at $p$ with eigenvalue~$-1$.
Therefore, under a positive elementary transformation, the eigenvalues are transformed as $(\nu^+, \nu^-) \mapsto (\nu^-, \nu^+ - 1)$.

\section{The pullback map} \label{sec:pullback-map}

Let $C\subset\P^2$ be an elliptic curve such that in some fixed affine chart it is given by the equation
\begin{equation*}
 y^2 = x(x-1)(x-\lambda), \qquad \lambda\in\C\setminus \{0,1\}.
\end{equation*}
This curve is endowed with the \textit{elliptic involution} $\iota \colon (x,y)\mapsto (x,-y)$. With respect to the group structure of $C$, this involution is precisely $\bfp \mapsto -\bfp$. The quotient of $C$ under this involution gives rise to the \textit{elliptic quotient} $\pi \colon C \to \P^1$. This is a $2:1$ cover ramified over the 2-torsion points $\bfw_0$, $\bfw_1$, $\bfw_\lambda$, $\bfw_\infty$, which are the points on $C$ that satisfy $x=0,1,\lambda,\infty$, respectively.

Let us choose a point $t\in\P^1\setminus \{0,1,\lambda,\infty\}$, and let $\pi^{-1}(t) = \{\bft_1, \bft_2\}$. We define the following divisors of~$\P^1$:
\[
 W = 0 + 1 + \lambda + \infty, \qquad D = W + t.
\]
We define analogous divisors for $C$:
\[
 \bfW = \bfw_0 + \bfw_1 + \bfw_\lambda + \bfw_\infty, \qquad T = \bft_1 + \bft_2, \qquad \bfD=\bfW+ T.
\]

Now, let us fix the spectral data and weights to use throughout the text. We remark that for the most part we will work with $\mathfrak{sl}_2(\C)$-connections. Therefore the spectral data will always satisfy $\nu^-_i = -\nu^+_i$.

\begin{Definition}\label{def:mu-nu}
 Let $\nu$ any complex number such that $2\nu\not\in\Z$, and choose $\mu$ a real number $0<\mu<1$.
 When working with parabolic bundles over $(C,T)$, we define the spectral data $\bfnu = (\pm\nu, \pm\nu)$ and the weight vector $\bfmu = (\mu, \mu)$.
 For working with bundles over $\big(\P^1,D\big)$, we define the vectors
 $\bar\nu = \big(\pm\frac{1}{4}, \pm\frac{1}{4}, \pm\frac{1}{4}, \pm\frac{1}{4}, \pm\nu\big)$
 and $\bar\mu = \big(\frac{1}{2}, \frac{1}{2}, \frac{1}{2}, \frac{1}{2}, \mu\big)$.
\end{Definition}

\subsection[The map phi on parabolic bundles]{The map $\boldsymbol{\phi}$ on parabolic bundles} \label{subsec:Phi-bundles}

The map $\phi \colon \Bun^{\bar\mu}\big(\P^1,D\big) \to \Bun^{\bfmu}(C,T)$ was originally introduced in \cite[Section~6.1]{FernandezVargas2016} (there it is denoted by an uppercase $\Phi$). We repeat here the construction for the reader's convenience. We refer the reader to the aforementioned paper for extra details and proofs.

Consider a semistable parabolic bundle $(E, \bar\ell)$ on $\big(\P^1,D\big)$, with weight vector as defined in Definition~\ref{def:mu-nu}, namely, $\bar\mu= \big(\frac{1}{2}, \frac{1}{2}, \frac{1}{2}, \frac{1}{2}, \mu\big)$, for some $0<\mu<1$.
We recall that through this text we deal with degree zero bundles, thus $\det(E)=\mathcal{O}_{\P^1}$.
In order to define the image $\phi(E,\bar\ell)$ we proceed as follows.
\begin{enumerate}\itemsep=0pt
 \item Pullback $(E,\bar\ell)$ to $C$ using the elliptic cover $\pi\colon C\to \P^1$. This defines a bundle $E'=\pi^*E$ with parabolic structure supported over $\bfD=\bfW+T$. Such bundle is semistable with respect to the weights $\bar\mu'=(1,1,1,1,\mu,\mu)$.
 \item Perform a positive elementary transformation (cf.~Definition~\ref{def:elem-transf}) at each of the parabolic points in the divisor $\bfW=\bfw_0 + \bfw_1 + \bfw_\lambda + \bfw_\infty$. This defines a new bundle $E''$ of degree 4 (in fact $\det(E'')=\mathcal{O}_C(\bfW) = \mathcal{O}_C(4\bfw_\infty)$). The parabolic bundle is semistable with respect to $\bar\mu''=(0,0,0,0,\mu,\mu)$.
 \item Tensor the previous bundle with the line bundle $\mathcal{O}_C(-2\bfw_\infty)$. This new parabolic bundle has trivial determinant, and continues to be semistable for $\bar\mu''=(0,0,0,0,\mu,\mu)$.
 \item Because of the nullity of the weights over $\bfW$, we may simply forget the parabolic structure over $\bfW$ to recover a parabolic bundle with trivial determinant over $(C,T)$, which is semistable for the weight vector $\bfmu = (\mu,\mu)$.
\end{enumerate}

We denote the final parabolic bundle by $\phi(E,\bar\ell)$. This defines a transformation from parabolic bundles of degree zero over $\big(\P^1,D\big)$, to parabolic bundles with trivial determinant over $(C,T)$. The weight vectors are precisely those in Definition~\ref{def:mu-nu}. This transformation induces a map between the corresponding moduli spaces, which is also denoted $\phi$.

\subsection[The map Phi on connections]{The map $\boldsymbol{\Phi}$ on connections} \label{subsec:Phi-connections}

The same construction, with minor additions, can be adapted to define a transformation between parabolic connections. Below we make such steps explicit.

Let $\nabla$ be a rank 2 connection on $\P^1$, defined over a degree-zero bundle $E$, having simple poles over the divisor~$D$, and with spectral data given by $\bar\nu=\big( \pm\frac{1}{4}, \pm\frac{1}{4}, \pm\frac{1}{4}, \pm\frac{1}{4}, \pm\nu\big)$.
The following series of transformations defines the map~$\Phi$.
We remark that the underlying bundles and weight vectors are the same as in Section~\ref{subsec:Phi-bundles}.
\begin{enumerate}\itemsep=0pt
 \item \label{item:i}Pullback $\nabla$ to $C$ using $\pi$. This gives a connection $\pi^*\nabla$ on $E' = \pi^*E$ with poles on $\bfD$. Locally, the connection near $\bft_1$, $\bft_2$ looks like $\nabla$ around~$t$. This is not the case around the ramification points~$\bfw_k$, but we know this construction multiplies the residual eigenvalues by a factor of two. Therefore the spectral data is given by $\big( \pm\frac{1}{2}, \pm\frac{1}{2}, \pm\frac{1}{2}, \pm\frac{1}{2}, \pm\nu, \pm\nu\big)$.

 \item \label{item:ii}Perform a positive elementary transformation (cf.~Definition~\ref{def:elem-transf}) for each pole in the divisor $\bfW$. This gives a new connection on some bundle $E''$ of degree 4. The spectral data over the points $\bft_i$ is unchanged, and the new spectral data at the $\bfw_k$ is $\nu^+_k=-\frac{1}{2}$, $\nu_k^-=-\frac{1}{2}$ (not an $\mathfrak{sl}_2(\C)$-connection).

 \item \label{item:iii} Tensor with the rank 1 connection $(\mathcal{O}_C(-2\bfw_\infty), \xi)$, where $\xi$ is a suitable rank~1 connection with simple poles on $\bfW$ having residue $\frac{1}{2}$ at each of them (no poles on $T$)\footnote{ The connection $\xi$ on $\mathcal{O}_C(-2\bfw_\infty)$ is the pullback under the elliptic cover $\pi\colon C \to \P^1$ of the rank 1 connection on $\mathcal{O}_{\P^1}(-\infty)$ given by ${\rm d} + \frac{1}{4}\big(\frac{{\rm d}P}{P}\big)$, where $P(x)=x(x-1)(x-\lambda)$. The former connection has residue $\frac{1}{4}$ at $x=0,1,\lambda,\infty$. Therefore $\xi$ has residue $\frac{1}{2}$ at the branch points $\bfw_0$, $\bfw_1$, $\bfw_\lambda$, $\bfw_\infty$, as desired. }.
 By construction, the bundle $E''\otimes \mathcal{O}_C(-2\bfw_\infty)$ has trivial determinant and the residual eigenvalues at $\bfw_k$ are all zero.
 The fact that the poles of the original connection over $\P^1$ have semisimple monodromy around each pole in $W$ implies that this new connection over $C$ has trivial monodromy around each point $\bfw_k$. Thus this connection is in fact holomorphic at each point in $\bfW$.

 \item \label{item:iv} Since the final connection is holomorphic at $\bfW$, we may forget these points from the divisor of poles and consider it as a connection defined on $(C,T)$ with spectral data $\bar\nu = (\pm\nu, \pm\nu)$.
\end{enumerate}

We denote the last connection by~$\Phi(\nabla)$.

\begin{Lemma} \label{lemma:definition-Phi}
 Let $\nabla$ be any connection on $\big(\P^1,D\big)$ with spectral data $\bar\nu$. Then $\Phi(\nabla)$ is $\bar\mu$-semistable $($and has spectral data~$\bar\nu)$.
 Therefore, the correspondence $\nabla \mapsto \Phi(\nabla)$ induces a map
 \[
 \Phi \colon \ \Con^{\bar\mu}_{\bar\nu}\big(\P^1,D\big) \longrightarrow \Con^{\bfmu}_{\bfnu}(C,T).
 \]
 We will use the same notation for the geometric transformation defined by steps {\rm (\ref*{item:i})--(\ref*{item:iv})} and the induced map between moduli spaces.
\end{Lemma}

\begin{proof}First, we remark that, by definition, a connection can only be unstable if it is reducible.
 By our choice of spectral data~$\bar\nu$, any connection $\nabla$ on $\big(\P^1,D\big)$ is irreducible and therefore $\bar\mu$-stable.
 On the other hand, $\Phi(\nabla)$ (or any connection having spectral data $\bfnu$ on~$C$) must be $\bfmu$-semistable.
 To see this, assume that $\Phi(\nabla)$ has an invariant line bundle~$L$.
 Because $2\nu\not\in\Z$, Fuchs' relation~(\ref{eq:Fuchs}) implies that the restriction to the invariant line bundle $L$ must have spectral data $\{+\nu,-\nu\}$, and $\deg(L)=0$, and so the $\bfmu$-parabolic degree of $L$ (cf.~Definition~\ref{def:stability}) is $\mu-\mu = 0$.
 This implies that the connection (and also the underlying parabolic bundle) is strictly $\bfmu$-semistable in the reducible case.
\end{proof}

\begin{figure}[t] \centering
 \includegraphics[width=\linewidth]{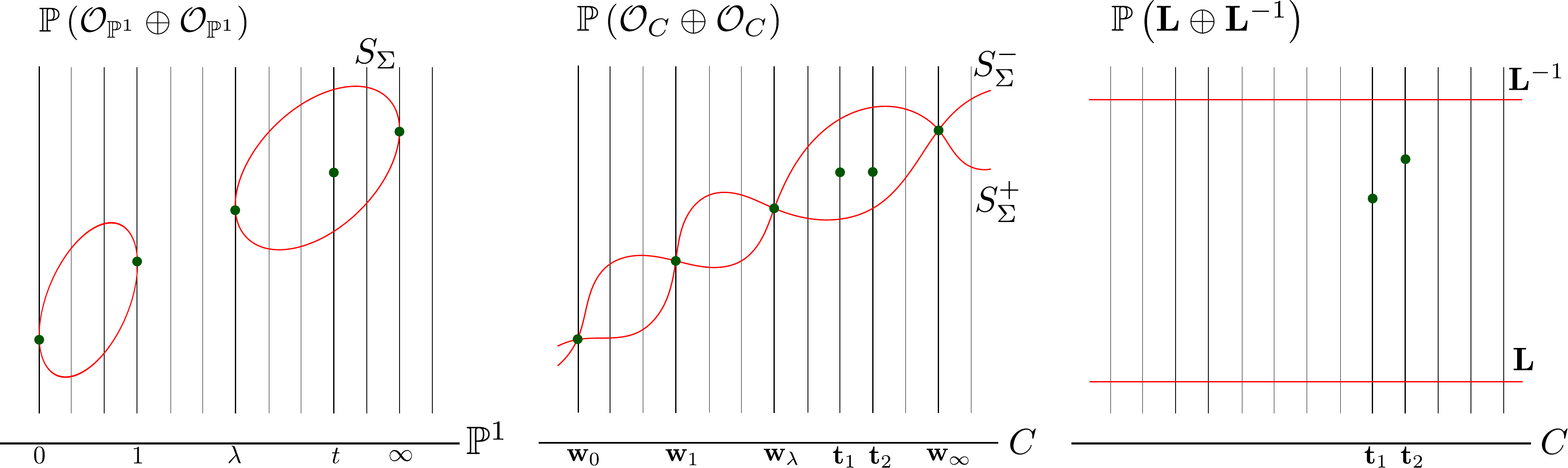}
 \caption{Steps of the transformation~$\Phi$. The canonical sections corresponding to $\mathbf{L}$, $\mathbf{L}^{-1}$ in $\P\big(\mathbf{L} \oplus \mathbf{L}^{-1}\big)$ come from a~\emph{multisection} $S_\Sigma$ in $\P\big(\mathcal{O}_{\P^1} \oplus \mathcal{O}_{\P^1}\big)$. This is further explained in Remark~\ref{rmk:multisection-S_Sigma}.} \label{fig:Phi-hat}
\end{figure}

We remark that a typical element of $\Bun^{\bar\mu}\big(\P^1,D\big)$ has a trivial underlying bundle (as will be discussed in Section~\ref{sec:genericity}).
On the other hand, a typical element of $\Bun^{\bfmu}(C,T)$ (for example a~$\bfmu$-stable parabolic bundle) is such that its underlying bundle can be written as $E'' = \mathbf{L} \oplus \mathbf{L}^{-1}$, where $\mathbf{L}$ is a rank~1 bundle on $C$ of degree zero \cite[Proposition~4.5]{FernandezVargas2016}.
For such case, Fig.~\ref{fig:Phi-hat} shows the effect of steps (\ref*{item:i}) and (\ref*{item:ii}) on the projectivization of the bundles involved.

\subsection[The map Phi-top on monodromy representations]{The map $\boldsymbol{\Phi^{{\rm top}}}$ on monodromy representations}
 \label{subsec:Phi-rep}

The map $\Phi$ is also defined in a more topological setting.
This was originally introduced in \cite{DiarraLoray2015}. We repeat here the construction presented in Section~4 of the cited paper.
Below we continue to use the spectral data $\bar\nu$ and $\bfnu$ introduced in Definition~\ref{def:mu-nu}.

Let us define $\Rep_{\bar\nu}\big(\P^1,D\big)$ as the subspace of
\[
 \Rep\big(\P^1,D\big) := \operatorname{Hom}\big(\pi_1\big(\P^1 \setminus D\big), \, \mathrm{SL}_2(\C)\big) \,/\, \mathrm{SL}_2(\C),
\]
consisting of those representations which are compatible with the spectral data $\bar\nu$ around the punctures $D$.
Explicitly, we define
\[
 \Rep_{\bar\nu}\big(\P^1,D\big):=\left\{(M_0,M_1,M_t,M_\lambda,M_\infty)\in\mathrm{SL}_2(\C)^5 ;\
 \begin{matrix}M_0M_1M_tM_\lambda M_\infty=I\\
 \mathrm{trace}(M_i)=0\ \text{for}\ i\not=t\\
 \mathrm{trace}(M_t)=2\cos(2\pi\nu)
 \end{matrix}\right\}/\sim \, ,
\]
where two representations are equivalent, $(M_i)\sim(M_i')$, if and only if there exists $M\in\mathrm{SL}_2(\C)$ such that
$M_i'=MM_iM^{-1}$ for all $i\in D$.
These representations are always irreducible, therefore the space defined above coincides with the GIT quotient of the action of $\mathrm{SL}_2(\C)$ on $\operatorname{Hom}\big(\pi_1\big(\P^1 \setminus D\big), \mathrm{SL}_2(\C)\big)$.
In this way, it admits a structure of smooth irreducible affine variety of complex dimension $4$.
In a similar fashion, we define for the twice-punctured torus
\[
 \widetilde{\Rep}_{\bfnu}(C, T):=\left\{(A,B,C_1,C_2)\in\mathrm{SL}_2(\C)^4 ;\
 \begin{matrix}AB=C_1 BA C_2\\
 \mathrm{trace}(C_i)=2\cos(2\pi\nu)
 \end{matrix}\right\}/\sim \, .
\]
This quotient however is not Hausdorff, as there are reducible representations; it is an algebraic stack.
A Hausdorff quotient
\[ \Pi_{{\rm GIT}}\colon \ \widetilde{\Rep}_{\bfnu}(C, T)\twoheadrightarrow\Rep_{\bfnu}(C,T),\]
is obtained by GIT, and $\Rep_{\bfnu}(C,T)$ is an irreducible affine variety of complex dimension $4$.
This quotient $\Pi_{{\rm GIT}}$ is obtained by identifying triangular representations with their diagonal part;
this gives rise to a singular locus for the affine variety $\Rep_{\bfnu}(C,T)$.
We now define a map between the above two spaces.
Given a representation $(M_i)\in\Rep_{\bar\nu}\big(\P^1,D\big)$, we can perform the following.
\begin{enumerate}\itemsep=0pt
 \item Pullback $(M_i)$ via the elliptic cover $\pi \colon C \to \P^1$. This defines a representation $\pi^*(M_i)$ of the fundamental group $\pi_1(C\setminus \bfW\cup T)$ with local monodromy $-I$ at the $\bfW$-punctures.
 \item Twist $\pi^*(M_i)$ by the central representation $\pi_1(C\setminus \bfW)\to\{\pm I\}$ with local monodromy $-I$ at each of the punctures.
 \item This last representation has trivial monodromy around the $\bfW$-punctures, and so we can regard it as an element of $\widetilde{\Rep}_{\bfnu}(C,T)$, or better of the quotient $\Rep_{\bfnu}(C,T)$.
\end{enumerate}

The above procedure provides a well-defined map
\[
 \Phi^{{\rm top}} \colon \ \Rep_{\bar\nu}\big(\P^1,D\big)\to\Rep_{\bfnu}(C,T),
\]
and is the analogue of the map $\Phi$ on the monodromy side. This claim is justified by the next remark.

\begin{Remark} \label{rmk:rep-con-equivalence}
 Let $(E,\nabla, \bar\ell)$ be a parabolic rank~2 connection on a marked curve $(X,D)$. An elementary transformation over a pole $p\in D$ does not change the monodromy representation associated with $\nabla$ (this defines an isomorphism outside the poles, where the monodromy is computed).
 On the other hand, twisting $\nabla$ by a rank~1 connection with residue $\frac{1}{2}$ at $p \in D$, scales the local monodromy around $p$ by a factor of $-1$.
 Indeed, with respect to a suitable local coordinate $z$ for which $p$ is given by $z=0$, the connection $\nabla$ is given by $\nabla = {\rm d} + A \frac{{\rm d}z}{z}$, for a~constant matrix~$A$.
 Such a twist transforms the connection into $\nabla' = {\rm d} + \big(A +\frac{1}{2}I\big) \frac{{\rm d}z}{z}$. Since the matrix $\frac{1}{2}I$ is scalar, it evidently commutes with~$A$, and so the local monodromy is multiplied by the matrix: $\operatorname{exp} \big(2\pi{\rm i} \frac{1}{2}I\big) = -I$.
\end{Remark}

The steps defining the map $\Phi$ in Section~\ref{subsec:Phi-connections} can now be clearly matched to those used to define $\Phi^{{\rm top}}$ in the present section.
Moreover, the Riemann--Hilbert correspondence provides us with a commutative diagram
\begin{equation} \label{eq:RH}
\begin{split}&
 \xymatrix{
 \Con^{\bar\mu}_{\bar\nu}\big(\P^1,D\big) \ar[d]_{\Phi} \ar[r]^{\rm RH} & \Rep_{\bar\nu}\big(\P^1,D\big) \ar[d]^{{\Phi}^{{\rm top}}} \\
 \Con^{\bfmu}_{\bfnu}(C,T) \ar[r]^{\boldsymbol{\rm RH}} & \Rep_{\bfnu}(C,T) .
 }\end{split}
\end{equation}
The horizontal maps are the Rieman--Hilbert correspondences associating to a connection its monodromy representation.
We remark that these are complex analytic isomorphisms, but they are transcendental.
Note that indeed, in the space $\Con^{\bfmu}_{\bfnu}(C,T)$, connections on decomposable bundles which are reducible but not diagonal are also identified to their diagonal reduction (these connections are said to be $s$-equivalent). Thus we've made equivalent identifications on each side and $\boldsymbol{\rm RH}$ gives an isomorphism of complex analytic spaces.

\section{A topological approach: the monodromy side} \label{sec:monodromy}

In this section we will discuss a few important facts about the map $\Phi^{{\rm top}} \colon \Rep_{\bar\nu}\big(\P^1,D\big)\to\Rep_{\bfnu}(C,T)$ defined in Section~\ref{subsec:Phi-rep}.
Through the Riemann--Hilbert correspondence (\ref{eq:RH}), we obtain analogous facts for the map $\Phi \colon \Con^{\bar\mu}_{\bar\nu}\big(\P^1,D\big) \to \Con^{\bfmu}_{\bfnu}(C,T)$.
The subsequent sections are devoted to further describing the map $\Phi$, and will not discuss monodromy representations any further.

The map $\Phi^{{\rm top}} \colon \Rep_{\bar\nu}\big(\P^1,D\big)\to\Rep_{\bfnu}(C,T)$ is studied in \cite[Section~4]{DiarraLoray2015}, and it is proved that it is dominant and generically a $2:1$ map.
Here, we want to provide a more precise result, namely:

\begin{Theorem}\label{thm:PhiMonodromySide}
Assume $2\nu\not\in\mathbb Z$. Then the map ${\Phi}^{{\rm top}} \colon \Rep_{\bar\nu}\big(\P^1,D\big)\to\Rep_{\bfnu}(C,T)$
is a branched covering of degree~$2$ $($in particular surjective$)$.
The branch locus consists of diagonal representations, forming the $2$-dimensional singular set of $\Rep_{\bfnu}(C,T)$.
The Galois involution of the covering
\[
 {\Psi}^{{\rm top}} \colon \ \Rep_{\bar\nu}\big(\P^1,D\big)\to\Rep_{\bar\nu}\big(\P^1,D\big)
\]
fixes the ramification locus, which consists of dihedral representations.
In particular we have an isomorphism $\Rep_{\bfnu}(C,T) \stackrel{\sim}{\longrightarrow} \Rep_{\bar\nu}\big(\P^1,D\big) / {\Psi}^{{\rm top}}$.
\end{Theorem}

\begin{Remark}\label{rem:PhiTopTilde} If we consider the analogous map
\[ \tilde{\Phi}^{{\rm top}} \colon \ \Rep_{\bar\nu}\big(\P^1,D\big)\to\widetilde{\Rep}_{\bfnu}(C,T)\]
at the level of isomorphism classes of representations, not onto the GIT quotient, then it is no longer surjective.
As we shall see in the proof of Theorem~\ref{thm:PhiMonodromySide} below, strictly triangular representation
(i.e., non diagonal ones) are precisely missing in the image of $\tilde{\Phi}^{{\rm top}}$; however,
this image defines a section of the GIT quotient $\Pi_{{\rm GIT}}$, so that ${\Phi}^{{\rm top}}$ is indeed surjective
on the GIT quotient.
\end{Remark}

\begin{proof}According to Remark~\ref{rem:PhiTopTilde}, we will consider the map $\tilde{\Phi}^{{\rm top}}$
between sets of isomorphism classes of representations.
It is proved in \cite[Section~4]{DiarraLoray2015} that the map $\tilde{\Phi}^{{\rm top}}$
is explicitly given by
\begin{gather*}
\tilde{\Phi}^{{\rm top}}\colon\ (M_0,M_1,M_t,M_\lambda,M_\infty) \longmapsto \begin{cases}
 A = M_1 M_t M_{\lambda},\\
 B = M_\lambda M_\infty,\\
 C_1 = M_{t},\\
 C_2 = M_\infty M_t M^{-1}_\infty.\end{cases}
 \end{gather*}
 and is generically two-to-one with Galois involution
\begin{equation}\label{eq:InvolCharacter}
{\Psi}^{{\rm top}}\colon \ (M_0,M_1,M_t,M_\lambda,M_\infty) \longmapsto (-M_0,-M_1,M_t,-M_\lambda,-M_\infty).
\end{equation}
More precisely, it is proved in \cite[Theorem~4.4]{DiarraLoray2015} that a representation $(A,B,C_1,C_2)$ is in the image of
$\tilde{\Phi}^{{\rm top}}$
provided that $A$ and $B$ generate an irreducible group: in that case, there exists a~matrix
$M\in\mathrm{SL}_2(\C)$, unique up to a sign, such that
\begin{equation}\label{eq:InvolutionCharacter}
M^{-1}AM=A^{-1},\qquad M^{-1}BM=B^{-1}\qquad \text{and}\qquad M^{-1}C_1M=C_2.
\end{equation}
Moreover, $M^2=-I$ and $(A,B,C_1,C_2)=\tilde{\Phi}^{{\rm top}}(M_0,M_1,M_t,M_\lambda,M_\infty)$ for
\begin{gather}
 M_0= -AM, \qquad
 M_1=ABC_2^{-1}M, \qquad
 M_t=C_1, \qquad
 M_\lambda=-BM, \qquad
 M_\infty=M .\label{eq:DescentCharacter}
 \end{gather}
Clearly, ${\Psi}^{{\rm top}}$ defines an involutive isomorphism of $\Rep_{\bar\nu}\big(\P^1,D\big)$
and we have $\tilde{\Phi}^{{\rm top}}\circ{\Psi}^{{\rm top}}=\tilde{\Phi}^{{\rm top}}$.
 Then to prove the statement of the theorem, it is enough to prove:
\begin{itemize}\itemsep=0pt
\item a representation $(A,B,C_1,C_2)$ is in the image of $\tilde{\Phi}^{{\rm top}}$ if, and only if,
it is either irreducible or diagonal;
therefore, the restriction $\Pi_{{\rm GIT}}\colon \text{image}\big(\tilde{\Phi}^{{\rm top}}\big) \to \Rep_{\bfnu}(C,T)$
will be one-to-one;
\item a diagonal representation has exactly one preimage by $\tilde{\Phi}^{{\rm top}}$ which is dihedral;
these are the fixed points of ${\Psi}^{{\rm top}}$;
\item an irreducible representation has two distinct preimages by $\tilde{\Phi}^{{\rm top}}$; equivalently,
a preimage cannot be fixed by ${\Psi}^{{\rm top}}$.
\end{itemize}
To complete the proof, we need to consider now the case where $A$ and $B$ generate a reducible, therefore triangular group.
A necessary condition to be in the image of $\tilde{\Phi}^{{\rm top}}$ is that the representation commutes (up to conjugacy) with
the action of the elliptic involution $\iota\colon C\to C$ which is as follows (see \cite[Section~4]{DiarraLoray2015}):
\[
 (A,B,C_1,C_2)\ \stackrel{\iota^*}{\mapsto}\ \big(A^{-1},B^{-1},C_2,C_1\big).
\]
We are going to prove that for all irreducible and diagonal representations, there exists an $M$ like above so that the representation is in the image of ${\Phi}^{{\rm top}}$; on the other hand, non-abelian triangular representations never commute with the involution $\iota^*$, and cannot be in the image.

Consider first the case where $A$ and $B$ (are triangular and) do not commute;
then we cannot conjugate simultaneously $A$ and $B$ to their inverse and the representation does not commute with~$\iota^*$.
Indeed, one of the two matrices $A$ and $B$ must have trace $\not=\pm2$, otherwise they commute;
if $\tr(A)\not=\pm2$, then the conjugacy with $A^{-1}$ must permute the two eigenvectors. One of them is common with $B$,
and if we can simultaneously conjugate $B$ to $B^{-1}$, then they share the two eigendirections, and therefore commute,
contradiction. So the representation $(A,B,C_1,C_2)$ is not in the image of ${\Phi}^{{\rm top}}$ in this case;
let us show that the representation is triangular in that case, i.e., $C_1$ and $C_2$ are also triangular.
Indeed, rewriting the relation
\[
 \underbrace{ABA^{-1}B^{-1}}_{[A,B]}=C_1\underbrace{(BA)C_2(BA)^{-1}}_{\tilde C_2}
\]
and, since traces satisfy $\tr\big(\tilde C_2\big)=\tr(C_1)$ and $\tr\big(C_1\tilde C_2\big)=2$, we deduce
(see for instance \cite[Section~6.2]{Loray2016}) that $C_1$ and $\tilde C_2$ generate also
a reducible group, i.e., there is a common eigenvector for $C_1$, $\tilde C_2$ and $[A,B]$,
which must coincide with the common eigenvector of $A$ and $B$.
We therefore conclude, when $A$ and $B$ do not commute, that the representation is triangular,
it does not commute with the involution~$\iota^*$, and is not in the image of ${\Phi}^{{\rm top}}$.

Assume now that $A$ and $B$ commute, and at least one of them has trace $\not=\pm 2$.
Therefore, we can assume that they are diagonal, and denote by $z_0,z_\infty\in\mathbb P^1$
the two common eigendirections. The relation writes
\[
 (BA)^{-1}C_1(BA)=(C_2)^{-1}.
\]
Denote by $z_i^+,z_i^-\in\mathbb P^1$ the two eigendirections of $C_i$, with $z_1^+$ and $z_2^+$
corresponding to the same eigenvalue. If $z_0$, $z_\infty$, $z_1^+$, $z_1^-$ are pairwise distinct,
then we have the following cross-ratio equivalences:
\[
 (z_0,z_\infty,z_1^+,z_1^-)\sim(z_0,z_\infty,z_2^-,z_2^+)\sim(z_\infty,z_0,z_2^+,z_2^-),
\]
where the left equivalence comes from the relation, and the second one is the cross-ratio's invariance under double transpositions. It follows that one can find a matrix $M\in\mathrm{SL}_2(\C)$, unique up to a sign, permuting $z_0$ and $z_\infty$, and sending $z_1^\pm$ to $z_2^\pm$, i.e., satisfying (\ref{eq:InvolutionCharacter}),
and therefore (\ref{eq:DescentCharacter}); the two choices provide the two preimages of the irreducible representation $(A,B,C_1,C_2)$
by ${\Phi}^{{\rm top}}$. If now $A$, $B$ and~$C_1$ share a common eigenvector, say $z_\infty=z_1^+$,
then the representation is triangular; the relation yields $z_\infty=z_1^+=z_2^-$ and a $M$
satisfying~(\ref{eq:InvolutionCharacter}) must permute
\[
 z_0\leftrightarrow z_\infty,\qquad z_1^+\leftrightarrow z_2^+,\qquad \text{and}\qquad z_1^-\leftrightarrow z_2^- .
\]
there exists such a $M$ if, and only if we also have $z_0=z_1^-=z_2^+$, i.e., the representation is diagonal.
We therefore conclude in this case that either the representation is irreducible and in the image of ${\Phi}^{{\rm top}}$,
or it is triangular and can be in the image of ${\Phi}^{{\rm top}}$ only if it is diagonal; we will check later that diagonal
representations are indeed in the image of ${\Phi}^{{\rm top}}$.

Assume now that $A$ and $B$ commute, and one of them is parabolic (i.e., $\tr=\pm2$ but not in the center).
Then we have only one common eigenvector $z_0=z_\infty$ for~$A$ and~$B$.
Again, if $z_\infty$, $z_1^+$, $z_1^-$ are pairwise distinct,
then one can find a matrix $M\in\mathrm{SL}_2(\C)$, unique up to a sign,
fixing~$z_\infty$, and sending $z_1^\pm$ to $z_2^\pm$, i.e., satisfying~(\ref{eq:InvolutionCharacter}).
Indeed:
\begin{itemize}\itemsep=0pt
\item If $AB=\pm I$, then $C_1=\pm C_2^{-1}$ and $(z_1^+,z_1^-)=(z_2^-,z_2^+)$ so that it suffices to choose
$M$ fixing $z_\infty$ and permuting $z_1^+$ and $z_1^-$: this characterize $M$ up to a sign, and $M$
must be a~projective involution, i.e., $M^2=\pm I$. Moreover, $M$ conjugates all parabolics fixing $z_\infty$
to their inverse, in particular~$A$ and~$B$.
\item If $AB\not=\pm I$, then we can choose~$M_0$ as before fixing $z_\infty$ and permuting $z_1^+$ and $z_1^-$
so that~$M_0$ is a projective involution, i.e., $M_0^2=-I$, conjugating~$A$,~$B$ and~$C_1$ to their inverse.
We deduce from the relation that $M^{-1}C_1M=C_2$ where $M:=M_0BA$.
Moreover, $A$,~$B$,~$M_0$ and therefore~$M$ fixes~$z_\infty$ and
\[ M^2=M_0BAM_0BA=\underbrace{M_0^2}_{-I}\underbrace{B^{-1}A^{-1}BA}_{I}=-I.\]
\end{itemize}
Again in these two cases, the representation $(A,B,C_1,C_2)$ is irreducible and has two preimages by ${\Phi}^{{\rm top}}$.
On the other hand, when $z_\infty$ is also fixed by $C_1$, say $z_\infty=z_1^+$, then the representation
is reducible non-abelian, the relation gives $z_0=z_1^-=z_2^+$, and we cannot find a $M$ fixing $z_\infty$
and sending $z_1^+$ to $z_2^+$ for instance, so we are not in the image of ${\Phi}^{{\rm top}}$.
Finally, if $A,B=\pm I$, then the relation shows that the representation is abelian (cyclic) and therefore diagonal;
it is in the image of ${\Phi}^{{\rm top}}$.

To finish the proof, note that the ramification locus of ${\Phi}^{{\rm top}}$ is given by the fixed points of the involution ${\Psi}^{{\rm top}}$, i.e., by those
$(M_0,M_1,M_t,M_\lambda,M_\infty)$ for which there exists $M\in\mathrm{SL}_2(\C)$ such that
\[
 M^{-1}M_tM=M_t,\qquad \text{and}\qquad M^{-1}M_iM=-M_i\qquad \text{for}\quad i=0,1,\lambda,\infty.
\]
This forces the representation to be dihedral, i.e., $(M_0,M_1,M_t,M_\lambda,M_\infty)$ of the form
\[ \left( \begin{pmatrix}0&a_0\\ -a_0^{-1}&0\end{pmatrix},\begin{pmatrix}0&a_1\\ -a_1^{-1}&0\end{pmatrix},\begin{pmatrix}a_t&0\\ 0&a_t^{-1}\end{pmatrix},\begin{pmatrix}0&a_\lambda\\ -a_\lambda^{-1}&0\end{pmatrix},\begin{pmatrix}0&a_\infty\\ -a_\infty^{-1}&0\end{pmatrix} \right)
\]
with $a_1a_\infty=a_0a_ta_\lambda$ and $M$ diagonal with eigenvalues $\pm\sqrt{-1}$. By diagonal conjugacy, one can normalize
$a_\infty=1$ and check that these representations are parametrized by $(a_0,a_\lambda)$ (we can fix $a_t={\rm e}^{2{\rm i}\pi\nu}$).
After lifting on $C$ we get the diagonal representation
\[ (A,B,C_1,C_2)=\left(\begin{pmatrix}-a_0&0\\ 0&-a_0^{-1}\end{pmatrix},\begin{pmatrix}-a_\lambda&0\\ 0&-a_\lambda^{-1}\end{pmatrix},\begin{pmatrix}a_t&0\\ 0&a_t^{-1}\end{pmatrix},\begin{pmatrix}a_t^{-1}&0\\ 0&a_t\end{pmatrix} \right)
\]
and we conclude that ${\Phi}^{{\rm top}}$ is one-to-one between above dihedral representations on~$\mathbb P^1$
and diagonal ones on~$C$.

To conclude, $\Rep_{\bfnu}(C,T)$ is precisely the quotient of the smooth variety
$\Rep_{\bar\nu}\big(\P^1,D\big)$ divided by the action of the involution ${\Psi}^{{\rm top}}$.
The locus of fixed points has codimension~$2$ and gives rise to a singular set in the quotient.
\end{proof}

\begin{Corollary}\label{cor:ConIsSingular}
The map constructed in Section~{\rm \ref{subsec:Phi-connections}},
\[
 {\Phi}\colon \ \Con^{\bar\mu}_{\bar\nu}\big(\P^1,D\big) \longrightarrow \Con^{\bfmu}_{\bfnu}(C,T) ,
\]
is a $2:1$ branched covering, branching over the locus of decomposable connections.
Moreover, this locus is the singular set for $\Con^{\bfmu}_{\bfnu}(C,T)$.
\end{Corollary}

\begin{proof}The Riemann--Hilbert correspondence is one-to-one between isomorphism classes of connections and conjugacy classes of representations.
On $\big(\P^1,D\big)$, all connections (representations) are irreducible and the GIT quotient therefore coincide with isomorphism (conjugacy) classes.
The Riemann--Hilbert map is an analytic isomorphism in this case.
On $(C,T)$, all connections (representations) are semistable,
as was shown in Lemma~\ref{lemma:definition-Phi}.
Again the Riemann--Hilbert map induces an homeomorphism between the GIT quotients which is analytic in restriction to the stable part.
Clearly, the strictly semistable locus of $\Con^{\bfmu}_{\bfnu}(C,T)$, which has codimension $2$,
must be singular, since, locally around a point in the strictly semistable locus, the complement of this locus is not simply connected.
\end{proof}

We now proceed to discuss the symplectic structure of the spaces involved.
There is a natural holomorphic symplectic structure on the spaces of representations arising from works of Atyiah-Bott and Goldman.
In the compact case,
it is defined as follows. Let $S$ be a compact surface (without boundary), and denote by $\Rep(S)$ the space
of representations $\rho\colon \pi_1(S)\to \mathrm{SL}_2(\mathbb C)$, up to conjugacy.
The tangent space at an irreducible representation $\rho\in\Rep(S)$ identifies with the cohomology group
$H^1(\pi_1(S),\mathfrak{sl}_2(\mathbb C)_{\mathrm{Ad}_\rho})$ where $\mathfrak{sl}_2(\mathbb C)_{\mathrm{Ad}_\rho}$
is the Lie algebra $\mathfrak{sl}_2(\mathbb C)$, viewed as a~$\pi_1(S)$-module via the adjoint action of $\rho$
(see \cite[Section~1.2]{Goldman1984}). Equivalently, one can define the tangent space via the de Rham cohomology
by ${H}^1_{\rm dR}(S,E_\rho)$, where $E_\rho$ is the flat $\mathfrak{sl}_2(\mathbb C)$-bundle whose monodromy
is given by the adjoint representation $\mathrm{Ad}_\rho$ on the Lie algebra $\mathfrak{sl}_2(\mathbb C)$.

Now, there is a natural bilinear map on the tangent space which can be defined, in the de~Rham setting, by
\begin{equation} \label{eq:dR-pairing1}
 H^1_{\rm dR}(S,E_\rho)\otimes H^1_{\rm dR}(S,E_\rho) \overset{\text{cup product}}{\longrightarrow} H^2_{\rm dR}(S,E_\rho\otimes E_\rho)
 \overset{\text{Killing form}}{\longrightarrow} H^2_{\rm dR}(S,\mathbb C) .
\end{equation}
Combining with the canonical isomorphism
\begin{equation} \label{eq:dR-pairing2}
 H^2_{\rm dR}(S,\mathbb C) \longrightarrow \C, \qquad\text{given by}\quad \eta \longmapsto \int_S \eta ,
\end{equation}
we obtain the pairing
\begin{align*}
 H^1_{\rm dR}(S,E_\rho)\otimes H^1_{\rm dR}(S,E_\rho) & \longrightarrow \mathbb C, \\
 (\alpha,\beta) & \longmapsto \int_S\mathrm{trace}(\alpha\wedge\beta) ,
\end{align*}
see \cite[Section~1.8]{Goldman1984} for more details.
In this way, we get a holomorphic $2$-form $\omega$ on $\Rep(S)$, which turns out to be closed and non-degenerate,
i.e.,
\[
 {\rm d}\omega=0\qquad \text{and}\qquad \underbrace{\omega\wedge\cdots\wedge\omega}_{\text{half the dimension of }\Rep(S)}\not=0 ,
\]
(cf.~main theorem in \cite[Section~1.7]{Goldman1984}).
The non-compact case, where $S$ has a non-empty boundary~$\partial S$, was considered by Iwasaki in \cite[Section~3]{Iwasaki1992}
(see also \cite[Section~5]{Iwasaki2002} for more details, including explicit computations for the punctured sphere case).
There we have to deal with relative cohomology, i.e., cohomology of pairs $(S,\partial S)$ in order to take into account
the fixed conjugacy classes at punctures.
We have a map $\partial\colon \Rep(S)\to\Rep(\partial S)$ associating to a representation~$\rho$,
the conjugacy classes at punctures; our space of representations correspond to the fibers of this map
and we have to restrict $H^1_{\rm dR}(S,E_\rho)$ (closed $1$-forms with compact support) to some subspace.
But at the end, the definition in the de~Rham setting is the same and Iwasaki proves that this is symplectic
in restriction to $\partial$-fibers.
Going back to our map $\Phi^{{\rm top}}$, let us simply denote by $\omega_{\P^1}^{\rm IG}$ and $\omega_{C}^{\rm IG}$
the Goldman--Iwasaki symplectic $2$-form on the corresponding side.

\begin{Proposition}\label{prop:PhiSymplecticMonodromy}
The map ${\Phi}^{{\rm top}}$ of Theorem~{\rm \ref{thm:PhiMonodromySide}} is symplectic, up to a constant factor of~$2$. Namely, it satisfies
\[ \big({\Phi}^{{\rm top}}\big)^*\omega_{C}^{\rm IG}=2 \omega_{\P^1}^{\rm IG}.\]
\end{Proposition}

\begin{proof}This fact is true in general whenever we have a covering map $\pi\colon X\to Y$ between Riemann surfaces.
In the construction of the Goldman--Iwasaki bracket (\ref{eq:dR-pairing1}), everything commutes with the base change.
In the last step (\ref{eq:dR-pairing2}), where we integrate on $Y$ instead of $X$, the result is multiplied by the degree $\deg(\pi \colon X\to Y)$.
The ramified case is analogous using relative cohomology.
In our case we have $\pi \colon C\to\P^1$ with $\deg(\pi) = 2$.
\end{proof}

From the main theorem of \cite{Iwasaki1992}, the Riemann--Hilbert map $\operatorname{RH} \colon \Con^{\bar\mu}_{\bar\nu}(X, D) \to \Rep_{\bar\nu}(X, D)$ is symplectic, i.e., pulls back the Goldman--Iwasaki bracket $\omega^{\rm IG}$ to the Atiyah--Bott symplectic structure $\omega^{AB}$ on the moduli space of connections.
In particular, our pullback map $\Phi$ is also symplectic (up to a multiplicative factor of~$2$)
and we will be able to compute the symplectic structure on $\Con^{\bfmu}_{\bfnu}(C,T)$ from the known symplectic structure on $\Con^{\bar\mu}_{\bar\nu}\big(\P^1,D\big)$ (for the latter, see~\cite{LoraySaito2015} and references therein).

\begin{Corollary}\label{cor:PhiIsSymplectic}
The map constructed in Section~{\rm \ref{subsec:Phi-connections}}
\[ {\Phi}\colon \ \Con^{\bar\mu}_{\bar\nu}\big(\P^1,D\big) \longrightarrow \Con^{\bfmu}_{\bfnu}(C,T)\]
is symplectic up to a factor of $2$:
\[
 {\Phi}^*\omega_{C}^{AB} = 2\, \omega_{\P^1}^{AB} .
\]
\end{Corollary}

\section{Genericity assumptions} \label{sec:genericity}

In this section we will briefly describe the geometry of the moduli spaces of parabolic bundles we work with. We will explain which families of bundles are particularly special, and define a~\emph{generic bundle} to be one not belonging to these families.
For our own convenience, we will first define these special families in terms of coordinate descriptions of the moduli spaces of parabolic bundles. Only later we describe such families intrinsically.

We begin with bundles over $\P^1$.
Let $\lambda, t \in \P^1 \setminus \{0,1,\infty\}$ be different points, and let $D$ be the divisor $D=0+1+\lambda+\infty+t$. We are interested in $\bar\mu$-semistable parabolic bundles $(E, \bar\ell)$ of degree zero over the marked curve $\big(\P^1,D\big)$.

According to \cite[Proposition~6.1]{FernandezVargas2016} (which is a direct consequence of the results in \cite[Section~6.1]{LoraySaito2015}), the moduli space $\Bun^{\bar\mu}\big(\P^1,D\big)$, for the weights fixed in Definition~\ref{def:mu-nu}, is isomorphic to a Del~Pezzo surface of degree~4 which we denote~$\S$.
We refer the reader to \cite[Section~3]{LoraySaito2015} for the explicit construction of moduli spaces of rank~2 parabolic bundles over $\P^1$ as projective varieties.
An explicit bridge between~$\mathcal{S}$ and parabolic bundles is given, in the generic case, by Remark~\ref{rmk:coordinates-Bun-P1P1}.

The surface $\S$ is a smooth projective surface that is obtained by blowing-up 5 particular points $D_i \in\P^2$. It is well-known that this surface $\S$ has exactly 16 rational curves of self-intersection~$-1$. Namely, the five exceptional divisors from the blow-up $E_i$, the strict transform of the conic~$\Pi$ passing through the five points, and the strict transform of the 10 lines $L_{i j}$ passing through every possible pair $(D_i,D_j)$.
The five points~$D_i$ in the smooth conic $\Pi$ are in the same position as the five points $0,1,\lambda,t,\infty \in \P^1$.

Let us take the above geometric description of $\Bun^{\bar\mu}\big(\P^1,D\big)$ as the base for our definition of genericity. An intrinsic interpretation will follow.

\begin{Definition}\label{def:generic-bunle-P1}
 We will say that a parabolic bundle is \textit{generic} in $\Bun^{\bar\mu}\big(\P^1,D\big) \cong \S$ if it lies outside the union of the 16 $(-1)$-curves $\{\Pi, E_i, L_{i j}\}$.
 A parabolic connection in $\Con^{\bar\mu}_{\bar\nu}\big(\P^1,D\big)$ will be called \textit{generic} if the underlying parabolic bundle is $\bar\mu$-semistable and generic.
 We denote by $\Bun^{\bar\mu}\big(\P^1,D\big)^0$ and $\Con^{\bar\mu}_{\bar\nu}\big(\P^1,D\big)^0$ the open subsets of generic bundles and connections, respectively.
\end{Definition}

\begin{Proposition}\label{prop:intrinsic-def-generic}
 A parabolic bundle $(E,\bar\ell)$ of degree zero belongs to $\Bun^{\bar\mu}\big(\P^1,D\big)^0$ if and only if all the following conditions are true:
 \begin{itemize}\itemsep=0pt
 \item $E = \mathcal{O}_{\P^1} \oplus \mathcal{O}_{\P^1}$.
 \item There exists no subbundle $L=\mathcal{O}_{\P^1} \hookrightarrow E$, such that $\ell_i, \ell_j \subset L$ for any pair $i \neq j$.
 \item There exists no subbundle $L=\mathcal{O}_{\P^1}(-1) \hookrightarrow E$, such that $\ell_i, \ell_j, \ell_k, \ell_m \subset L$ for pairwise different $i$, $j$, $k$, $m$.
 \end{itemize}
\end{Proposition}

The proof of this proposition follows immediately from \cite[Table~1, Section~6.1]{LoraySaito2015}, where the non-generic families $\Pi$, $E_i$, $L_{i j}$ are explicitly described. In the work just cited these families are described for parabolic bundles of degree $-1$, but translation to bundles of degree zero is done in a straightforward manner by performing one positive elementary transformation.

We know from \cite{LoraySaito2015} that the \emph{coarse moduli space} of indecomposable parabolic bundles is a~non-separated variety obtained by gluing together a finite number of spaces $\Bun^{\bar\mu}\big(\P^1,D\big)$ for suitable choices of weight vectors $\bar\mu$. As the weights vary, the bundles in the special families $\{\Pi, E_i, L_{i j}\}$ may become unstable, and new bundles that were previously unstable are now semistable. However, the bundles represented in $\S \setminus \{\Pi, E_i, L_{i j}\}$ are always stable and thus common to every chart.
In particular, this means that the set $\Bun^{\bar\mu}\big(\P^1,D\big)^0$ coincides with the set of \textit{closed points} in the non-separated coarse moduli space of indecomposable bundles.

Let us now move on to parabolic bundles over $(C,T)$.
As shown in \cite[Theorem~A]{FernandezVargas2016}, the moduli spa\-ce $\Bun^{\bfmu}(C,T)$ is isomorphic to $\P^1 \times \P^1$.
Indeed, an explicit isomorphism $\Bun^{\bfmu}(C,T) \rightarrow$ \mbox{$\P^1_z \times \P^1_w$} is constructed in \cite[Section~4.3]{FernandezVargas2016}.
We briefly recall the construction here.
The moduli space of semistable rank~2 vector bundles with trivial determinant, $\Bun(C)$ is canonically isomorphic to the quotient of $\operatorname{Jac}(C)$ by the elliptic involution~\cite{Tu1993}.
The map $\Bun^{\bfmu}(C,T) \rightarrow \P^1_z$, defining the first coordinate, is given by forgetting the parabolic structure on $(E,\bar\ell)$, and mapping $E$ to $\operatorname{Jac}(C) / \pi \cong \P^1$ via the above isomorphism. Let us call this map $\mathrm{Tu}$.
The second coordinate is obtained by precomposing~$\mathrm{Tu}$ with the following automorphism $\varphi$ of $\Bun^{\bfmu}(C,T)$: Given $(E,\bar\ell)$, perform a positive elementary transformation above each of the two poles $\bft_i\in T$, and a twist by $\mathcal{O}_C(-\bfw_\infty)$. The resulting parabolic bundle $\varphi(E,\bar\ell)$ has trivial determinant and thus belongs to $\Bun^{\bfmu}(C,T)$. This map is in fact an involution, and in coordinates~$(z,w)$ it is, by definition, precisely the map $\varphi \colon (z,w) \mapsto (w,z)$.

The map $\phi \colon \Bun^{\bar\mu}\big(\P^1,D\big) \to \Bun^{\bfmu}(C,T)$ has been described using the above coordinate system.
It is proved in \cite[Section~6.3]{FernandezVargas2016} that it transforms the special $(-1)$-curves of~$\S$ to either horizontal or vertical lines defined by $z=0,1,\lambda,\infty$, and $w=0,1,\lambda,\infty$.
Our definition of a~generic bundle in $\Bun^{\bfmu}(C,T)$ will exclude these special lines, together with the curve described below.

\begin{Definition}\label{def:Sigma}
 We denote by $\boldsymbol{\Sigma} \subset \Bun^{\bfmu}(C,T)$ the set of strictly $\bar\mu$-semistable parabolic bundles.
 We also denote $\Sigma = \phi^{-1}(\boldsymbol{\Sigma}) \subset \Bun^{\bar\mu}\big(\P^1,D\big)$.
\end{Definition}

\begin{Definition}\label{def:generic-bunle-C}
 We will say that a parabolic bundle is \textit{generic} in $\Bun^{\bfmu}(C,T) \cong \P^1_z \times \P^1_w$ if it lies outside the following loci:
 \begin{itemize}\itemsep=0pt
 \item The union of the 8 lines $z=0,1,\lambda,\infty$, and $w=0,1,\lambda,\infty$,
 \item The strictly $\bar\mu$-semistable locus $\boldsymbol{\Sigma}$,
 \end{itemize}
 A parabolic connection in $\Con^{\bfmu}_{\bfnu}(C,T)$ will be called \textit{generic} if the underlying parabolic bundle is $\bar\mu$-semistable and generic.
 We denote by $\Bun^{\bfmu}(C,T)^0$ and $\Con^{\bfmu}_{\bfnu}(C,T)^0$ the open subsets of generic bundles and connections, respectively.
\end{Definition}

It is natural to exclude the curve $\boldsymbol{\Sigma}$ from our definition of genericity, since such bundles are not $\bar\mu$-stable.
In fact, because of the following theorem, this curve plays a crucial role in the description of the moduli space of parabolic connections.

\begin{Theorem}[Theorem B and Theorem~6.4 in~\cite{FernandezVargas2016}]\label{thm:Nestor6.4}
 The map $ \phi \colon \Bun^{\bar\mu}\big(\P^1,D\big) \to \Bun^{\bfmu}(C,T)$ is a $2:1$ ramified cover. The branch locus is precisely the strictly semistable locus $\boldsymbol{\Sigma}$, and the ramification locus coincides with $\Sigma$. Moreover, both curves are isomorphic to the elliptic curve~$C$.
\end{Theorem}

The vertical and horizontal lines $z=0,1,\lambda,\infty$, and $w=0,1,\lambda,\infty$, in $\P^1_z \times \P^1_w$ are excluded mostly for technical reasons, since we are interested in the map $\phi \colon \Bun^{\bar\mu}\big(\P^1,D\big) \to \Bun^{\bfmu}(C,T)$, and these lines are the images of the families of non-generic bundles over $\P^1$. However, these families do represent parabolic bundles with very specific properties. For example, if a parabolic bundle $(\mathbf{E}, \bar\ell)$ satisfies $z \in \{0,1,\lambda,\infty\}$, then either
 $\mathbf{E} = \mathbf{L} \oplus \mathbf{L}$, with $\mathbf{L}$ a torsion bundle, or
 $\mathbf{E} = \mathbf E_0 \otimes \mathbf L$, where $\mathbf E_0$ is the unique non-trivial extension of $\mathcal{O}_C$ by $\mathcal{O}_C$, and $\mathbf L$ is a torsion bundle.
See \cite{Tu1993} and \cite[Section~4.3]{FernandezVargas2016}.

\begin{figure}[t] \centering
 \includegraphics[width=0.9\linewidth]{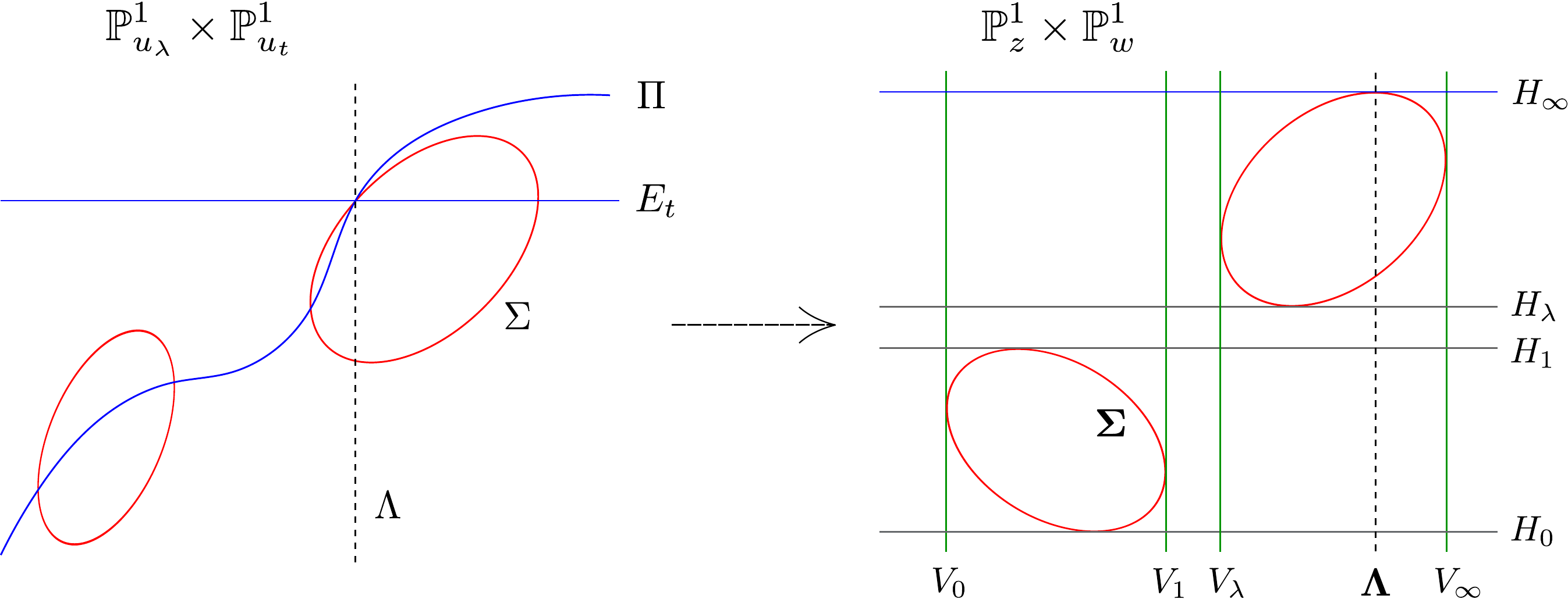}
 \caption{Non-generic curves on $\P^1_z \times \P^1_w$ and their counterparts in $\P^1_{u_\lambda} \times \P^1_{u_t}$.} \label{fig:non-generic}
\end{figure}

Fig.~\ref{fig:non-generic} shows some of the curves in $\Bun^{\bar\mu}\big(\P^1,D\big)$ which are excluded from the definiton of a~generic bundle, and how these curves are transformed under $\phi$.

Later on, once we perform computations in coordinates, we will sometimes need to exclude the vertical line $\boldsymbol{\Lambda} = \{ z=t \}$, which appears as a polar divisor in some of our formulas. The pre-image $\phi^{-1}(\boldsymbol{\Lambda})$ corresponds to another vertical line $\Lambda \subset \P^1_{u_\lambda} \times \P^1_{u_t}$.

\section{Recap of previously known results} \label{sec:recap}

In this section we will further recall several facts from \cite{FernandezVargas2016, LoraySaito2015} in order to make our results precise and to put them into context. We restrict ourselves to the cases that are relevant to us. We refer the reader to the original papers cited for a detailed treatment and for more general cases.
Definitions and basic results about parabolic bundles and connections have already been presented in Section~\ref{sec:moduli-spaces}.

\subsection{Moduli spaces of parabolic bundles}

In the previous section we have mentioned that the moduli space $\Bun^{\bar\mu}\big(\P^1,D\big)$ is isomorphic to a Del~Pezzo surface, which we continue to denote~$\S$. Below we present two coordinate systems that can be used to describe the set $\Bun^{\bar\mu}\big(\P^1,D\big)^0$ of generic parabolic bundles. For later convenience we present them as two remarks.

\begin{Remark}\label{rmk:coordinates-Bun-P1P1}
 The first coordinate system relies on Proposition~\ref{prop:intrinsic-def-generic}. According to this proposition, a generic parabolic bundle has a trivial underlying bundle, and no two parabolics lie on a~same trivial subline bundle.
 Consider such a generic bundle and let us introduce an affine coordinate~$\zeta$ on the projectivized fibers of the (trivial) bundle. After a fractional linear transformation, we may assume that the parabolic structures over the points $0$, $1$, $\infty$ are given by $\zeta=0,1,\infty$, respectively. Under this situation, any parabolic bundle is completely determined by two parameters $u_\lambda, u_t \in\P^1$, which represent the value of the $\zeta$-coordinate for the parabolic structures over $\lambda$ and $t$, respectively. This assignment defines a birational map $\Bun^{\bar\mu}\big(\P^1,D\big) \dashrightarrow \P^1_{u_\lambda} \times \P^1_{u_t}$, which provides a coordinate system.
 The surface $\mathcal{S}$ is recovered by blowing-up the points $(0,0)$, $(1,1)$, $(\lambda,t)$ and $(\infty,\infty)$ in $\P^1_{u_\lambda} \times \P^1_{u_t}$.
\end{Remark}

\begin{Remark}\label{rmk:coordinates-Bun-P2}
 The second system is based on the fact that the Del~Pezzo surface $\S$ is defined as the blow-up of $\P^2$ at five points, so there is a canonical birational map $\P^2\dashrightarrow\S$. Since $\Bun^{\bar\mu}\big(\P^1,D\big)^0$ excludes (together with other curves) the exceptional divisors of the blow-up, this map defines a one-to-one map between an open subset of $\P^2$ and $\Bun^{\bar\mu}\big(\P^1,D\big)^0$. Fixing homogeneous coordinates $[b_0:b_1:b_2]$, the space $\P^2_b$ defines a coordinate system for $\Bun^{\bar\mu}\big(\P^1,D\big)^0$.
 An explicit interpretation of these coordinates can be found (for bundles of degree $-1$) in \cite[Section~3.6]{LoraySaito2015}.
\end{Remark}

Now, let $C\subset\P^2$ be an elliptic curve such that in some fixed affine chart it is given by the equation
\[
 y^2 = x(x-1)(x-\lambda), \qquad \lambda\in\C\setminus \{0,1\}.
\]
As described in Section~\ref{subsec:Phi-bundles}, a parabolic bundle $(E,\bar\ell)$ on $\P^1$ can be lifted to the elliptic curve using the elliptic covering $\pi\colon C \to \P^1$.
After a series of birational transformations, we obtain a parabolic bundle on $C$ with parabolic structure supported over the divisor $T=\pi^*(t)$. This defines a map
\[
 \phi \colon \ \Bun^{\bar\mu}\big(\P^1,D\big) \longrightarrow \Bun^{\bfmu}(C,T)
\]
between moduli spaces.
The map $\phi$ is a $2:1$ ramified covering (as was already stated in~Theorem~\ref{thm:Nestor6.4}).
The domain space is the Del~Pezzo surface $\S$ discussed above, and the target space is proved to be isomorphic to $\P^1\times\P^1$ in \cite[Theorem~A]{FernandezVargas2016} (see Section~\ref{sec:genericity} for a brief description of this coordinate system).

\begin{Definition}\label{def:psi}
 We define $\psi \in \operatorname{Aut}(\S)$ as the involution of $\S\cong \Bun^{\bar\mu}\big(\P^1,D\big) $ which permutes the two sheets of the map $\phi$ and fixes every point in the ramification divisor $\Sigma$.
\end{Definition}

The above involution is a lift of a de~Jonqui\`{e}res automorphism of~$\P^2_b$ (a birational automorphism of degree 3 that preserves a pencil of lines through a point, and a pencil of conics through four other points).
The curve $\Sigma\subset\S$ is, according to Theorem~\ref{thm:Nestor6.4}, the ramification locus of $\phi$. Therefore, it corresponds to the curve of fixed points of~$\psi$. This involution is further discussed in \cite[Section~6.4]{FernandezVargas2016}.

\begin{Remark}\label{rmk:quotient-psi}
 The moduli space $\Bun^{\bar\mu}\big(\P^1,D\big)$ is endowed with an involution $\psi$, in such a~way that the quotient of $\Bun^{\bar\mu}\big(\P^1,D\big) / \psi$ is isomorphic to $\Bun^{\bar\mu}(C,T)$.
 We have a~similar situation in Theorem~\ref{thm:PhiMonodromySide}, where $\Rep_{\bar\nu}\big(\P^1,D\big) / \Psi^{{\rm top}}$ is identified to $\Rep_{\bfnu}(C,T)$.
 Our description of $\Con^{\bfmu}_{\bfnu}(C,T)$ will be analogous, and so it is clear the Galois involution~$\psi$ (and its counterpart~$\Psi$) will play a crucial role in the present work.
\end{Remark}

The involution $\psi$, defined above as an automorphism of the moduli space $\Bun^{\bar\mu}\big(\P^1,D\big)$, can also be defined by an explicit geometric transformation acting directly on bundles (i.e., not just on equivalence classes of bundles). Because of the importance of $\psi$ in the sequel, we describe such geometric construction below.

\begin{Remark}\label{rmk:definition-psi}
 Consider a parabolic bundle $(E, \bar\ell)$ on $\big(\P^1,D\big)$ of degree zero. We recall that the divisor $D$ is given by $D=0+1+\lambda+\infty+t$. Consider the following transformations.
 \begin{enumerate}\itemsep=0pt
 \item \label{item:psi1} Perform a positive elementary transformation (cf.~Definition~\ref{def:elem-transf}) at each of the parabolic points in the divisor $W=0+1+\lambda+\infty$. This defines a new bundle of degree~4.
 \item \label{item:psi2} Tensor the previous bundle with the line bundle $\mathcal{O}_{\P^1}(-2)$. This bundle now has degree zero.
 \end{enumerate}
 The bundle obtained by the above steps is again a parabolic bundle on $\big(\P^1,D\big)$ of degree zero. It represents a class in the moduli space $\Bun^{\bar\mu}\big(\P^1,D\big)$, and so this transformation defines a self-map of the moduli space.
 It was proved in \cite[Proposition~6.5]{FernandezVargas2016} that such map is precisely the involution $\psi$.
 In other words, the transformation defined by steps (\ref*{item:psi1}) and (\ref*{item:psi2}) above coincides with the Galois involution of the covering $\phi \colon \Bun^{\bar\mu}\big(\P^1,D\big) \to \Bun^{\bfmu}(C,T)$.
\end{Remark}

\begin{Remark}\label{rmk:definition-Psi}
 The above construction naturally extends to parabolic connections if we twist by the rank~1 connection $\big(\mathcal{O}_{\P^1}(-2), \eta\big)$ in the last step, where~$\eta$ is the unique rank~1 connection having poles at $W$ with residue $\frac{1}{2}$ at each pole.
 We remark that the stability condition is always preserved. Indeed, if a parabolic point $p$ has weight $\mu$, the stability condition is preserved as long as the new weight of~$p$ in $\operatorname{elem}^+_p(E, \bar\ell)$ is defined to be $1-\mu$.
 Since in step (\ref*{item:psi1}) we perform elementary transformations only at the points in~$W$, whose weights have been fixed to be $\mu=\frac{1}{2}$, we conclude that the stability is indeed preserved for \textit{the same} weight vector $\bar\mu$.
\end{Remark}

We denote this transformation by $\Psi$. As the lemma below shows, this transformation corresponds to the Galois involution of the ramified cover $\Phi$ between moduli spaces of connections.

\begin{Lemma} \label{lemma:psi-Galois}
 The transformation $\Psi$ defined by steps {\rm (\ref*{item:psi1})} and {\rm (\ref*{item:psi2})} above, twisting by the rank~1 connection $(\mathcal{O}_{\P^1}(-2), \eta)$ in the second step, coincides with the Galois involution of the covering $\Phi \colon \Con^{\bar\mu}_{\bar\nu}\big(\P^1,D\big) \to \Con^{\bfmu}_{\bfnu}(C,T)$.
\end{Lemma}

\begin{proof} This is a direct consequence of the analogous result for the map $\Phi^{{\rm top}}$ between the spaces of monodromy representations discussed in Section~\ref{sec:monodromy}.
 The Galois involution of $\Phi^{{\rm top}}$ is given in~(\ref{eq:InvolCharacter}). It corresponds to the transformation that preserves $M_t$ and transforms $M_i \mapsto -M_i$, for $i=0,1,\lambda,\infty$.
 In the definition of $\psi$ above, step (\ref*{item:psi1}) does not change the monodromy representation, while step (\ref*{item:psi2}) changes the local monodromy precisely by a factor of $-1$ at the points $0$, $1$, $\lambda$, $\infty$ (cf.~Remark~\ref{rmk:rep-con-equivalence}).
\end{proof}

\subsection[Moduli spaces of connections over P1]{Moduli spaces of connections over $\boldsymbol{\P^1}$}\label{sec:connections-P1}

Recall that the space $\Con^{\bar\mu}_{\bar\nu}\big(\P^1,D\big)$ carries a natural symplectic structure in such a way that the map $\operatorname{Bun} \colon \Con^{\bar\mu}_{\bar\nu}\big(\P^1,D\big) \to \Bun^{\bar\mu}\big(\P^1,D\big)$ is a Lagrangian fibration. In \cite{LoraySaito2015} it is shown that the so-called \textit{apparent map} defines a dual Lagrangian fibration. Given a connection $\nabla$ on a~bundle $E$ and a rank 1 subbundle $L\subset E$, the apparent map is defined by the zero divisor of the composite map
\[\begin{tikzcd}[column sep=small]
 L \arrow[r, hook]
 &
 E \arrow[r, "\nabla"]
 &
 E\otimes\Omega^1_{\P^1}(D) \arrow[r]
 &
 (E/L)\otimes\Omega^1_{\P^1}(D).
\end{tikzcd}\]
Note that the apparent map is defined geometrically as the set of points of tangency between the Riccati foliation defined by $\nabla$ on $\P(E)$ and the section induced by the subbundle $L$.

For a generic connection of degree $-1$, the underlying bundle is $E=\mathcal{O}_{\P^1} \oplus \mathcal{O}_{\P^1}(-1)$. This bundle has a unique trivial subbundle $L=\mathcal{O}_{\P^1}$, which provides a canonical choice for the apparent map. In this case we obtain a rational map
\begin{equation}\label{eq:App-P1}
 \operatorname{App} \colon \Con^{\bar\mu}_{\bar\nu}\big(\P^1,D\big) \dashrightarrow |\mathcal{O}_{\P^1}(n-3)| \cong \P^{n-3},
\end{equation}
where $n$ denotes the number of singularities (in our particular case $n=5$).
For generic connections of degree zero the underlying bundle is $\mathcal{O}_{\P^1} \oplus \mathcal{O}_{\P^1}$, and we may perform an elementary transformation to replace it by $\mathcal{O}_{\P^1} \oplus \mathcal{O}_{\P^1}(-1)$. After this, we may proceed as above. This extends the definition of the apparent map to bundles of degree zero.

The Lagrangian fibrations provide a description of the geometric structure of the space of connections. Indeed, the map
\begin{equation}\label{eq:AppBun-P1}
 \operatorname{App} \times \operatorname{Bun} \colon \ \Con^{\bar\mu}_{\bar\nu}\big(\P^1,D\big) \dashrightarrow \P^2 \times \P^2 ,
\end{equation}
when restricted to the space of generic connections $\Con^{\bar\mu}_{\bar\nu}\big(\P^1,D\big)^0$, and under a simple assumption on the residual eigenvalues, is an open embedding \cite[Theorem~4.2]{LoraySaito2015}. Moreover, a suitable compactification of the space of generic bundles makes the above map an isomorphism.

\subsection{A universal family of connections} \label{subsec:universal-family-P1}

Another result that we try to imitate is the construction of an explicit \textit{universal family} for $\Con^{\bar\mu}_{\bar\nu}\big(\P^1,D\big)$.
By universal family we mean an algebraic family $\mathcal{U} = \{(E_\theta, \nabla_\theta, \ell_\theta)\}$ with the property that the natural map $\mathcal{U} \rightarrow \Con^{\bar\mu}_{\bar\nu}\big(\P^1,D\big)$, which assigns to each parabolic connection its isomorphism class, is dominant and injective. Thus every generic element of the moduli space is represented by a unique connection $(E_\theta, \nabla_\theta, \ell_\theta) \in \mathcal{U}$.
We describe such a family below.

As pointed out in Remark~\ref{rmk:coordinates-Bun-P1P1}, a generic parabolic bundle $(E,\bar\ell)$ of degree zero and polar divisor $D=0+1+\lambda+\infty+t$ has $E=\mathcal{O}_{\P^1} \oplus \mathcal{O}_{\P^1}$ as underlying bundle, and we may assume the parabolic structure is given by $\bar\ell = (0,1,u_\lambda,\infty,u_t)$, for some $u_\lambda, u_t\in\P^1$.
This defines a~coordinate system on $\Bun^{\bar\mu}\big(\P^1,D\big)$, which we continue to use below.

For each pair $(u_\lambda, u_t)$ with $u_\lambda,u_t \neq \infty$, we define a parabolic connection $\nabla_0 (u_\lambda, u_t)$ and two Higgs bundles $\Theta_i (u_\lambda, u_t)$, $i=1,2$, compatible with the parabolic structure $(u_\lambda, u_t)$ by the following explicit formulas:
\begin{gather*}
 \nabla_0 (u_\lambda, u_t) = \frac{1}{4}
 \begin{pmatrix}
 -1 & 0 \\
 -2 - 4\nu & 1
 \end{pmatrix}
 \frac{{\rm d}x}{x} + \frac{1}{4}
 \begin{pmatrix}
 1 + 4\nu & -4\nu \\
 2 + 4\nu & -1 - 4\nu
 \end{pmatrix}
 \frac{{\rm d}x}{x-1} \\
\hphantom{\nabla_0 (u_\lambda, u_t) =}{} + \frac{1}{4}
 \begin{pmatrix}
 -1 & 2u_\lambda \\
 0 & 1
 \end{pmatrix}
 \frac{{\rm d}x}{x-\lambda} + \nu
 \begin{pmatrix}
 -1 & 2 u_t \\
 0 & 1 \\
 \end{pmatrix}
 \frac{{\rm d}x}{x-t} ,
\\
 \Theta_1 (u_\lambda, u_t) =
 \begin{pmatrix}
 0 & 0 \\
 1-u_t & 0 \\
 \end{pmatrix}
 \frac{{\rm d}x}{x} +
 \begin{pmatrix}
 u_t & -u_t \\
 u_t & -u_t \\
 \end{pmatrix}
 \frac{{\rm d}x}{x-1} +
 \begin{pmatrix}
 -u_t & u_t^2 \\
 -1 & u_t \\
 \end{pmatrix}
 \frac{{\rm d}x}{x-t} ,
\\
 \Theta_2 (u_\lambda, u_t) =
 \begin{pmatrix}
 0 & 0 \\
 1-u_\lambda & 0 \\
 \end{pmatrix}
 \frac{{\rm d}x}{x} +
 \begin{pmatrix}
 u_\lambda & -u_\lambda \\
 u_\lambda & -u_\lambda \\
 \end{pmatrix}
 \frac{{\rm d}x}{x-1} +
 \begin{pmatrix}
 -u_\lambda & u_\lambda^2 \\
 -1 & u_\lambda \\
 \end{pmatrix}
 \frac{{\rm d}x}{x-\lambda} .
\end{gather*}

\begin{Remark} \label{rmk:nabla0-Lagrangian}
 There are many possible choices for the rational section $\nabla_0 \colon \Bun^{\bar\mu}\big(\P^1,D\big) \dashrightarrow \Con^{\bar\mu}_{\bar\nu}\big(\P^1,D\big)$. The important property of this one is that it is Lagragian with respect to the natural symplectic form.
 This will be important later for the explicit computation of the symplectic structure on $\Con^{\bfmu}_{\bfnu}(C,T)$.
 The apparent map $\operatorname{App}$ is a Lagrangian fibration, which is proved in \cite{LoraySaito2015}
 to be transversal to the bundle map $\operatorname{Bun}$. Its fibers therefore provide Lagragian sections,
 and the above section is one of them. Precisely, $\nabla_0$ is characterized by the fact that it is the unique connection
 (compatible with the given parabolic structure) such that the divisor of the apparent map $\operatorname{App}$
 in (\ref{eq:App-P1}) is precisely $\lambda + t$.
\end{Remark}

In \cite[Section~5.1]{LoraySaito2015}, the authors show that any connection on a generic parabolic bundle defined by parameters $u_\lambda$, $u_t$ can be written uniquely as
\begin{equation}\label{eq:univ-family}
 \nabla = \nabla_0(u_\lambda, u_t) + c_1\Theta_1(u_\lambda, u_t) + c_2\Theta_2(u_\lambda, u_t),
\end{equation}
for some $(c_1,c_2)\in\C^2$. This follows from the fact that for any $(u_\lambda, u_t)$ the Higgs bundles $\Theta_i(u_\lambda, u_t)$ are linearly independent over $\C$.
Note that the above description defines a birational map
\begin{equation}\label{eq:coords-Con-P1}
 \Con^{\bar\mu}_{\bar\nu}\big(\P^1,D\big) \dashrightarrow \P^1_{u_\lambda} \times \P^1_{u_t} \times \C^2_{(c_1,c_2)} .
\end{equation}

\begin{Definition}\label{def:universal-family-P1}
 We define the \textit{universal family over} $\P^1$ as the explicit algebraic family of connections $\mathcal{U} = \{(E_\theta, \nabla_\theta, \ell_\theta)\}$, where $\theta = (u_\lambda, u_t, c_1, c_2) \in \C^4$ and
 \begin{itemize}\itemsep=0pt
 \item $E_\theta = \mathcal{O}_{\P^1} \oplus \mathcal{O}_{\P^1}$ for all $\theta$,
 \item $\nabla_\theta = \nabla_0(u_\lambda, u_t) + c_1\Theta_1(u_\lambda, u_t) + c_2\Theta_2(u_\lambda, u_t)$,
 \item $\ell_\theta = (0,1,u_\lambda,\infty,u_t)$.
 \end{itemize}
\end{Definition}

Note that we now have three isomorphic spaces that can be naturally identified to each other:
the parameter space $\C_{u_\lambda} \times \C_{u_t} \times \C^2_{(c_1,c_2)}$, the above family~$\mathcal{U}$, and the open set $U \subset \Con^{\bar\mu}_{\bar\nu}\big(\P^1,D\big)$ which is the image of the natural map $\mathcal{U} \rightarrow \Con^{\bar\mu}_{\bar\nu}\big(\P^1,D\big)$ assigning to each connection its isomorphism class in the moduli space.

The above identifications give a coordinate chart on a Zariski~open subset of the moduli space of connections.
Moreover, we obtain, over the space of generic bundles, a trivialization of the affine $\C^2$-bundle $\Con^{\bar\mu}_{\bar\nu}\big(\P^1,D\big)^0 \to \Bun^{\bar\mu}\big(\P^1,D\big)^0$.
Indeed, if we consider the space of generic bundles $\Bun^{\bar\mu}\big(\P^1,D\big)^0$ as a subset of $\C_{u_\lambda} \times \C_{u_t}$, equation (\ref{eq:univ-family}) implies that $\Con^{\bar\mu}_{\bar\nu}\big(\P^1,D\big)^0 \cong \Bun^{\bar\mu}\big(\P^1,D\big)^0 \times \C^2_{(c_1,c_2)}$.

In the coordinates defined by the map (\ref{eq:coords-Con-P1}), the Goldman--Iwasaki symplectic structure is given by the 2-form
\[
 \omega_{\P^1} = {\rm d}c_1\wedge {\rm d}u_t + {\rm d}c_2\wedge {\rm d}u_\lambda.
\]
(see \cite[p.~1031]{LoraySaito2015}).
In general, the moduli space of parabolic Higgs bundles is naturally identified with the total space of the cotangent bundle to the moduli space of parabolic bundles~\cite{ArinkinLysenko1997A, ArinkinLysenko1997B}.
As explained in \cite[Section~5.1]{LoraySaito2015}, under this correspondence we have $\Theta_1 \mapsto {\rm d}u_t$, $\Theta_2 \mapsto {\rm d}u_\lambda$.
Note that the 1-form $\omega_{\P^1}$ given above corresponds with the canonical symplectic form on the cotangent bundle $T^*\Bun^{\bar\mu}\big(\P^1,D\big)$, using coordinates $(u_\lambda, u_t)$ to describe $\Bun^{\bar\mu}\big(\P^1,D\big)$, and coordinates $(c_1,c_2)$ to describe the cotangent space at a given point $(u_\lambda, u_t)$.

\section{Statement of the main results} \label{sec:main-body}

We begin with the map $\phi \colon \Bun^{\bar\mu}\big(\P^1,D\big) \to \Bun^{\bfmu}(C,T)$ between moduli spaces of parabolic bundles.
Our main objective is to describe the space $\Con^{\bfmu}_{\bfnu}(C,T)$ and the map obtained by extending $\phi$ to the moduli spaces of connections (which we denote $\Phi$).
We know little about the space $\Con^{\bfmu}_{\bfnu}(C,T)$, except it is an affine $\C^2$-bundle over the space $\Bun^{\bfmu}(C,T)$.
However, we know beforehand from Proposition~\ref{prop:PhiSymplecticMonodromy} that the extended map $\Phi$ is dominant and generically~$2:1$, and if follows from Theorem~\ref{thm:PhiMonodromySide} that the space $\Con^{\bfmu}_{\bfnu}(C,T)$ can be obtained as a quotient of $\Con^{\bar\mu}_{\bar\nu}\big(\P^1,D\big)$ by the Galois involution $\Psi$ associated to the map $\Phi$.

\subsection{A ``double'' universal family of connections}\label{sec:double-family} \label{subsec:universal-family}

Recall that $\Psi$ acts on parabolic connections on $\big(\P^1,D\big)$, as is defined in Remark~\ref{rmk:definition-Psi}, by performing positive elementary transformations at $0$, $1$, $\lambda$, $\infty$, and a twist by a rank~1 connection. Moreover, Lemma~\ref{lemma:psi-Galois} shows that such transformation corresponds to the Galois involution of the map $\Phi$ defined in Section~\ref{subsec:Phi-connections}.
Namely, $\Psi$ is not the identity and it satisfies $\Phi\circ\Psi = \Phi$.

We can use the transformation $\Phi$ to pull the universal family of connections $\mathcal{U}$, described in Definition~\ref{def:universal-family-P1}, from $\P^1$ to $C$.
We do this by simply applying $\Phi$ to every element $(E_\theta, \nabla_\theta, \ell_\theta) \in \mathcal{U}$. Explicitly, we pullback each connection to $C$ using the elliptic cover $\pi\colon C \to \P^1$, we then perform elementary transformations at the torsion points~$\bfw_k$, and finally we twist by a rank~1 connection on $\mathcal{O}_C(-2\bfw_\infty)$.
This defines a family of connections on $(C,T)$, which is also parametrized by $\theta \in \C_{u_\lambda} \times \C_{u_t} \times \C^2_{(c_1,c_2)}$.

\begin{Definition}\label{def:universal-family-C}
 We define the \textit{double universal family over} $C$ as the explicit algebraic family of connections $\mathcal{U}_C = \{\Phi(E_\theta, \nabla_\theta, \ell_\theta)\}$, where $(E_\theta, \nabla_\theta, \ell_\theta) \in \mathcal{U}$.
\end{Definition}

This family of connections represents every generic class in the moduli space $\Con^{\bfmu}_{\bfnu}(C,T)$, but each class has two representatives in~$\mathcal{U}_C$.
This follows from the fact that we already know from Corollary~\ref{cor:ConIsSingular} that~$\Phi$, as a map between moduli spaces of connections, is dominant and generically $2:1$.
In this way, the family $\mathcal{U}_C$ is not birational to the space $\Con^{\bfmu}_{\bfnu}(C,T)$, but a~double cover of it. Thus we call it a~``double'' universal family.

It is important to note that (the image of) the original basis $\nabla_0(u_\lambda, u_t)$, $\Theta_1(u_\lambda, u_t)$, $\Theta_2(u_\lambda, u_t)$ in (\ref{eq:univ-family}) is not the most suitable for describing the family~$\mathcal{U}_C$.
Indeed, these are not equivariant with respect to the Galois involutions $\Psi$, $\psi$, which are fundamental in the description of $\Con^{\bfmu}_{\bfnu}(C,T)$ (cf.~Remark~\ref{rmk:quotient-psi}).
Let us clarify the above claim.
Given a parabolic bundle defined by $u=(u_\lambda, u_t)$, it is not true in general that $\Psi\nabla_0(u) = \nabla_0(\psi u)$ nor $\Psi\Theta_1(u) = \Theta_1(\psi u)$ (incidentally, $\Theta_2$ is always equivariant).
We now seek for a basis that is equivariant with respect to these involutions. This will descend as a basis to the quotient $\mathcal{U} / \Psi \stackrel{\sim}{\dashrightarrow} \Con^{\bfmu}_{\bfnu}(C,T)$.

\begin{Remark}
 The connection $\Psi\nabla_0(u)$ is a parabolic connection whose underlying parabolic structure is given by~$\psi u$. The latter connection can be expressed uniquely as a combination $\nabla_0(\psi u) + c_1\Theta_1(\psi u) + c_2\Theta_2(\psi u)$, for some coefficients $c_1$, $c_2$.
 The claim that the basis is not equivariant means that the coefficients $c_1$, $c_2$ are not both zero. Equivalently, we could state that $\Psi\nabla_0(\psi u) \neq \nabla_0(u)$, even though both left hand side and right hand side are connections with parabolic structure given by $\psi (\psi u) = u$.
\end{Remark}

Note that because $\Psi\nabla_0(\psi u)$ is a connection with underlying parabolic structure~$u$, the mean $\frac{1}{2} (\nabla_0(u) + \Psi\nabla_0(\psi u))$ also has $u$ as its parabolic structure.
Moreover, this average is equivariant with respect to~$(\Psi,\psi)$.

\begin{Definition}\label{def:equivariant-basis}
 We define the following elements of $\mathcal{U}$:
 \[
 \nabla_0^\psi(u) = \tfrac{1}{2} (\nabla_0(u) + \Psi\nabla_0(\psi u)), \qquad
 \Theta_i^\psi(u) = \tfrac{1}{2} (\Theta_i(u) + \Psi\Theta_i(\psi u)), \qquad i=1,2,
 \]
 and of $\mathcal{U}_C$:
 \[
 \bfnabla_0^\psi(u) = \Phi(\nabla_0^\psi(u)), \qquad\text{and}\qquad \bfTheta_i^\psi(u) = \Phi(\Theta_i^\psi(u)), \qquad i=1,2,
 \]
 and we call them the \textit{equivariant bases} for $\mathcal{U}$ and $\mathcal{U}_C$, respectively.
\end{Definition}

By construction, all the above are equivariant with respect to $(\Psi,\psi)$. The underlying parabolic structure of $\nabla_0^\psi(u)$ is precisely~$u$, and the parabolic structure of $\bfnabla_0(u)$ is~$\Phi(u)$.

By combining $\bfTheta_1^\psi$ and $\bfTheta_2^\psi$ in a suitable manner, we arrive to a new equivariant basis $\bfnabla_0^\psi$, $\bfTheta_z$, $\bfTheta_w$, which is more convenient. The exact expression for these elements can be found in Definition~\ref{def:canonical-basis}. We shall refer to this as the \emph{canonical basis} (it comes from the canonical identification of the moduli space of parabolic Higgs bundles and the cotangent bundle of~$\Bun^{\bfmu}(C,T)$, cf.~Section~\ref{subsec:canonical-basis}).
This is the basis we will use to describe the family $\mathcal{U}_C$ (it is only well defined for generic connections). Each generic element of $\mathcal{U}_C$ can be written uniquely as a~linear combination $\bfnabla_0^\psi(u) + \kappa_1 \bfTheta_z(u) + \kappa_2 \bfTheta_w(u)$.
The change of basis acts on the coefficients in such a way that the map $(c_1,c_2) \mapsto (\kappa_1,\kappa_2)$ is given by an affine transformation which depends rationally on $u$.

As pointed out earlier, the map $\mathcal{U}_C \to \Con^{\bfmu}_{\bfnu}(C,T)$ is dominant and generically $2:1$.
Let us denote by $\bfPsi$ the transformation that permutes the fibers.
Note that $\Psi$ and $\bfPsi$ represent essentially the same transformation, except that $\Psi$ is acting on connections on $\big(\P^1,D\big)$, and $\bfPsi$ on connections on $(C,T)$.
We obtain the following commutative diagram
\[
 \begin{tikzcd}
 \mathcal{U} \arrow[r, "\cong"] \arrow[d, "\Psi "]
 &
 \mathcal{U}_C \arrow[d, "\bfPsi"] \arrow[r, "2:1"]
 &
 \Con^{\bfmu}_{\bfnu}(C,T) \arrow[d, "id"]
 \\
 \mathcal{U} \arrow[r, swap, "\cong"]
 &
 \mathcal{U}_C \arrow[r, swap, "2:1"]
 &
 \Con^{\bfmu}_{\bfnu}(C,T) .
 \end{tikzcd}
\]

We conclude that
\begin{equation}\label{eq:birat-equiv0}
 \mathcal{U}_C / \bfPsi \stackrel{\sim}{\dashrightarrow} \Con^{\bfmu}_{\bfnu}(C,T) .
\end{equation}
Since $\mathcal{U}_C$ is isomorphic to $\mathcal{U}$, and the latter is birational to $\Con^{\bar\mu}_{\bar\nu}\big(\P^1,D\big)$, we can further write
\begin{equation*}
 \Con^{\bar\mu}_{\bar\nu}\big(\P^1,D\big) / \Psi \stackrel{\sim}{\dashrightarrow} \Con^{\bfmu}_{\bfnu}(C,T) ,
\end{equation*}
as anounced at the beginning of this section.
This fact is analogous to the description of $\Bun^{\bfmu}(C,T)$ as the quotient $\Bun^{\bar\mu}\big(\P^1,D\big) / \psi$ in Remark~\ref{rmk:quotient-psi}. However, we remark that, unlike the case of parabolic bundles, the fixed-point set of $\Psi$ has codimension bigger than one, making the space $\Con^{\bfmu}_{\bfnu}(C,T)$ singular. This is discussed in Theorem~\ref{thm:fixed-points-psi} below, and was previously announced in Corollary~\ref{cor:ConIsSingular}.

Because of equivariance of the basis $\bfnabla_0^\psi$, $\bfTheta_z$, $\bfTheta_w$, we have that $\bfPsi$ acts only on the first factor of $\Bun^{\bar\mu}\big(\P^1,D\big) \times \C^2_{(\kappa_1,\kappa_2)}$.
Namely,
\[
 \bfPsi(\bfnabla_0^\psi(u) + \kappa_1\bfTheta_z(u) + \kappa_2\bfTheta_w(u)) = \bfnabla_0^\psi(\psi u) + \kappa_1\bfTheta_z(\psi u) + \kappa_2\bfTheta_w(\psi u),
\]
and so, in coordinates, $\bfPsi \colon (u, \kappa_1, \kappa_2) \mapsto (\psi u, \kappa_1, \kappa_2)$.
We conclude that
\begin{equation}\label{eq:birat-equiv1}
 \mathcal{U}_C/\bfPsi \stackrel{\sim}{\dashrightarrow} \Bun^{\bar\mu}\big(\P^1,D\big) / \psi \times \C^2_{(\kappa_1,\kappa_2)} .
\end{equation}
From the description of the moduli space of parabolic bundles in \cite{FernandezVargas2016}, we have that
\begin{equation}\label{eq:birat-equiv2}
 \Bun^{\bar\mu}\big(\P^1,D\big)/\psi \cong \Bun^{\bfmu}(C,T) \cong \P^1_z\times\P^1_w .
\end{equation}
Finally, combining (\ref{eq:birat-equiv0})--(\ref{eq:birat-equiv2}) we obtain
\begin{equation} \label{eq:coords-ConCT}
 \Con^{\bfmu}_{\bfnu}(C,T) \stackrel{\sim}{\dashrightarrow} \P^1_z\times\P^1_w \times \C^2_{(\kappa_1,\kappa_2)} .
\end{equation}
Note that this gives, over some open and dense subset of $\Bun^{\bfmu}(C,T)$, a coordinate system, as well as a local trivialization of the affine $\C^2$-bundle $\Con^{\bfmu}_{\bfnu}(C,T) \to \Bun^{\bfmu}(C,T)$.

We recall the reader that $\Con^{\bfmu}_{\bfnu}(C,T)^0$ denotes the open subset of \emph{generic} connections, as specified in Definition~\ref{def:generic-bunle-C} (see also Definition~\ref{def:generic-bunle-P1} for connections on $\P^1$).

\begin{Theorem}\label{thm:Bun-C}
 The family of connections $\bfnabla_0^\psi$ in Definition~{\rm \ref{def:equivariant-basis}} defines a global section of the affine bundle $\operatorname{Bun}$, namely, $\bfnabla_0^\psi \colon \Bun^{\bfmu}(C,T)^0 \to \Con^{\bfmu}_{\bfnu}(C,T)^0$.
 This identifies the affine bundle $\Con^{\bfmu}_{\bfnu}(C,T)^0 \to \Bun^{\bfmu}(C,T)^0$ to the vector bundle $\Higgs^{\bfmu}(C,T)^0 \to \Bun^{\bfmu}(C,T)^0$.
 Moreover, the section is Lagrangian and the above identification is symplectic with respect to the natural symplectic structures on the moduli spaces of connections and Higgs bundles.
 Finally, the latter vector bundle is algebraically trivial, namely
 \[
 \Con^{\bfmu}_{\bfnu}(C,T)^0 \cong \Bun^{\bfmu}(C,T)^0 \times \C^2.
 \]
\end{Theorem}

That the section $\bfnabla_0^\psi$ is Lagrangian and that the identification is symplectic is proved in Section~\ref{subsec:canonical-basis}.
The proof that such a section is well defined over the space $\Bun^{\bfmu}(C,T)^0$, and that the affine bundle $\operatorname{Bun}$ is algebraically trivial is given in Section~\ref{sec:final-steps}.

\begin{Corollary} \label{coro:omega-C}
 Under the canonical identification between $T^*\Bun^{\bfmu}(C,T)$ and $\Higgs^{\bfmu}(C,T)$, we have the correspondence ${\rm d}z \mapsto \bfTheta_z$, ${\rm d}w \mapsto \bfTheta_w$, and so the symplectic structure of $\Con^{\bfmu}_{\bfnu}(C,T)$ is given by the $2$-form
 \[
 \omega_C = {\rm d}\kappa_1\wedge {\rm d}z + {\rm d}\kappa_2\wedge {\rm d}w.
 \]
\end{Corollary}

The proof of this corollary can be found in Section~\ref{subsec:canonical-basis}.

\begin{Remark} \label{rmk:Phi_hat-in-coords}
 The open space of generic parabolic connections $\Con^{\bfmu}_{\bfnu}(C,T)^0$ is smooth.
 Therefore the map $\Phi$ restricted to $\Phi^{-1}\big(\Con^{\bfmu}_{\bfnu}(C,T)^0\big) \subset \Con^{\bar\mu}_{\bar\nu}\big(\P^1,D\big)^0$ is smooth.
 An important contribution of the present paper is providing explicit expressions for this map in terms of the coordinates defined by (\ref{eq:coords-ConCT}).
 The formulas are given in Section~\ref{sec:comp-coords}.
\end{Remark}

The full space $\Con^{\bfmu}_{\bfnu}(C,T)$ is singular.
The coordinate system used above fails to describe the space $\Con^{\bfmu}_{\bfnu}(C,T)$ precisely because the latter is singular over $\boldsymbol\Sigma$.
Recall that our description of the moduli space $\Con^{\bfmu}_{\bfnu}(C,T)$ is based on the fact that such space can be realized as the quotient of $\Con^{\bar\mu}_{\bar\nu}\big(\P^1,D\big)$ by the involution $\Psi$.

\begin{Theorem}\label{thm:fixed-points-psi}
 The set of fixed points $\Fix(\Psi) \subset \Con^{\bar\mu}_{\bar\nu}\big(\P^1,D\big)^0$ is a codimension~$2$ sub\-va\-riety.
 In fact, it defines a subbundle of rank~$1$ of the affine $\C^2$ bundle $\Bun\vert_{\Sigma}$.
 A connection $(\bfnabla, \mathbf{E}, \{\ell_1, \ell_2\}) \in \Con^{\bfmu}_{\bfnu}(C,T)$ belongs to $\Phi(\operatorname{Fix} (\Psi))$ if and only if it decomposes as a direct sum
 \[
 \big(\bfnabla, \mathbf{E}, \bar\ell\big) =
 \big(\eta, \mathbf{L}, \ell\big) \oplus \iota^* \big(\eta, \mathbf{L}, \ell\big),
 \]
 where $\eta$ is a rank $1$ connection with a single pole at either $\bft_1$ or $\bft_2$, $\mathbf{L}$ is a line bundle of degree zero, and $\iota\colon C\to C$ is the elliptic involution.
\end{Theorem}

Evidently, the connection $\big(\bfnabla, \mathbf{E}, \bar\ell\big)$ has each direct summand of $\mathbf{E}$ as an invariant subbundle. However, the parabolic degree of $(\mathbf{L}, \ell)$ is zero (cf.~Definition~\ref{def:stability}). There are no more invariant subbundles and so the connection is (strictly) semistable, which means that it does belong to the moduli space $\Con^{\bfmu}_{\bfnu}(C,T)$.

We remark that the fact that the fixed-point set of $\Psi$ has codimension bigger than one causes the space $\Con^{\bfmu}_{\bfnu}(C,T)$ to have singularities at $\Phi(\operatorname{Fix}(\Psi))$. The next theorem gives the local description of the singular set.

\begin{Theorem}\label{thm:Sing-locally} Around a generic point of the singular locus $\Phi(\operatorname{Fix}(\Psi))$, the space $\Con^{\bfmu}_{\bfnu}(C,T)$ is locally isomorphic to the hypersurface in $\C^5_{\bar x}$ given by the equation $x_0^2 = x_1 x_2$.
 Thus, locally, the singularities look like the product of a quadratic conic singularity and a~bidisk.
 Moreover, the symplectic structure is given by
 \[
 \tilde{\omega} = \frac{{\rm d}x_1\wedge {\rm d}x_2}{4x_0}+{\rm d}x_3\wedge {\rm d}x_4 .
 \]
\end{Theorem}

Since, in the above coordinates, the singular locus is given by $x_0=x_1=x_2=0$, we conclude that the symplectic form restricted to the singular locus is ${\rm d}x_3\wedge {\rm d}x_4$. Thus the singular locus is a symplectic subvariety of $\Con^{\bfmu}_{\bfnu}(C,T)$.

The proof of the last two theorems can be found in Section~\ref{subsec:singular-locus}.

\section{Geometric description of a generic parabolic bundle} \label{sec:geom-pic}

In this section we want to describe some geometric aspects of parabolic bundles in $\Bun^{\bfmu}(C,T)^0$.

We begin with a generic parabolic connection $(E,\nabla,\bar\ell)$ of degree zero and polar divisor $D=0+1+\lambda+\infty+t$ over $\P^1$. As pointed out in Proposition~\ref{prop:intrinsic-def-generic}, in the generic case we can assume $E=\mathcal{O}_{\P^1} \oplus \mathcal{O}_{\P^1}$. Introducing an affine coordinate $\zeta$ on the fibers of $\P(E)$, the parabolic structure is given by $\bar\ell = (0,1,u_\lambda,\infty,u_t)$, for some $u_\lambda, u_t\in\P^1\setminus \{0,1,\lambda,\infty\}$. We denote by $\sigma_\infty$ the constant horizontal section $\zeta=\infty$.
This section plays a central role, since it is used to define the apparent map and the universal family $\mathcal{U}$ (cf.~Sections~\ref{sec:connections-P1} and~\ref{subsec:universal-family-P1}). The birational involution~$\psi $ defined in Remark~\ref{rmk:definition-psi} preserves the trivial bundle $E$ and transforms $\sigma_\infty$ into a~section~$\sigma_\psi$ of self-intersection $+2$. This section passes through the parabolic points over $0$, $1$, $\lambda$ but not~$\infty$. The sections $\sigma_\infty$ and $\sigma_\psi$ intersect (transversally) over a unique point $x=p$.

Now we apply the transformation $\Phi$ to obtain a parabolic connection on $(C,T)$ with trivial determinant.
This time the underlying vector bundle $\mathbf E$ will not be trivial, and it depends on the parabolic structure of $(E,\bar\ell)$. Generically, it is of the form $\mathbf E = \mathbf L \oplus \mathbf L^{-1}$, where~$\mathbf L$ is a line bundle of degree zero. As such, there exists a unique point $\bfp_1\in C$ such that $\mathbf L=\mathcal{O}_C(\bfp_1-\bfw_\infty)$.
The other summand is given by $\mathbf L^{-1} = \mathcal{O}_C(\bfp_2-\bfw_\infty)$, in such a way that $\bfp_1$, $\bfp_2$ is a pair of points in involution. Moreover, these points project to the point $p$ defined at the end of the last paragraph (i.e., $\pi^{-1}(p) = \{ \bfp_1, \bfp_2\}$).
The bundle $\mathbf E$ contains two subbundles $S_\infty$, $S_\psi$ defining sections on $\mathbb P(\mathbf E)$,
which are the images of $\sigma_\infty$, $\sigma_\psi$, respectively.
These sections are exchanged by an automorphism of the bundle $\mathbf E$, which we denote $\bfPsi$, making the following diagram commute
\[
 \begin{tikzcd}
 \mathbf{E} \arrow[r, dashed, "\Phi"] \arrow[d, swap, "\bfPsi "]
 &
 E \arrow[d, "\Psi"]
 \\
 \mathbf{E} \arrow[r, swap, dashed, "\Phi"]
 &
 E .
 \end{tikzcd}
\]

\begin{Remark}\label{rmk:multisection-S_Sigma}
 The subbundles $\mathbf L$, $\mathbf L^{-1}$ do not come from subbundles of~$E$ in~$\P^1$. Rather, there exists a~curve~$S_\Sigma$ in $\P(E) \cong \P^1_x\times\P^1_\zeta$ of bidegree~$(2,2)$ projecting down to~$\P^1_x$ as a double cover ramified at the points $x=0,1,\lambda,\infty$ (hence isomorphic to the elliptic curve~$C$ itself).
 In fact, this curve is defined by the fixed points of the transformation~$\psi $.
 See \cite[Section~3.3, Theorem~7]{Loray2016} for details on this construction.
 Once pulled back to $C$ via $\pi$, this curve splits into two different sections~$S_\Sigma^+$, $S_\Sigma^-$ of $\P(\pi^*(E))$ that intersect over the points $\bfw_0$, $\bfw_1$, $\bfw_\lambda$, $\bfw_\infty$. After the elementary transformations dictated by $\phi$ (step (\ref*{item:ii}) in Section~\ref{subsec:Phi-bundles}), the sections no longer intersect in~$\P(\mathbf{E})$.
 The action of $\phi$ on~$S_\Sigma$ can be seen in~Fig.~\ref{fig:Phi-hat} at the end of Section~\ref{sec:pullback-map}.
\end{Remark}

\begin{figure}[t] \centering
 \includegraphics[width=0.9\linewidth]{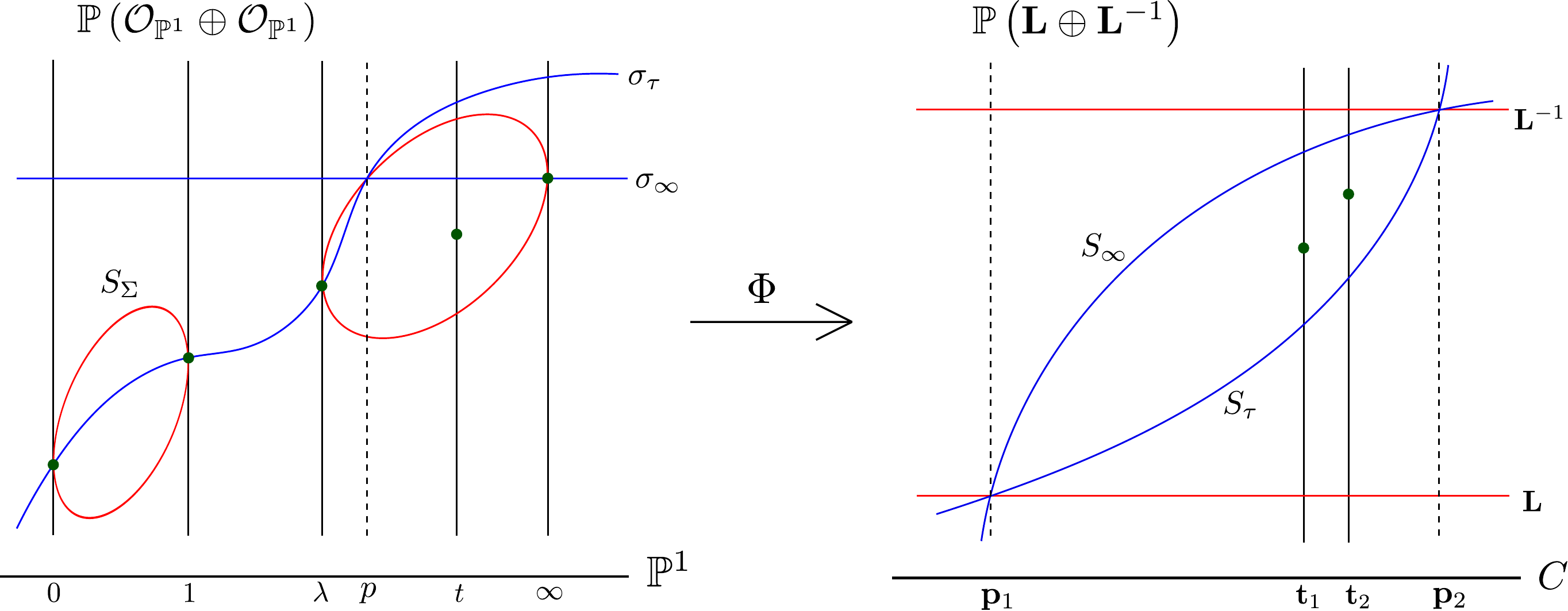}
 \caption{The sections $\sigma_\infty$ and $\sigma_\psi$ under $\Phi$.} \label{fig:4sections}
\end{figure}

The curve $S_\Sigma$ passes through the parabolic points above $0$, $1$, $\lambda$, $\infty$. In general, it does not pass through the parabolic above~$t$. When it does, the parabolic structure is unchanged by~$\psi $.
Therefore the fixed points of $\psi \in \operatorname{Birat} \big(\P^1_{u_\lambda} \times \P^1_{u_t}\big)$ are precisely those for which \mbox{$(t,u_t) \in S_\Sigma$}.
Moreover, if the curve $S_\Sigma$ passes through the parabolic point above~$t$ in $\P^1$, then, after performing the transformation $\phi$, the parabolics over~$T$ are in either~$\mathbf L$ or~$\mathbf L^{-1}$.
The elliptic involution permutes the summands $\mathbf L \oplus \mathbf L^{-1}$.
Since the parabolic bundles that come from~$\P^1$ via~$\phi$ are invariant under the elliptic involution, we conclude that each direct summand contains exactly one parabolic point. In such case we have a direct sum decomposition $(\mathbf{E}, \{\ell_1, \ell_2\}) = (\mathbf{L}, \{\ell_1\}) \oplus \big(\mathbf{L}^{-1}, \{\ell_2\}\big)$.
Thus we conclude that the generic elements of $\boldsymbol\Sigma$ are precisely the decomposable parabolic bundles.

\section{Computations in coordinates} \label{sec:comp-coords}

Recall that in Section~\ref{subsec:universal-family-P1} we discussed how to define a coordinate sytem $\P^1_{u_\lambda} \times \P^1_{u_t} \times \C^2_{(c_1,c_2)} \stackrel{\sim}{\dashrightarrow} \Con^{\bar\mu}_{\bar\nu}\big(\P^1,D\big)$, using the basis $\nabla_0$, $\Theta_1$, $\Theta_2$ of the universal family $\mathcal{U}$.

In this section we will first define the canonical basis $\bfnabla_0^\psi$, $\bfTheta_z$, $\bfTheta_w$, and use it to establish a~coordinate system $\P^1_z \times \P^1_w \times \C^2_{(\kappa_1,\kappa_2)} \stackrel{\sim}{\dashrightarrow} \Con^{\bfmu}_{\bfnu}(C,T)$. We obtain the following diagram:
\[
 \begin{tikzcd}
 \Con^{\bar\mu}_{\bar\nu}\big(\P^1,D\big) \arrow[r, dashrightarrow] \arrow[d, swap, "\Phi"]
 &
 \P^1_{u_\lambda} \times \P^1_{u_t} \times \C^2_{(c_1,c_2)} \arrow[d, "\Phi"]
 \\
 \Con^{\bfmu}_{\bfnu}(C,T) \arrow[r, dashrightarrow]
 &
 \P^1_z\times\P^1_w \times \C^2_{(\kappa_1,\kappa_2)} .
 \end{tikzcd}
\]

The horizontal arrows are birational isomorphisms. In order to describe the left vertical map~$\Phi$ between moduli spaces, we will compute explicitly the right arrow $\Phi$ in the given coordinates.

First, let us split the map $\Phi$ as follows: $\Phi = (\phi, \mathrm{T}_\phi)$, where $\phi \colon \P^1_{u_\lambda} \times \P^1_{u_t} \to \P^1_z \times \P^1_w$ is the map between moduli spaces of parabolic bundles, and $\mathrm{T}_\phi$ defines fiberwise an affine map $\C_{(c_1,c_2)} \to \C_{(\kappa_1, \kappa_2)}$.

\begin{Proposition}\label{prop:formula-Nestor}
 The map $\phi \colon \P^1_{u_\lambda} \times \P^1_{u_t} \dashrightarrow \P^1_z\times\P^1_w$ is given by $(u_\lambda, u_t) \mapsto (z, w)$, where
 \begin{equation}\label{eq:formula-Nestor}
 z = \frac{\lambda (u_{\lambda} - 1)}{u_{\lambda} - \lambda} , \qquad
 w = \frac{\lambda u_{t}(\lambda u_{t} - t u_{\lambda} + t - \lambda - u_{t} + u_{\lambda})}
 {t \lambda u_t - t \lambda u_\lambda - t u_t u_\lambda + \lambda u_t u_\lambda - \lambda u_t + t u_\lambda} .
 \end{equation}
\end{Proposition}

\begin{Remark}
 The above map $\phi$ has been explicitly described both in coordinates and in geometric terms in \cite{FernandezVargas2016}.
 Unfortunately, there is a mistake in the final expression for this map in the coordinates used above. Namely, the equations given in \cite[Proposition~6.3]{FernandezVargas2016} are wrong.
 We have recomputed the formulas following the steps in the cited paper to obtain (\ref{eq:formula-Nestor}).
 An alternative approach yielding the same formulas can be found in our repository \cite{github}.
\end{Remark}

\subsection[The involution psi and families of non-generic bundles in coordinates]{The involution $\boldsymbol{\psi}$ and families of non-generic bundles in coordinates} \label{sec:involution-exceptional-curves}

The next proposition describes the involution $\psi$ in terms of the $(u_\lambda, u_t)$ coordinates.
In particular, we have that $\psi$ fixes the variable $u_\lambda$ thus acting only on~$u_t$. We will use the notation $(u_\lambda, \overline{u}_t)$ to denote the image of $(u_\lambda, u_t)$ under $\psi$.
We remark that the value of $\overline{u}_t$ does depend on~$u_\lambda$, but as long as $u_\lambda$ is fixed, the map $u_t \mapsto \overline{u}_t$ is an involution (i.e., $\overline{\overline{u}}_t = u_t$), hence the choice of this notation.

\begin{Proposition}\label{prop:u_t-bar}
 The involution $\psi\in\operatorname{Birat}\big(\P^1_{u_\lambda} \times \P^1_{u_t}\big)$ is given by $(u_\lambda, u_t) \mapsto (u_\lambda, \overline{u}_t)$, where~$\overline{u}_t$ is given by
 \begin{equation}\label{eq:u_t-bar}
 \overline{u}_t = \frac{t u_{\lambda}(\lambda u_{t} - t u_{\lambda} + t - \lambda - u_{t} + u_{\lambda})}
 {t \lambda u_t - t \lambda u_\lambda - t u_t u_\lambda + \lambda u_t u_\lambda - \lambda u_t + t u_\lambda} .
 \end{equation}
\end{Proposition}

\begin{proof} Let $(z,w) = \phi(u_\lambda, u_t)$. We know that $(z,w)$ has another preimage under~$\phi$, which by definition is $\psi(u_\lambda, u_t)$. Therefore, we need to solve $\phi(u'_\lambda, u'_t) = (z,w)$ for~$u'_\lambda$,~$u'_t$. From~(\ref{eq:formula-Nestor}) we can see that $z$ is uniquely determined by $u_\lambda$, and that fixing the value for $w$ imposes a quadratic condition on $u_t$. Solving the equation we recover~(\ref{eq:u_t-bar}).
\end{proof}

\begin{Remark} The polynomial
 \begin{equation}\label{eq:Pi}
 P_\Pi = t \lambda u_t - t \lambda u_\lambda - t u_t u_\lambda + \lambda u_t u_\lambda - \lambda u_t + t u_\lambda,
 \end{equation}
 which appears as the denominator of~(\ref{eq:u_t-bar}) defines a rational curve in $\P^1_{u_\lambda} \times \P^1_{u_t}$, which corresponds to the rational curve $\Pi \subset \mathcal{S}$ introduced in Section~\ref{sec:genericity} (and corresponds to the conic $b_1^2 - b_0b_2$ in the birational model~$\P^2_b$ discussed in Remark~\ref{rmk:coordinates-Bun-P2}).
 The involution $\psi$ permutes this curve with the line $u_t=\infty$, thus the former appears as a pole of $\overline{u}_t$.
 These rational curves correspond to the two $(-1)$-curves in the Del Pezzo surface $\S$ that are mapped by $\phi$ to the horizontal line $w=\infty$. Therefore, the polynomial $\Pi$ also appears in the denominator of the formula for $w$ in~(\ref{eq:formula-Nestor}).
 Indeed, comparing (\ref{eq:formula-Nestor}) to (\ref{eq:u_t-bar}), we can see that $w = \frac{\lambda}{tu_\lambda} u_t \overline{u}_t$.
\end{Remark}

The last equation above shows that we can write the product $u_t \overline{u}_t$ in terms of~$w$.
Moreover, we can also write $u_t + \overline{u}_t$ in terms of~$z$,~$w$.
Explicitly,
\[
 u_t + \overline{u}_t = \frac{zw + (z-w)t - \lambda}{z-\lambda}, \qquad u_t\overline{u}_t = \frac{wt(z-1)}{z-\lambda}.
\]
Finally, we can also express $u_\lambda$ as a fractional linear transformation of $z$.
This shows that any rational function on $u_\lambda$,~$u_t$ that is invariant under $\psi $ can be expressed as a rational function on~$z$,~$w$.

\begin{Remark}
 In the chart $\P^1_{u_\lambda} \times \P^1_{u_t}$ the curve $\Sigma$ which was defined as the ramification locus of $\phi$ is defined by the zeros of the following polynomial:
 \begin{equation}\label{eq:Sigma}
 P_\Sigma = t \lambda u_t^2 - 2 t \lambda u_t u_\lambda - t u_t^2 u_\lambda + \lambda u_t^2 u_\lambda + t^2 u_\lambda^2 - \lambda u_t^2 - t^2 u_\lambda + t \lambda u_\lambda + 2 t u_t u_\lambda - t u_\lambda^2.
 \end{equation}
\end{Remark}

\begin{Remark}\label{rmk:Lambda}
 In Section~\ref{sec:genericity} we have defined $\boldsymbol\Lambda \subset \P^1_z \times \P^1_w$ to be the vertical line given by $z=t$, and $\Lambda$ its preimage $\Lambda = \phi^{-1}(\boldsymbol\Lambda) \subset \Bun^{\bar\mu}\big(\P^1,D\big)$.
 In our chart $\P^1_{u_\lambda} \times \P^1_{u_t}$, the curve $\Lambda$ is defined by the vertical line $u_\lambda = \lambda(1 - t)/(\lambda - t)$, or equivalently, by the zeros of the polynomial
 \begin{equation}\label{eq:Lambda}
 P_\Lambda = \lambda u_\lambda - t u_\lambda + \lambda t - \lambda .
 \end{equation}
\end{Remark}

The special curves discussed in the above remarks can be seen in Fig.~\ref{fig:non-generic}, Section~\ref{sec:genericity}.

\subsection[Geometry of the apparent map on $\P^1$]{Geometry of the apparent map on $\boldsymbol{\P^1}$}

\looseness=-1 Before we move on, let us recall the geometric picture of the universal family $\mathcal{U}$ on $\P^1$. As usual, we consider parabolic bundles on $\big(\P^1,D\big)$ with trivial underlying bundle.
We assume that, with respect to an affine coordinate $\zeta$ on the vertical fibers of $\P\big(\mathcal{O}_{\P^1} \oplus \mathcal{O}_{\P^1}\big)$, the parabolic structure is given by $\bar\ell = (0,1,u_\lambda,\infty,u_t)$ (cf.~Remark~\ref{rmk:coordinates-Bun-P1P1}).
For each parabolic bundle represented in $\P^1_{u_\lambda} \times \P^1_{u_t}$, we define a connection $\nabla_0(u_\lambda, u_t)$.
This connection is characterized as the unique connection compatible with the parabolic structure such that the divisor of the apparent map is given by $\operatorname{App}(\nabla_0) = \lambda + t$.
We recall that the apparent map is defined by the tangencies of the Riccati foliation with the section $\sigma_\infty = \{\zeta=\infty\}$.
We will exploit the fact that every connection $\nabla\in\mathcal{U}$ is completely determined by its parabolic structure and its image under the apparent map.

In Remark~\ref{rmk:definition-psi}, $\psi$ is defined as acting on parabolic bundles as follows: Given $(E,\bar\ell)$, we obtain $(E', \bar\ell') = \psi(E, \bar\ell)$ by performing elementary transformations and twisting by a rank~1 bundle. This can be seen as a birational transformation on the total space $E\dasharrow E'$. Passing to the projectivized bundles we obtain a birational fibered map $\psi\colon\P(E)\dasharrow \P(E')$.
Since in the generic case both $E$ and $E'$ are trivial, we can interpret $\psi$ as a birational automorphism of $\P^1_x\times\P^1_\zeta$.

For generic bundles, $\psi $ acts on $\P^1_x\times\P^1_\zeta$ in such a way that it exchanges $\sigma_\infty = \{\zeta=\infty\}$ with a section $\sigma_\psi$ of self-intersection $+2$.
More precisely, we have the following.

\begin{Proposition}
 Assume $(u_\lambda, u_t)$ defines a generic parabolic bundle $(E,\bar\ell)$.
 Then $\psi \colon \P^1_x\times\P^1_\zeta \dasharrow \P^1_x\times\P^1_\zeta $ transforms $\sigma_\infty$ into the section $\sigma_\psi \colon \P^1_x \to \P^1_x \times \P^1_\zeta$ defined by
 \[
\sigma_\psi(x) = \left(x,\, \frac{u_\lambda(1 - \lambda)x}{(u_\lambda - \lambda)x - \lambda(u_\lambda - 1)} \right) .
 \]
\end{Proposition}

\begin{proof}
 Recall that $\psi $ is the transformation obtained by performing elementary transformations at the parabolic points above the divisor $W=0+1+\lambda+\infty$, and subsequently performing a~twist by the bundle $\mathcal{O}_{\P^1}(-2)$ (cf.~Remark~\ref{rmk:definition-psi}). Since $\sigma_\infty$ is a section of degree zero (constant) passing only through the parabolic above $x=\infty$, we conclude that $\sigma_\psi$ must be a section of self-intersection~$+2$ passing through the parabolics above~$0$, $1$, $\lambda$ (but not~$\infty$). A simple computation shows that there is a unique such section and it is given by the expression above.
\end{proof}

\begin{Remark}\label{rmk:characterization-psi-nabla}
 Let $\nabla\in\mathcal{U}$ be defined by the parabolic structure $(u_\lambda, u_t)$ and divisor $Z$ for the apparent map. Then $\Psi\nabla$ is the unique connection with parabolic structure $\psi(u_\lambda, u_t) = (u_\lambda, \overline{u}_t)$ and whose tangencies with the section $\sigma_\psi$ are exactly given by the divisor $Z$.
\end{Remark}

In order to be more explicit, let us denote $\operatorname{App}_\infty$ the usual apparent map with respect to the constant section $\sigma_\infty$. Now that we have a formula for $\sigma_\psi$, we can compute the tangencies of this section with a given connection. Let us detail this construction.

\begin{Definition}
 Given a connection $\nabla\in\mathcal{U}$, define the vector $v(x) = (1, \sigma_\psi(x))^\top$.
 Let $v_1 = v(x)$ and $v_2 = \nabla v_1$.
 We define $\operatorname{App}_\psi(\nabla)$ as the numerator of the rational expression
 $\det \big( v_1 , v_2 \big)$.
 We call the map $\nabla \mapsto \operatorname{App}_\psi(\nabla)$ the \emph{apparent map} with respect to $\sigma_\psi$.
\end{Definition}

Explicit expressions for $\operatorname{App}_\infty$ in terms of the variables $u_\lambda$, $u_t$, $c_1$, $c_2$ are given in \cite[Section~6]{LoraySaito2015}. We omit those for $\operatorname{App}_\psi$ here since they are considerably more intricate. Explicit formulas can be found in the repository~\cite{github}.

\subsection{The base change map} \label{subsec:base-change}

Recall that we have split the map $\Phi$ between moduli spaces of connections as $\Phi = (\phi, \mathrm{T}_\phi)$, where~$\phi$ is the map between moduli spaces of parabolic bundles and $\mathrm{T}_\phi$ defines maps $\C^2_{(c_1,c_2)} \to \C^2_{(\kappa_1,\kappa_2)}$, between the fibers.

In this and the next section we describe the map $\mathrm{T}_\phi$. Since it an affine map between the fibers, it can be expressed as a $3\times3$ matrix. Furthermore, it can be understood as an affine change of coordinates.
We do this in two steps. First, in this section, we go from the \textit{original basis} (used to describe $\Con^{\bar\mu}_{\bar\nu}\big(\P^1,D\big)$ in~\cite{LoraySaito2015} and given explicitly in Section~\ref{subsec:universal-family-P1}) to the intermediate basis called the \textit{equivariant basis}.
In the next section we define the \textit{canonical basis} (which is the one we will use to describe $\Con^{\bfmu}_{\bfnu}(C,T)$), and compute the base change from the equivariant basis to the canonical one.

For the reader's convenience we repeat Definition~\ref{def:equivariant-basis} here.
The equivariant basis $\nabla_0^\psi(u)$, $\Theta_1^\psi(u)$, $\Theta_2^\psi(u)$ are the elements of $\mathcal{U}$ given by
\[
 \nabla_0^\psi(u) = \tfrac{1}{2} (\nabla_0(u) + \Psi\nabla_0(\psi u)), \qquad
 \Theta_i^\psi(u) = \tfrac{1}{2} (\Theta_i(u) + \Psi\Theta_i(\psi u)), \qquad i=1,2 .
\]
These satisfy the equivariant property $\Psi\nabla_0^\psi(u) = \nabla_0^\psi(\psi u)$, $\Psi\Theta_i^\psi(u) = \Theta_i^\psi(\psi u)$.

In order to compute the base change, we must describe in $(c_1, c_2)$-coordinates the action of the involution $\Psi$ on $\Con^{\bar\mu}_{\bar\nu}\big(\P^1,D\big)$. To do so we will exploit the idea presented in Remark~\ref{rmk:characterization-psi-nabla}.

Let $\nabla\in\mathcal{U}$ have parabolic structure $\operatorname{Bun}(\nabla) = (u_\lambda, u_t)$, and assume that with respect to the original basis for $\mathcal{U}$ it is written as
$\nabla = \nabla_0(u_\lambda, u_t) + c_1\Theta_1(u_\lambda, u_t) + c_2\Theta_2(u_\lambda, u_t)$. Following Remark~\ref{rmk:characterization-psi-nabla}, we seek for the unique connection $\nabla'$ such that
\[
 \operatorname{Bun}(\nabla') = (u_\lambda, \overline{u}_t), \qquad\text{and}\qquad
 \operatorname{App}_\psi(\nabla') = \operatorname{App}_\infty(\nabla).
\]
The connection $\nabla'$ is precisely the image of $\nabla$ under $\Psi $. A straightforward computation allows us to find coefficients $c'_1$, $c'_2$ as functions of $u_\lambda$, $u_t$, $c_1$, $c_2$, such that
\[
 \nabla' = \nabla_0(u_\lambda, \overline{u}_t) + c'_1\Theta_1(u_\lambda, \overline{u}_t) + c'_2\Theta_2(u_\lambda, \overline{u}_t) .
\]
If we hold the parabolic structure fixed, the coefficients $c'_1$, $c'_2$ are affine functions of $c_1$, $c_2$. This means that there exists a $3\times3$ matrix $\mathrm{T}_\psi(u_\lambda, \overline{u}_t)$ such that
\[
 \begin{pmatrix} 1 \\ c'_1 \\ c'_2 \end{pmatrix}
 = \mathrm{T}_\psi(u_\lambda, \overline{u}_t) \begin{pmatrix} 1 \\ c_1 \\ c_2 \end{pmatrix} .
\]

The involution $\Psi$ decomposes as $\Psi = (\psi, \mathrm{T}_\psi)$.
The symbol $\mathrm{T}_\psi$ has been chosen, in analogy to $\mathrm{T}_\phi$. The former gives the fiberwise action of the involution $\Psi$, while the latter is the fiberwise action of the map~$\Phi$.

\begin{Proposition}\label{prop:formulas-T}
 The matrix $\mathrm{T}_\psi = \mathrm{T}_\psi(u_\lambda, u_t)$ is given as follows:
 \[
 \mathrm{T}_\psi =
 \begin{pmatrix}
 1 & 0 & 0 \\
 \mathrm{T}_{10}/\delta & \mathrm{T}_{11}/\delta & 0 \\
 \mathrm{T}_{20}/\delta & \mathrm{T}_{21}/\delta & 1
 \end{pmatrix},
 \]
 where $\delta$, $\mathrm{T}_{ij}$ are functions of $u_\lambda$, $u_t$ given by
 \begin{gather*}
 \delta = t (t - 1) (t - \lambda) u_\lambda (u_\lambda - 1) (u_\lambda - \lambda), \\
 \mathrm{T}_{10} = - 2 \nu (t \lambda u_t - t \lambda u_\lambda - t u_t u_\lambda + \lambda u_t u_\lambda - \lambda u_t + t u_\lambda) (t \lambda - t u_\lambda + \lambda u_\lambda - \lambda), \\
 \mathrm{T}_{20} = \nu t (t - 1) \big(2 \lambda^2 u_t - 2 t \lambda u_\lambda + t u_\lambda^2 - \lambda u_\lambda^2 + t \lambda - \lambda^2 - 2 \lambda u_t + 2 \lambda u_\lambda\big), \\
 \mathrm{T}_{11} =- (t \lambda u_t - t \lambda u_\lambda - t u_t u_\lambda + \lambda u_t u_\lambda - \lambda u_t + t u_\lambda)^2, \\
 \mathrm{T}_{21} = - t (t - 1) \big({-}\lambda^2 u_t^2 + 2 t \lambda u_t u_\lambda - t \lambda u_\lambda^2 - t u_t u_\lambda^2 + \lambda u_t u_\lambda^2 - t \lambda u_t \\
 \phantom{ \mathrm{T}_{21} =} + \lambda^2 u_t + \lambda u_t^2 - 2 \lambda u_t u_\lambda + t u_\lambda^2\big).
 \end{gather*}
 Note that the factors that appear in $\mathrm{T}_{10}$ are $P_\Pi$ and $P_\Lambda$ $($defined in Section~{\rm \ref{sec:involution-exceptional-curves})}, and that $\mathrm{T}_{11} = - P_\Pi^2$.
\end{Proposition}

The above proposition is just the result of the computations mentioned earlier.
Looking at the expression for $\mathrm{T}_\psi$, the top row and the last column stand out.
Recall that a connection $\nabla = \nabla_0 +c_1\Theta_1 + c_2\Theta_2$ is represented by the vector $(1, c_1, c_2)$, and so $\psi(\nabla)$ is the connection $\nabla_0 +c'_1\Theta_1 + c'_2\Theta_2$ represented by the vector
\[
 \begin{pmatrix} 1 \\ c'_1 \\ c'_2 \end{pmatrix} = \mathrm{T}_\psi \begin{pmatrix} 1 \\ c_1 \\ c_2 \end{pmatrix}.
\]
The first row of $\mathrm{T}_\psi$ is $\begin{pmatrix}1&0&0\end{pmatrix}$ so that multiplication of $\mathrm{T}_\psi$ with a vector of the form $(1, \ast, \ast)$, will return another vector of the same form.
The third column of $\mathrm{T}_\psi$ is exactly $(0,0,1)$ since the Higgs bundle $\Theta_2(u_\lambda, u_t)$ is equivariant with respect to $(\Psi,\psi)$. This in turn is a consequence of the fact that $u_\lambda$ is unaffected by $\psi$ (cf.~Proposition~\ref{prop:u_t-bar}).

\begin{Remark}\label{T-inverse}
 The fact that $\Psi$ is an involution can be translated to the identity $\mathrm{T}_\psi(u_\lambda, u_t) = \mathrm{T}_\psi(u_\lambda, \overline{u}_t)^{-1}$, which can be easily verified from the above expressions.
\end{Remark}

The equivariant basis is defined by the conditions:
\[
 \nabla_0^\psi(u) = \tfrac{1}{2} (\nabla_0(u) + \Psi\nabla_0(\psi u)), \qquad
 \Theta_i^\psi(u) = \tfrac{1}{2} (\Theta_i(u) + \Psi\Theta_i(\psi u)), \qquad i=1,2.
\]
From these, we deduced that the matrix
$\frac{1}{2} (\mathrm{Id} + \mathrm{T}_\psi(u_\lambda, \overline{u}_t) )$
dictates the base change from the equivariant basis to the original basis (it tells us how $\nabla_0^\psi$, $\Theta_i^\psi$ are written in terms of~$\nabla_0$,~$\Theta_i$). The base change that we seek to define is its inverse. Namely
\[
 \mathrm{B}(u_\lambda, u_t) = 2\left(\mathrm{Id} + \mathrm{T}_\psi(u_\lambda, \overline{u}_t)\right)^{-1} .
\]

\begin{Proposition}\label{prop:formulas-B}
 The explicit expressions for $\mathrm{B}$ are as follows:
 \[
 \mathrm{B} =
 \begin{pmatrix}
 1 & 0 & 0 \\
 \mathrm{B}_{10}/\beta & \mathrm{B}_{11}/\gamma\beta & 0 \\
 \mathrm{B}_{20}/\alpha\beta & \mathrm{B}_{21}/\gamma\beta & 1 \\
 \end{pmatrix},
 \]
 where the numerators are
 \begin{gather*}
 \mathrm{B}_{10} = 2 \nu (t \lambda u_t - t \lambda u_\lambda - t u_t u_\lambda + \lambda u_t u_\lambda - \lambda u_t + t u_\lambda), \\
 \mathrm{B}_{20} = - \nu \big({-}2 t \lambda^2 u_t^2 u_\lambda + 3 t \lambda u_t^2 u_\lambda^2 + 2 t^2 \lambda u_\lambda^3 - 2 t \lambda u_t u_\lambda^3 - t u_t^2 u_\lambda^3 + \lambda u_t^2 u_\lambda^3 - t^2 u_\lambda^4 + t \lambda^2 u_t^2 \\
\hphantom{\mathrm{B}_{20} =}{} + 2 t \lambda^2 u_t u_\lambda - t \lambda u_t^2 u_\lambda + \lambda^2 u_t^2 u_\lambda - 3 t^2 \lambda u_\lambda^2 - 3 \lambda u_t^2 u_\lambda^2 + t^2 u_\lambda^3 - t \lambda u_\lambda^3 + 2 t u_t u_\lambda^3 \\
\hphantom{\mathrm{B}_{20} =}{} + t u_\lambda^4 - \lambda^2 u_t^2 + t^2 \lambda u_\lambda - t \lambda^2 u_\lambda - 2 t \lambda u_t u_\lambda + 2 \lambda u_t^2 u_\lambda + 3 t \lambda u_\lambda^2 - 2 t u_\lambda^3\big), \\
 \mathrm{B}_{11} = 2 (t \lambda u_t - t \lambda u_\lambda - t u_t u_\lambda + \lambda u_t u_\lambda - \lambda u_t + t u_\lambda)^2, \\
 \mathrm{B}_{21} = - t (t - 1) \big(\lambda^2 u_t^2 - 2 t \lambda u_t u_\lambda + t \lambda u_\lambda^2 + t u_t u_\lambda^2 - \lambda u_t u_\lambda^2 + t \lambda u_t - \lambda^2 u_t \\
\hphantom{\mathrm{B}_{21} =}{} - \lambda u_t^2 + 2 \lambda u_t u_\lambda - t u_\lambda^2\big),
 \end{gather*}
 and the denominators are given by
 \begin{gather*}
 \alpha = 2 u_\lambda (u_\lambda - 1) (u_\lambda - \lambda), \\
 \beta = t \lambda u_t^2 - 2 t \lambda u_t u_\lambda - t u_t^2 u_\lambda + \lambda u_t^2 u_\lambda + t^2 u_\lambda^2 - \lambda u_t^2 - t^2 u_\lambda + t \lambda u_\lambda + 2 t u_t u_\lambda - t u_\lambda^2, \\
 \gamma = \lambda u_\lambda - t u_\lambda + \lambda t - \lambda.
 \end{gather*}
\end{Proposition}

This expressions were also obtained by a direct computation. Comparing with the polynomials introduced in Section~\ref{sec:involution-exceptional-curves}, we have $\beta=P_\Sigma$, and $\gamma=P_\Lambda$. The polynomial $P_\Pi$ appears again in $\mathrm{B}_{10}$ and $\mathrm{B}_{11}$.

\subsection{The canonical basis} \label{subsec:canonical-basis}

As explained in Section~\ref{sec:double-family}, in order to have good coordinates on the quotient $\mathcal{U} / \Psi $ we need to replace the basis $(\nabla_0, \Theta_1, \Theta_2)$ by one that is equivariant with respect to $(\Psi, \psi)$.

There is a canonical identification between the spaces $\Higgs^{\bfmu}(C,T)$ and the cotangent bundle $T^*\Bun^{\bfmu}(C,T) \stackrel{\sim}{\dasharrow} T^*\big(\P^1_z \times \P^1_w\big)$, so we seek for the parabolic Higgs bundles $\bfTheta_z$, $\bfTheta_w$ that correspond to ${\rm d}z$, ${\rm d}w$ under this identification.
In fact, we know that the elements $\Theta_1$, $\Theta_2$ in the original basis correspond to ${\rm d}u_t$, ${\rm d}u_\lambda$ under the analogous identification.
However, the identifications between the moduli spaces of parabolic Higgs bundles and the corresponding cotangent bundles do not commute with the map $\Phi$.

Let us make the above paragraph more precise.
On one hand, given a stable parabolic bundle $(E, \bar\ell)$ on the pair $(X,D)$, a parabolic Higgs field is given by some element $\Theta\in H^0\big(\mathrm{End}(E,\bar\ell)\otimes\Omega_X^1(D)\big)$.
On the other hand, we identify the tangent space $T_{(E,\bar\ell)}\Bun^{\bar\mu}(X,D)$ to $H^1(\mathrm{End}(E,\bar\ell))$.
Indeed, a small deformation of $(E,\bar\ell)$ is given by some element $(g_{ij})\in H^1(\mathrm{Aut}(E,\bar\ell))$
(we compare local isomorphisms with $(E,\bar\ell)$ on overlapping open sets on the curve).
Therefore, an infinitesimal deformation of $(E,\bar\ell)$ is given by some element $(\eta_{ij})\in H^1(\mathrm{End}(E,\bar\ell))$.

Consider the following pairing:
\begin{equation} \label{eq:pairing}
\begin{matrix}
 H^0 \! \left(\mathrm{End}(E,\bar\ell)\otimes\Omega_X^1(D)\right) \times H^1 \! \left(\mathrm{End}(E,\bar\ell)\right) & \longrightarrow & H^1\big(\Omega_X^1\big) & \stackrel{\simeq}{\longrightarrow} & \C ,
\end{matrix}
\end{equation}
where the first arrow maps $(\Theta,(\eta_{ij}))$ to the cocycle $(\omega_{ij}) = \tr(\Theta\cdot\eta_{i,j})\in H^1\big(\Omega_X^1\big)$; while the isomorphism on the right $H^1\big(\Omega_X^1\big)\simeq\mathbb C$ associates to a cocycle $(\omega_{ij})$ the sum of its residues.

We remark two things.
First, that the Higgs fields in question represent endomorphisms in the Lie algebra $\mathfrak{sl}_2(\C)$, which is self-dual via the Killing form $(A,B) \mapsto \tr(AB)$.
Second, that $\tr(\Theta\cdot\eta_{i,j})$ has no poles, and so it belongs indeed to $ H^1\big(\Omega_X^1\big)$. The latter claims follows from the fact that, around each point $p$ of the divisor~$D$, $\Theta$ has a simple pole and~$(\eta_{ij})$ is holomorphic; thus the product has a pole of order at most one with residue given by $(\Res_p\Theta \cdot \eta_{ij}(p))$.
Invariance of the parabolic direction~$\ell_p$ translates to the fact that $\Res_p\Theta$ is nilpotent with null space $\ell_p$, and the constant part $(\eta_{ij})(p)$ has $\ell_p$ as an invariant space. The product of such matrices is itself nilpotent, and so the trace has no pole.

Note that these imply that the pairing (\ref{eq:pairing}) is in fact a~version of Serre duality and thus a~perfect pairing.
This is in general how the space $\Higgs^{\bar\mu}(X,D)$ is identified with $T^*\Bun^{\bar\mu}(X,D)$. This identification gives rise to a rational map (it is only defined over stable bundles) which we denote $I_X \colon \Higgs^{\bar\mu}(X,D) \dasharrow T^*\Bun^{\bar\mu}(X,D)$.

Let us now show that this construction does not commute with a base change given by a~ramified cover $Y\to X$.
When we pullback $\Theta$ and $(\eta_{i,j})$ by a finite cover $Y\to X$, the duality pairing on $Y$ associates to this pair the pullback of $(\omega_{ij})$.
However, on $Y$ the sum of residues is multiplied by the degree of the cover.
This is why the identification between $\Higgs^{\bar\mu}(X,D)$ and the contangent bundle $T^*\Bun^{\bar\mu}(X,D)$ commutes with ramified covers only up to a constant factor (given by the degree of the cover).

In order to get a commutative diagram, we introduce a scalar automorphism of $T^*\!\Bun^{\bar\mu}\!\big(\P^1{,}D\big)$ that multiplies each fiber by a factor of 2. The following diagram commutes:
\[\begin{tikzcd}
 \Higgs^{\bar\mu}\big(\P^1,D\big) \arrow[r, "\Phi"] \arrow[d, swap, "I_{\P^1}"]
 &
 \Higgs^{\bfmu}(C,T) \arrow[dd, "I_C"]
 \\
 T^*\Bun^{\bar\mu}\big(\P^1,D\big) \arrow[d, swap, "\times 2"]
 &
 \\
 T^*\Bun^{\bar\mu}\big(\P^1,D\big) \arrow[r, swap, "\hat{\phi}"]
 &
 T^*\Bun^{\bfmu}(C,T) .
\end{tikzcd}\]

We seek for $\bfTheta_z = I_C^{-1}({\rm d}z)$, and $\bfTheta_w = I_C^{-1}({\rm d}w)$.
At the level of cotangent bundles we know that, fiberwise, $\hat{\phi}$ acts as the inverse of the differential of $\phi$, thus
\begin{equation*} 
 \hat{\phi}^{-1}({\rm d}z) = \frac{\partial z}{\partial u_t}{\rm d}u_t + \frac{\partial z}{\partial u_\lambda}{\rm d}u_\lambda, \qquad\text{and}\qquad \hat{\phi}^{-1}({\rm d}w) = \frac{\partial w}{\partial u_t}{\rm d}u_t + \frac{\partial w}{\partial u_\lambda}{\rm d}u_\lambda .
\end{equation*}
Following the commutative diagram, we conclude that $\bfTheta_z$, $\bfTheta_w$ are the images under $\Phi$ of
\begin{equation} \label{eq:pullback-Thetaz-Thetaw}
 \Theta_z = \frac{1}{2}\left( \frac{\partial z}{\partial u_t} \Theta_1 + \frac{\partial z}{\partial u_\lambda} \Theta_2 \right),
 \qquad
 \Theta_w = \frac{1}{2}\left(\frac{\partial w}{\partial u_t} \Theta_1 + \frac{\partial w}{\partial u_\lambda} \Theta_2 \right) ,
\end{equation}
respectively.

\begin{Definition}\label{def:canonical-basis}
 We define Higgs bundles $\Theta_z,\Theta_w \in \Higgs^{\bar\mu}\big(\P^1,D\big)$ by the formulas in (\ref{eq:pullback-Thetaz-Thetaw}).
 We denote their images under $\Phi$ by $\bfTheta_z, \bfTheta_w\in \Higgs^{\bfmu}(C,T)$. We will refer to the triples $\big(\nabla_0^\psi, \Theta_z, \Theta_w\big)$ and $\big(\bfnabla_0^\psi, \bfTheta_z, \bfTheta_w\big)$ as the \emph{canonical bases} for $\mathcal{U}$ and $\mathcal{U}_C$, respectively.
\end{Definition}

We remark that the canonical form on $T^*\Bun^{\bfmu}(C,T)$, and so also on $\Higgs^{\bfmu}(C,T)$, is given by ${\rm d}\kappa_1\wedge {\rm d}z + {\rm d}\kappa_2\wedge {\rm d}w$.
The reduction
\begin{equation}\label{eq:reduction}
 \bfnabla_0^\psi + \kappa_1\bfTheta_z + \kappa_2\bfTheta_w \longmapsto \kappa_1\bfTheta_z + \kappa_2\bfTheta_w
\end{equation}
defines a map $\Con^{\bfmu}_{\bfnu}(C,T) \dasharrow \Higgs^{\bfmu}(C,T)$, which in coordinates is nothing but the identity map.
The symplectic form on $\Con^{\bfmu}_{\bfnu}(C,T)$ will be thus given by
\[
 \omega_C = {\rm d}\kappa_1\wedge {\rm d}z + {\rm d}\kappa_2\wedge {\rm d}w,
\]
provided that the map in (\ref{eq:reduction}) is symplectic.
This will be established at the end of the present section, cf.~Lemma~\ref{lemma:nabla-Lagrangian}.

\begin{Proposition}\label{prop:formula-canonical-basis}
 The parabolic Higgs bundles $\Theta_z$, $\Theta_w$ in the canonical basis are equivariant with respect to the pair $(\Psi,\psi)$, i.e.,~$\Psi\Theta_i(u) = \Theta_i(\psi u)$, for $i=z,w$, and $u \in \P^1_{u_\lambda} \times \P^1_{u_t}$. Moreover, with respect to the equivariant basis in Definition~{\rm \ref{def:equivariant-basis}}, they can be expressed as
 \begin{equation}\label{eq:formula-canonical-basis}
 \Theta_z = \frac{(z-\lambda)^2}{2\lambda(1-\lambda)} \Theta_2^\psi,
 \qquad
 \Theta_w = \frac{z-\lambda}{z-t} \Theta_1^\psi + \frac{(wt-w\lambda-t\lambda+\lambda) (z-\lambda)}{2(z-t)(\lambda-1)\lambda} \Theta_2^\psi,
 \end{equation}
 where $z$, $w$ are the functions given in~\eqref{eq:formula-Nestor}.
\end{Proposition}

\begin{proof}
 We have three bases to describe the universal family $\mathcal{U}$, and with them the following base-change matrices:
 \begin{enumerate}\itemsep=0pt
 \item[1)] $\mathrm{B}$ going from the original basis to the equivariant one,
 \item[2)] $\mathrm{J}$ going from the canonical basis to the original one,
 \item[3)] $\mathrm{C}$ going from the canonical basis to the equivariant one.
 \end{enumerate}
 If formula (\ref{eq:formula-canonical-basis}) is true, it would define the change of coordinates $\mathrm{C}$. Thus we shall compute the matrix $\mathrm{C}=\mathrm{BJ}$ and from it deduce (\ref{eq:formula-canonical-basis}).
 The matrix $\mathrm{B}$ is given in Proposition~\ref{prop:formulas-B}. The matrix $\mathrm{J}$ is given by putting together the first column of $\mathrm{B}^{-1} = \frac{1}{2}\left(\mathrm{Id} + \mathrm{T}_\psi(u_\lambda, \overline{u}_t)\right)$ (which defines the connection $\nabla_0^\psi$), and the (transposed) Jacobian matrix $\frac{\partial (z,w)}{\partial (u_t, u_\lambda)}$ multiplied by a factor of $\frac{1}{2}$, which defines the parabolic Higgs bundles $\Theta_z$, $\Theta_w$. More precisely,
 \begin{equation} \label{eq:J}
 \mathrm{J} =
 \begin{pmatrix}
 1 & 0 & 0 \\
 c_1^0 & \dfrac{1}{2}\dfrac{\partial z}{\partial u_t} & \dfrac{1}{2}\dfrac{\partial w}{\partial u_t}\vspace{1mm}\\
 c_2^0 & \dfrac{1}{2}\dfrac{\partial z}{\partial u_\lambda} & \dfrac{1}{2}\dfrac{\partial w}{\partial u_\lambda}
 \end{pmatrix},
 \end{equation}
 where $\nabla_0^\psi = \nabla_0 + c_1^0\Theta_1 + c_2^0\Theta_2$.
 The matrix $\mathrm{C} = \mathrm{BJ}$ is the one we are at the moment interested in. Having explicit expressions for both $\mathrm{B}$ and~$\mathrm{J}$, we can compute the explicit expressions for~$\mathrm{C}$. It is straightforward to check that the entries of this matrix are invariant under~$\psi$, and in fact can be rewritten in terms of the $\psi$-invariant functions~$z$,~$w$, see~(\ref{eq:formula-Nestor}).
 These expressions coincide with the coefficients given in~(\ref{eq:formula-canonical-basis}).
 In particular, this proves that the parabolic Higgs bundles~$\Theta_z$ and~$\Theta_w$ are equivariant with respect to~$\Psi$,~$\psi$, since they are expressed as a linear combination of the equivariant bundles $\Theta_1^\psi$, $\Theta_2^\psi$, with $\psi$-invariant coefficients.
\end{proof}

For future reference, we provide here an explicit expression of the matrix $\mathrm{C}^{-1}$,
\[
 \mathrm{C}^{-1} =
 \begin{pmatrix}
 1 & 0 & 0 \\
 0 & \mathrm{K}_{11} & \mathrm{K}_{12}\\
 0 & \mathrm{K}_{21} & 0
 \end{pmatrix},
\]
where
\begin{equation}\label{eq:entries-C-inv}
 \mathrm{K}_{11} = \frac{wt-w\lambda-t\lambda+\lambda}{(z-\lambda)^2}, \qquad
 \mathrm{K}_{21} = \frac{z-t}{z-\lambda}, \qquad
 \mathrm{K}_{12} = \frac{2\lambda(1-\lambda)}{(z-\lambda)^2}.
\end{equation}

We arrive finally to a complete description of the map $\Phi$ in coordinates:

\begin{Proposition}\label{prop:final-expression-Phi}
 The map
 \[
 \Phi \colon \ \P^1_{u_\lambda} \times \P^1_{u_t} \times \C^2_{(c_1,c_2)} \dashrightarrow \P^1_z\times\P^1_w \times \C^2_{(\kappa_1,\kappa_2)}
 \]
 decomposes as $\Phi = (\phi, \mathrm{T}_\phi)$ where $\phi \colon \P^1_{u_\lambda} \times \P^1_{u_t} \dashrightarrow \P^1_z\times\P^1_w$ is given by \eqref{eq:formula-Nestor}, and $\mathrm{T}_\phi \colon \P^1_{u_\lambda} \times \P^1_{u_t} \times \C^2_{(c_1,c_2)} \dashrightarrow \C^2_{(\kappa_1,\kappa_2)}$ is defined by
 \[
 \begin{pmatrix} 1 \\ \kappa_1 \\ \kappa_2 \end{pmatrix}
 = \mathrm{J}(u_\lambda, u_t)^{-1} \begin{pmatrix} 1 \\ c_1 \\ c_2 \end{pmatrix} ,
 \]
 where the matrix $\mathrm{J}(u_\lambda, u_t)$ is defined in \eqref{eq:J} and the inverse $\mathrm{J}^{-1}$ is given by the product $\mathrm{J}^{-1} = \mathrm{C}^{-1}\mathrm{B}$. The explicit expression of the matrix $\mathrm{B}$ appears in Proposition~{\rm \ref{prop:formulas-B}}, and the entries of $\mathrm{C}^{-1}$ are given in \eqref{eq:entries-C-inv}.
\end{Proposition}

Having explicit formulas for the map $\Phi$, we can finally prove the following lemma.

\begin{Lemma} \label{lemma:nabla-Lagrangian}
 The section $\bfnabla_0^\psi \colon \Bun^{\bfmu}(C,T) \dasharrow \Con^{\bfmu}_{\bfnu}(C,T)$ is Lagrangian, and the reduction of $\Con^{\bfmu}_{\bfnu}(C,T)$ to $\Higgs^{\bfmu}(C,T)$ using $\bfnabla_0^\psi$ is a symplectic identification.
\end{Lemma}

\begin{proof} At this point, we don't know the exact expression of $\omega_C$, we only know from Corollary~\ref{cor:PhiIsSymplectic} that $\Phi^*\omega_C = 2\omega_{\P^1}$. We are now going to prove that $\omega_C = {\rm d}\kappa_1\wedge {\rm d}z + {\rm d}\kappa_2\wedge {\rm d}w$.
 Let us define $\delta = \omega_C - ({\rm d}\kappa_1\wedge {\rm d}z + {\rm d}\kappa_2\wedge {\rm d}w)$. Using the explicit formulas that we have computed above for the map $\Phi$, we can readily verify that $\Phi^*(\kappa_1\wedge {\rm d}z + {\rm d}\kappa_2\wedge {\rm d}w) = 2\omega_{\P^2}$. This means that $\Phi^*\delta = \Phi^*\omega_C - \Phi^*({\rm d}\kappa_1\wedge {\rm d}z + {\rm d}\kappa_2\wedge {\rm d}w) = 0$.
 Note that around a generic point of $\Con^{\bar\mu}_{\bar\nu}\big(\P^1,D\big)$, the map $\Phi$ is a local biholomorphism. This implies that $\Phi^*\delta$ can only vanish identically if $\delta\equiv0$.
 This establishes the fact that $\omega_C = {\rm d}\kappa_1\wedge {\rm d}z + {\rm d}\kappa_2\wedge {\rm d}w$.
 Since the same expression defines the symplectic form on $\Higgs^{\bfmu}(C,T)$, we conclude that the identification between these two spaces is symplectic.

 From the formula for $\omega_C = {\rm d}\kappa_1\wedge {\rm d}z + {\rm d}\kappa_2\wedge {\rm d}w$ it follows immediately that the section $\bfnabla_0^\psi$, whose image is given in coordinates by $\{(z,w,\kappa_1,\kappa_2) \,\vert\, \kappa_1=\kappa_2=0 \}$, is Lagrangian.
\end{proof}

\subsection[Final details on the map Phi]{Final details on the map $\boldsymbol{\Phi}$}\label{sec:final-steps}

We conclude this section by filling in details to establish the second part of Theorem~\ref{thm:Bun-C}.

\begin{proof}[Proof of Theorem~\ref{thm:Bun-C}]
 The coordinate system $\Con^{\bfmu}_{\bfnu}(C,T) \dashrightarrow \P^1_z\times\P^1_w \times \C^2_{(\kappa_1,\kappa_2)}$ and the map $\Phi$ are defined in terms of the canonical family $\bfnabla_0^\psi$, $\bfTheta_z$, $\bfTheta_w$.
 In order to prove the theorem it is enough to establish that, when the underlying parabolic bundle belongs to $\Bun^{\bfmu}(C,T)^0$, the canonical family $\bfnabla_0^\psi$, $\bfTheta_z$, $\bfTheta_w$ is well defined, and that the parabolic Higgs bundles $\bfTheta_z$, $\bfTheta_w$ are linearly independent.
 We will prove these properties for the canonical family $\nabla_0^\psi$, $\Theta_z$, $\Theta_w$.
 According to Proposition~\ref{prop:final-expression-Phi}, we need to focus on the base change $\mathrm{J}^{-1} = \mathrm{C}^{-1}\mathrm{B}$ which converts the original family $\nabla_0$, $\Theta_t$, $\Theta_\lambda$ used in \cite{LoraySaito2015} into the canonical one.

 First, let us note the following. The map $\S \dashrightarrow \P^1_{u_\lambda} \times \P^1_{u_t}$ is obtained by contracting four of the sixteen $(-1)$-curves in $\S$. Therefore, the space of generic bundles $\Bun^{\bar\mu}\big(\P^1,D\big)^0$ is in one-to-one correspondence with an open subset of $\P^1_{u_\lambda} \times \P^1_{u_t}$ (the complement of the remaining twelve rational curves we did not contract). This means that every generic bundle is represented in $\P^1_{u_\lambda} \times \P^1_{u_t}$. We remark that our genericity assumptions also exclude the lines $u_\lambda = \infty$, $u_t = \infty$ (these are mapped by $\phi$ to $z=\lambda$ and $w=\infty$, respectively). Therefore the affine coordinates we are using in $\P^1_{u_\lambda} \times \P^1_{u_t}$ suffice for all computations.

 Next, we recall that the map $\phi \colon \S \to \P^1_z\times\P^1_w$ is in fact everywhere well defined (regular). Our formula for~$w$ in (\ref{eq:formula-Nestor}) has points of indeterminacy precisely at the four points obtained by contracting $(-1)$-curves in~$\S$, but these points are not in $\Bun^{\bar\mu}\big(\P^1,D\big)^0$.

 Finally, we analyze the poles of the entries of the matrices $\mathrm{B}$ and $\mathrm{C}^{-1}$, as well as the zeros of their determinants.
 Let us start with~$\mathrm{B}$.
 The polynomial $\alpha$ in Proposition~\ref{prop:formulas-B} vanishes at the lines $u_\lambda = 0$, $u_\lambda=1$, $u_\lambda=\lambda$. These are among the rational curves excluded by $\Bun^{\bar\mu}\big(\P^1,D\big)^0$. They are mapped by $\phi$ to the curves $z=1$, $z=0$, $z=\infty$, respectively. The polynomial~$\beta$ is exactly the polynomial $P_\Sigma$ defining the ramification locus $\Sigma$ in~(\ref{eq:Sigma}). Lastly, the polynomial~$\gamma$ coincides with $P_\Lambda$ which defines the special line $\Lambda$ in~(\ref{eq:Lambda}).
 Since the matrix is triangular we evidently have $\det(\mathrm{B}) = \mathrm{B}_{11}/\gamma\beta$. The denominators we have already discussed, and the numerator is precisely $2P_\Pi^2$, where $P_\Pi$ was given in (\ref{eq:Pi}). Since this is another rational curve excluded by $\Bun^{\bar\mu}\big(\P^1,D\big)^0$, we conclude that $\det(\mathrm{B})$ is never zero for a generic parabolic bundle, and it has a simple pole over $\Lambda$.
 From the equations~(\ref{eq:entries-C-inv}) it is straightforward that the poles of~$\mathrm{C}^{-1}$ are given by $z=\lambda$ and the determinant only vanishes at $z=t$.
 We remark that the line~$\Lambda$ is equivalently defined by the zeroes of~$P_\Lambda$ in~(\ref{eq:Lambda}) or by the condition $z=t$. Therefore the simple pole $\det(\mathrm{B})$ is cancelled by the simple zero of $\det (\mathrm{C}^{-1})$.
 A straight forward computation shows that indeed the entries of~$\mathrm{J}^{-1}$ have no poles over $\Lambda$.
 Thus the coordinate change $\mathrm{J}^{-1}$ is well defined over~$\Lambda$.

 We conclude that over $\Bun^{\bar\mu}\big(\P^1,D\big)^0 \setminus \Sigma$ the canonical family $\nabla_0^\psi$, $\Theta_z$, $\Theta_w$ is well defined, and~$\Theta_z$,~$\Theta_w$ are linearly independent.
 As a consequence, the canonical family $\bfnabla_0^\psi$, $\bfTheta_z$, $\bfTheta_w$ enjoys the same properties over $\Bun^{\bfmu}(C,T)^0$.
\end{proof}

\begin{Remark}
 All our constructions are well defined for parabolic bundles in $\Bun^{\bar\mu}\big(\P^1,D\big)^0$, except that the base-change matrix $\mathrm{J}^{-1}$ has poles over $\Sigma$.
 This is because the Jacobian matrix $\frac{\partial (z,w)}{\partial (u_t, u_\lambda)}$ drops rank over $\Sigma$. Thus, the canonical basis introduced in Definition~\ref{def:canonical-basis} is not a true basis: it spans only a one-dimensional linear space when the underlying parabolic bundle belongs to $\Sigma$ (see also Remark~\ref{rmk:eq-family-over-Sigma}).
 Finally, this means that the coordinate system $\Con^{\bfmu}_{\bfnu}(C,T) \dashrightarrow \P^1_z\times\P^1_w \times \C^2_{(\kappa_1,\kappa_2)}$ cannot be used to describe connections with such underlying parabolic bundles.
\end{Remark}

\subsection{The singular locus on the space of connections}
 \label{subsec:singular-locus}

We have described $\Con^{\bfmu}_{\bfnu}(C,T)$ as a quotient of $\Con^{\bar\mu}_{\bar\nu}\big(\P^1,D\big)$ by the involution $\Psi$. As announced in Section~\ref{sec:main-body}, the set of fixed points of $\Psi$ is a codimension~2 subvariety, which causes the quotient to be singular at the image of the fixed-point locus.
We shall first characterize the fixed-point set in $\Con^{\bar\mu}_{\bar\nu}\big(\P^1,D\big)$, and then describe the singularities of $\Con^{\bfmu}_{\bfnu}(C,T)$ locally.

\begin{proof}[Proof of Theorem~\ref{thm:fixed-points-psi}]
 We analyze the fixed points of $\Psi $ in $\mathcal{U}$. This is equivalent to the fixed points of $\bfPsi$ in $\mathcal{U}_C$. Any such fixed point must be a connection defined over a parabolic bundle that is fixed by $\psi $, namely a parabolic bundle $(u_\lambda, u_t) \in \Sigma$. Those connections that are fixed by $\Psi $ are the solutions to the linear inhomogeneous system
 \begin{equation}\label{eq:linear-system}
 \left( \mathrm{T}_\psi(u_\lambda, \overline{u}_t) - \mathrm{Id} \right) \begin{pmatrix} 1 \\ c_1 \\ c_2 \end{pmatrix}
 = 0.
 \end{equation}
 Because of the particular shape of $\mathrm{T}_\psi$ (cf.~Proposition~\ref{prop:formulas-T}), the last column and the first row of the matrix on the left-hand side of~(\ref{eq:linear-system}) are zero. Therefore we are left with two inhomogeneous equations on~$c_1$ only (and so~$c_2$ is free to take any value).
 Recall that the central entry in $\mathrm{T}_\psi(u_\lambda,u_t)$ is given by $\mathrm{T}_{11}/\delta$.
 We can easily verify that $(\mathrm{T}_{11}/\delta)\vert_\Sigma \equiv -1$. Indeed, a quick computation shows that
 \[
 \frac{\mathrm{T}_{11}}{\delta} = -1 - \frac{P_\Sigma P_\Lambda}{\delta},
 \]
 where $P_\Pi$, $P_\Sigma$, $P_\Lambda$ are defined in (\ref{eq:Pi})--(\ref{eq:Lambda}).

 \begin{Remark}\label{rmk:eigenvalue-1}
 The above equations shows that $\Sigma \cup \Lambda$ is precisely the locus where the involution~$\Psi $ acts as $\Theta_1 \mapsto -\Theta_1 + k\Theta_2$, where $k$ is a scalar (and recall that $\Theta_2$ is always equivariant). Because of this, the matrix $\frac{1}{2}(\mathrm{Id} + \mathrm{T}_\psi(u_\lambda, \overline{u}_t))$ used to define the equivariant basis in Definition~\ref{def:canonical-basis} drops rank and we are unable to define such basis above these parabolic bundles.
 \end{Remark}

 The determinant of the bottom-left $2\times2$ minor of the matrix in (\ref{eq:linear-system}) vanishes at $\Sigma$, and the middle row imposes an equation
 \[
 \mathrm{T}_{10}/\delta -2c_1 = 0.
 \]
 We conclude that the linear space of solutions is one-dimensional and defined by
 \begin{equation}\label{eq:c1-val}
 c_1 = \mathrm{T}_{10}/2\delta = \nu \frac{P_\Pi P_\Lambda}{\delta}.
 \end{equation}

 In Section~\ref{sec:geom-pic}, parabolic bundles in $\Sigma$ are characterized as those for which the parabolic direction belongs to the curve $S_\Sigma \subset \P\big(\mathcal{O}_{\P^1} \oplus \mathcal{O}_{\P^1}\big)$.
 A quick analysis shows that connections that moreover satisfy~(\ref{eq:c1-val}) are precisely those for which the second eigenvector of its residue over $t$ (corresponding to the eigenvalue~$-\nu$) also belongs to the curve $S_\Sigma$. Further analysis reveals that in this case the Riccati foliation is indeed tangent to $S_\Sigma$. After performing the transformation $\Phi$ we recover a connection on a bundle of the form $\mathbf{L} \oplus \mathbf{L}^{-1}$ for which each summand is invariant. Invariance by the elliptic involution (which exchanges $\mathbf{L}$ and $\mathbf{L}^{-1}$) implies that the connection must be of the form
 \[
 \big(\bfnabla, \mathbf{E}, \bar\ell\big) =
 \big(\eta, \mathbf{L}, \ell\big) \oplus \iota^* \big(\eta, \mathbf{L}, \ell\big).\tag*{\qed}
 \]\renewcommand{\qed}{}
\end{proof}

\begin{Remark}\label{rmk:eq-family-over-Sigma}
 We have shown that the equivariant basis fails to describe the universal family $\mathcal{U}$ whenever the underlying parabolic bundle belongs to $\Sigma$. This is in fact true for \emph{any} equivariant basis, including the canonical one. Indeed, since $\Sigma$ defines the fixed points of $\psi $, we have that any equivariant connection $\nabla$, parametrized by $u\in\Bun^{\bar\mu}\big(\P^1,D\big)$, satisfies
 \[
 \Psi\nabla(u) = \nabla(\psi u) = \nabla(u), \qquad \text{if} \quad u\in\Sigma.
 \]
 Therefore, if $u\in\Sigma$ equivariance implies invariance under $\Psi $. As shown above, for each fixed parabolic bundle $u\in\Sigma$, the space of $\Psi$-invariant connections is one-dimensional.
 Thus any equivariant basis drops rank over $\Sigma$.
\end{Remark}

\begin{proof}[Proof of Theorem~\ref{thm:Sing-locally}]
 Let us choose a generic point $u^0\in\Sigma$ (by generic we mean that $u_\lambda \neq 0,1,\lambda,\infty)$, and $\nabla(u^0)$ a~connection fixed by~$\Psi $. We will describe the space $\Con^{\bfmu}_{\bfnu}(C,T)$ around $\Phi\big(\nabla\big(u^0\big)\big)$.

 Consider the curve $\Sigma \subset \P^1_{u_\lambda} \times \P^1_{u_t}$, which is the set of fixed points of the involution $\psi $. Up to a linear factor, the curve is given by hyperelliptic equation
 \begin{equation*}
 P_\Pi = \sqrt{\delta},
 \end{equation*}
 where $P_\Pi$ and $\delta$ are given in (\ref{eq:Pi}) and Proposition~\ref{prop:formulas-T}, respectively.
 We propose the following \emph{local} change of variables
 \begin{equation}\label{eq:def-U}
 (u_t, u_\lambda) \longmapsto (U, u_\lambda), \qquad\text{where}\quad
 U = \frac{P_\Pi - \sqrt{\delta}}{P_\Pi + \sqrt{\delta}}.
 \end{equation}
 Let us assume that $u^0\mapsto U^0$, and choose a branch of the square root such that $\Sigma$ is given by $U = 0$ around the point $U^0$.
 A straightforward computation shows that in these coordinates the involution $\psi $ is given by $(U, u_\lambda) \mapsto (-U, u_\lambda)$.

 As discussed in Remark~\ref{rmk:eigenvalue-1}, the matrix $\mathrm{T}_\psi$ has a an eigenvalue equal to $-1$ whenever the underlying parabolic bundle belongs to $\Sigma$. Therefore, with respect to a suitably chosen (non-equivariant) basis $\widetilde{\Theta}_1$, $\Theta_2$, the involution $\Psi $ acts as $\widetilde{\Theta}_1 \mapsto -\widetilde{\Theta}_1$, and $\Theta_2$ is unchanged by $\Psi$.

 Let us make precise who the above Higgs bundle $\widetilde{\Theta}_1$ is.
 For this we identify, over an open subset, the moduli space $\Con^{\bar\mu}_{\bar\nu}\big(\P^1,D\big)$ to the space of parabolic Higgs bundles $\Higgs^{\bar\mu}\big(\P^1,D\big)$ using the Lagrangian section $\nabla_0^\psi$ (cf.~Lemma~\ref{lemma:nabla-Lagrangian}).
 Because this section is equivariant with respect to $(\Psi,\psi)$, the identification $\Con^{\bar\mu}_{\bar\nu}\big(\P^1,D\big) \dasharrow \Higgs^{\bar\mu}\big(\P^1,D\big)$ commutes with the action of $\Psi$.
 Now, it follows from the identity $\Phi\circ\Psi = \Phi$, that the map $\Psi$ is symplectic.
 This can also be verified by a direct computation using the formulas in Proposition~\ref{prop:formulas-T}.
 Moreover, the automorphism $\psi$ of $\Bun^{\bar\mu}\big(\P^1,D\big)$ lifts to a canonical symplectic automorphism $\hat{\psi}$ of $T^*\Bun^{\bar\mu}\big(\P^1,D\big)$ which covers $\psi$.
 We can verify from our formulas that, under the identification $\Higgs^{\bar\mu}\big(\P^1,D\big) \dasharrow T^*\Bun^{\bar\mu}\big(\P^1,D\big)$, we have $\hat{\psi} = \Psi$.
 The coordinate change (\ref{eq:def-U}) induces a~symplectic coordinate change
 \[
 (u_t, {\rm d}u_t, u_\lambda, {\rm d}u_\lambda) \mapsto (U, {\rm d}U, u_\lambda, {\rm d}u_\lambda).
 \]
 In these coordinates, $\psi$ is given by the linear action $(U, u_\lambda) \mapsto (-U, u_\lambda)$ and so $\Psi$, which is the lift of $\psi$ to the cotangent bundle, is given by $(U, {\rm d}U, u_\lambda, {\rm d}u_\lambda) \mapsto (-U, -{\rm d}U, u_\lambda, {\rm d}u_\lambda)$.
 Thus the action $\Psi$ has been linearized in these coordinates.

 \begin{Remark} Unwrapping the above arguments we see that the Higgs bundle $\widetilde{\Theta}_1$ we were looking for (which corresponds to ${\rm d}U$) is given explicitly by
 $\widetilde{\Theta}_1 = \frac{\partial U}{\partial u_t} \Theta_1 + \frac{\partial U}{\partial u_\lambda} \Theta_2$. Thus the basis in which the action of $\Psi$ is linear is given by $\nabla_0^\psi, \widetilde{\Theta}_1$, $\Theta_2$.
 \end{Remark}

 Let us rename the coordinates as $(y_1, y_2, y_3, y_4) = (U, {\rm d}U, u_\lambda, {\rm d}u_\lambda)$.
 In this way $\Psi$ is given by $\Psi(y)=(-y_1,-y_2,y_3,y_4)$, and the symplectic structure is given by $\omega={\rm d}y_1\wedge {\rm d}y_2+{\rm d}y_3\wedge {\rm d}y_4$. The map
 \[\mathbb C^4\to\mathbb C^5 ;\qquad y\mapsto x=\big(y_1y_2,y_1^2,y_2^2,y_3,y_4\big)\]
 realizes the quotient $\C^4_{\underline{y}} / \Psi $, whose image is locally the product of a bidisc in variables $(x_3,x_4)$
 by a quadratic conic singularity given in $\C^3$ given by $x_0^2 = x_1x_2$.
 Finally, we remark that the $2$-form $\omega={\rm d}y_1\wedge {\rm d}y_2+{\rm d}y_3\wedge {\rm d}y_4$ becomes $\tilde{\omega} = \frac{{\rm d}x_1\wedge {\rm d}x_2}{4x_0}+{\rm d}x_3\wedge {\rm d}x_4$ in the quotient.
 This completes the proof of the theorem.
\end{proof}

\section{The apparent map in the elliptic case} \label{sec:app}

Consider the projectivization of the trivial bundle $\P\big(\mathcal{O}_{\P^1} \oplus \mathcal{O}_{\P^1}\big)$. In~\cite{LoraySaito2015} the apparent map is defined with respect to the constant horizontal section $\sigma_\infty$, defined by $\zeta=\infty$ with respect to an affine coordinate $\zeta$ on the fiber.
After applying the transformation~$\Psi $, the section $\sigma_\infty$ becomes a section of self-intersection~$+2$, which we denote $\sigma_\psi$.
In general, the tangencies of a connection with $\sigma_\infty$, and with $\sigma_\psi$ occur at different points. This means that the apparent map is not invariant under~$\Psi $.
If we want to use $\Phi$ to push the concept of the apparent map from $\Con^{\bar\mu}_{\bar\nu}\big(\P^1,D\big)$ to $\Con^{\bfmu}_{\bfnu}(C,T)$, we need to redefine the apparent map in such a way that it becomes invariant under $\Psi $.
To do this, we consider both the tangency loci of the connection with $\sigma_\infty$ and with $\sigma_\psi$ simultaneously. This defines an element of $\P^2\times\P^2$. The action of $\Psi $ permutes these factors, so we need to pass to the symmetric product $\SP$.

\begin{Definition}\label{def:app}
 We define the \textit{apparent map} as the unique map $\operatorname{App}_C \colon \Con^{\bfmu}_{\bfnu}(C,T) \dashrightarrow \SP$ that completes the following commutative diagram:
 \[
 \begin{tikzcd}[column sep=large, row sep=large]
 \Con^{\bar\mu}_{\bar\nu}\big(\P^1,D\big) \arrow[rr, dashrightarrow, "\operatorname{App} \times (\operatorname{App}\circ\Psi)"] \arrow[d, swap, "\Phi"]
 &
 &
 \P^2 \times \P^2 \arrow[d, "\operatorname{Sym}"]
 \\
 \Con^{\bfmu}_{\bfnu}(C,T) \arrow[rr, dashrightarrow, swap, "\exists!\operatorname{App}_C"]
 &
 &
 \SP .
 \end{tikzcd}
 \]
\end{Definition}

Note that, unlike the case for $\P^1$, the above map is defined between spaces of the same dimension. Therefore, the map cannot be Lagrangian. It is a generically finite map, but the correspondence is not one-to-one.

In order to study the apparent map $\operatorname{App}_C$ defined above, we begin with connections on $\P^1$ and switch from the birational model $\Con^{\bar\mu}_{\bar\nu}\big(\P^1,D\big) \stackrel{\sim}{\dashrightarrow} \P^1_{u_\lambda} \times \P^1_{u_t} \times \C^2_{(c_1,c_2)}$, to the model $\Con^{\bar\mu}_{\bar\nu}\big(\P^1,D\big) \stackrel{\sim}{\dashrightarrow} \P^2_a \times \P^2_b$ defined by the map (\ref{eq:AppBun-P1}).
This birational model of $\Con^{\bar\mu}_{\bar\nu}\big(\P^1,D\big)$ is studied in detail in \cite{LoraySaito2015}.
The first factor, $\P^2_a$ defines the image of the apparent map. Indeed, the tangencies of a generic connection $\nabla(a,b)$ with the section $\sigma_\infty$ are precisely at the roots of the polynomial $\operatorname{App}_\infty (a,b) = a_2x^2 + a_1x + a_0$. The second factor, $\P^2_b$ defines the underlying parabolic bundle.
Explicit formulas to go from one coordinate system to the other are given in \cite[Section~6]{LoraySaito2015}, and so we omit them here.

Recall from Section~\ref{sec:genericity} that the Del~Pezzo surface $\S$ can be identified with the blow-up of $\P^2_b$ at five points, which we call $D_0$, $D_1$, $D_\lambda$, $D_\infty$, $D_t$. As shown in \cite[Section~6]{FernandezVargas2016}, the involution~$\psi$ is the lift of the de~Jonqui\`{e}res automorphism of $\P^2_b$ preserving the pencil of lines through $D_t$ and the pencil of conics through $D_0$, $D_1$, $D_\lambda$, $D_\infty$.
The following five rational curves in $\P^2_b$ are important for the upcoming discussion: the conic~$\Pi$ through all five points~$D_i$, and the lines~$L_{i t}$ passing through the points $D_i$ and $D_t$, for $i=0,1,\lambda,\infty$.
These become $(-1)$-curves in~$\S$, and so they are excluded from $\Con^{\bar\mu}_{\bar\nu}\big(\P^1,D\big)^0$ (cf.~Definition~\ref{def:generic-bunle-P1}).
Some of these rational curves have already appeared in previous sections. They correspond to the curves defined by~$P_\Pi$ in~(\ref{eq:Pi}), and the lines $u_\lambda = 0,1,\infty$. Only the line $L_{\lambda t}$ is absent from the $\P^1_{u_\lambda} \times \P^1_{u_t}$ model. In fact, the birational isomorphism~$\P^2_b \dashrightarrow \P^1_{u_\lambda} \times \P^1_{u_t}$ is obtained by blowing-up $D_\lambda$ and $D_t$, and contracting the line $L_{\lambda t}$ through them.

\begin{Remark}
 As always, the involution $\Psi$ plays a crucial role in the passage from connections over $\P^1$ to connections over $C$. In these coordinates, the involution acts as $\Psi \colon (a,b) \mapsto (s,\bar{b})$. The action on the parabolic bundles $b \mapsto \bar{b}$ is the de~Jonqui\`{e}res automorphism of $\P^2_b$ discussed above. Below we seek to understand the correspondence $(a,b) \mapsto s$. This is given by a~matrix~$\mathrm{M}_b$, whose entries are parametrized by $b$. This matrix will be the main object of study in this section. Note that since $a$ is given by the apparent map $\operatorname{App}_\infty$, we have that $s = \operatorname{App}_\infty\!\circ\,\Psi$.
\end{Remark}

We now analyze the map $\operatorname{App}_\infty \times (\operatorname{App}_\infty\!\circ\,\Psi) \colon \Con^{\bar\mu}_{\bar\nu}\big(\P^1,D\big) \dashrightarrow \P^2\times\P^2$, which was introduced in Definition~\ref{def:app}.
Recall we have defined $\operatorname{App}_\psi = \operatorname{App}_\infty\!\circ\,\Psi$. Let us introduce homogeneous coordinates $s=[s_0:s_1:s_2]$ on $\P^2$ so that $\operatorname{App}_\infty$ takes values on $\P^2_a$, and $\operatorname{App}_\psi$ takes values on $\P^2_s$. Under these coordinate systems the former map is nothing but{\samepage
\begin{equation}\label{eq:ab-as}
 \operatorname{pr}_1 \times \operatorname{App}_\psi \colon \ \P^2_a\times\P^2_b \dashrightarrow \P^2_a\times\P^2_s ,
\end{equation}
where $\operatorname{pr}_1$ denotes projection onto the first factor $\P^2_a$.}

We remark that, if we fix $b$ a generic bundle, the map $\operatorname{App}_\psi(\_, b) \colon \P^2_a \to \P^2_s$ is holomorphic and invertible. Therefore it defines an element of $\operatorname{PGL}(3,\C)$. It is a straightforward computation to translate the formula for $\operatorname{App}_\psi$ mentioned in Section~\ref{sec:comp-coords} to these new coordinates.

\begin{Proposition}\label{prop:app_psi}
 The map $\operatorname{App}_\psi \colon \P^2_a\times\P^2_b \dashrightarrow \P^2_s$ is defined by a $3\times3$ matrix $\mathrm{M}_b$ via
 \[
 (a,b) \longmapsto \mathrm{M}_b
 \begin{pmatrix}
 a_0 \\ a_1 \\ a_2
 \end{pmatrix}.
 \]
 The entries of $\mathrm{M}_b$ are homogeneous polynomials on $b$ of degree four. Moreover, its determinant only vanishes along the curves $\Pi$ and $L_{i t}$, for $i=0,1,\lambda,\infty$ $($each divisor with multiplicity two$)$.
\end{Proposition}

This proposition is proved by a direct inspection of the explicit formulas for $\mathrm{M}_b$, which we exclude here.

\begin{Remark}\label{rmk:ab-s_well-defined}
 The above proposition implies that the map $\operatorname{App}_\infty \times (\operatorname{App}_\infty\!\circ\,\,\Psi)$ is everywhere well defined on $\Con^{\bar\mu}_{\bar\nu}\big(\P^1,D\big)^0$.
\end{Remark}

That (\ref{eq:ab-as}) is dominant is straightforward: it is a map between irreducible spaces of the same dimension, and we can check at a generic point that the derivative has maximal rank. The fact that the determinant of $\mathrm{M}_b$ factors as a product of rational curves makes it easy to analyze the behavior of this matrix along such divisor (since we can parametrize each branch). Below we describe the kernel of $\mathrm{M}_b$ at every point where the determinant vanishes. These are linear subspaces of $\C^3_a$ which are naturally identified with subsets of $\P^2_a = \P\big(\C^3_a\big)$.
\begin{enumerate}\itemsep=0pt
 \item At a generic point of the conic $\Pi$ the matrix $\mathrm{M}_b$ has rank 2 and its kernel is a fixed point $\tilde D_t \in \P^2_a$. \label{item:1}
 \item At a generic point of the line $L_{i t}$, for $i=0,1,\lambda,\infty$, the matrix $\mathrm{M}_b$ has rank 1 with a~one-dimensional kernel that can be identified with a fixed point $\tilde D_i \in \P^2_a$. \label{item:2}
 \item At the point $D_i$, for $i=0,1,\lambda,\infty$, the matrix $\mathrm{M}_b$ has rank 1 with a two-dimensional kernel which we identify with the line $L_{i t}$ through $\tilde D_t$ and $\tilde D_i$ in $\P^2_a$.
 \item At the point $D_t$ the matrix $\mathrm{M}_b$ vanishes identically.
\end{enumerate}

The above description implies that for a generic choice of $a\in\P^2_a$ (namely distinct from the points $\tilde{D}_i$ and not on any line $\tilde{\Pi}_{it}$), the only way in which $a\in \P^2_a$ could be in the kernel of $\mathrm{M}_b$ is if $b=D_t$. In particular, fixing $a$, the map $\operatorname{App}_\psi (a,\_) \colon \P^2_b \dashrightarrow \P^2_s$ is a rational map of degree four with a single point of indeterminacy at $D_t$. A simple analysis shows that the generic fiber of this map consists of 12 points. This has the following consequence.

\begin{Theorem}\label{thm:App-fibers} The map $\operatorname{App}_C\colon \Con^{\bfmu}_{\bfnu}(C,T) \dashrightarrow \SP$ is a rational dominant map whose generic fiber consists of exactly $12$ points. This map is everywhere well defined over the space $\Con^{\bfmu}_{\bfnu}(C,T)^0$ of generic connections.
\end{Theorem}

\begin{proof}
 This theorem readily follows from the fact that the map $\P^2_a\times\P^2_b \dashrightarrow \P^2_a\times\P^2_s$ in (\ref{eq:ab-as}) is dominant and generically $12:1$.
 Consider the following diagram
 \begin{equation}\label{diag:App_C}
 \begin{tikzcd}
 \Con^{\bar\mu}_{\bar\nu}\big(\P^1,D\big) \arrow[r, dashrightarrow, "1:1"] \arrow[d, swap, "2:1"]
 &
 \P^2_a \times \P^2_b \arrow[rr, dashrightarrow, "12:1"]
 &&
 \P^2_a \times \P^2_s \arrow[d, "2:1"]
 \\
 \Con^{\bfmu}_{\bfnu}(C,T) \arrow[rrr, dashrightarrow, swap, "12:1"]
 &&&
 \SP .
 \end{tikzcd}
 \end{equation}
 In principle, the bottom arrow, which represents $\operatorname{App}_C$, \textit{should} be generically $12:1$. There is one place where we need to be careful: the first two maps on the top row are not surjective. This could decrease the cardinality of the fibers once we trace preimages from right to left.

 Take a point $[(a,s)] \in \SP$ in the image of $\operatorname{App}_C$. Then either $(a,s)$ or $(s,a)$ is in the image of the map $\operatorname{App}_\infty \times \operatorname{App}_\psi$. However, the involution $\Psi$ acts in such a way that if $\operatorname{App}_\infty \times \operatorname{App}_\psi (\nabla) = (a,s)$ then $\operatorname{App}_\infty \times \operatorname{App}_\psi (\Psi\nabla) = (s,a)$.
 Hence the image of this map is invariant under permuting the factors in $\P^2_a\times\P^2_s$. Thus both $(a,s)$ and $(s,a)$ are in the image of such map. From the above discussion, it follows that each of these two points has 12 preimages in $\P^2_a\times\P^2_b$. It is shown in \cite[Theorem~1.1]{LoraySaito2015} that the image of the map $\Con^{\bar\mu}_{\bar\nu}\big(\P^1,D\big) \dashrightarrow \P^2_a\times\P^2_b$ coincides with the complement of the incidence variety defined by $a_0b_0 + a_1b_1 + a_2b_2 = 0$.
 By computing a particular example, we are able to confirm that generically none of the 24 points coming from our original $[(a,s)]$ lie on the incidence variety.
 Therefore the composition of the horizontal arrows with the vertical arrow on the right gives a map which is invariant under $\Psi$ and generically $24:1$.
 Such map descends to the quotient $\Con^{\bar\mu}_{\bar\nu}\big(\P^1,D\big) / \Psi \dashrightarrow \Con^{\bfmu}_{\bfnu}(C,T)$. The resulting map, which by definition is $\operatorname{App}_C$, is thus generically $12:1$.
\end{proof}

We conclude by showing that the map $\operatorname{App}_C\times\operatorname{Bun}$ is generically injective. This means that, generically, a connection is completely defined by its underlying parabolic bundle and the image of the apparent map $\operatorname{App}_C$. Note however that the domain of this map is four-dimensional while the target space has dimension six.

\begin{Theorem}\label{thm:App-embedding}
 The map
 \[
 \operatorname{App}_C\times\operatorname{Bun} \colon \ \Con^{\bfmu}_{\bfnu}(C,T) \dashrightarrow \SP\times\P^1_z\times\P^1_w
 \]
 is a generically injective.
\end{Theorem}

\begin{proof}
 As usual, let us start by analyzing the corresponding map for connections on $\P^1$.
 Let $A \colon \P^2_a \times \P^2_b \dashrightarrow \SP$ be the map obtained by the composition of the last two horizontal arrows in (\ref{diag:App_C}). Let us denote by $B \colon \P^2_a \times \P^2_b \to \P^1_z\times\P^1_w$ the map $B = \phi \circ \operatorname{Bun}$.
 Consider the map $A\times B \colon \P^2_a \times \P^2_b \dashrightarrow \SP\times\P^1_z\times\P^1_w$.
 In order to prove the theorem, we are going to show that the generic fibers of $A\times B$ consist of two points that are in involution with respect to $\Psi$.
 Therefore, in the quotient $\Con^{\bar\mu}_{\bar\nu}\big(\P^1,D\big) / \Psi \stackrel{\sim}{\dashrightarrow} \Con^{\bfmu}_{\bfnu}(C,T)$, the induced map $\operatorname{App}_C$ will be generically injective.

 Let $\nabla$ be a generic connection defined by $(a,b)\in\P^2_a \times \P^2_b$. Suppose $A(a,b) = [(a,s)] \in \SP$. Let $\phi(b) = (z,w)$, and denote $\bar{b} = \psi(b)$.
 We have the following straightforward constraints for a point $(a',b')$ to be on the fiber of $A\times B$ over the point $([(a,s)], \, (z,w))$. First, $b'$ must equal either $b$ or $\bar{b}$, since the map $\phi$ is $2:1$ and the fiber over $(z,w)$ is precisely $\{b, \bar{b}\}$. For the second, let $s' = \operatorname{App}_\psi(a',b')$, such that $A(a',b') = [(a',s')]$. Since $A(a',b') = A(a,b)$, we must have that either $a'=a$ and $s'=s$, or $a'=s$ and $s'=a$.

 At this point, we have shown that the fibers of $A\times B$ consist of at most four points $(a',b')$ that satisfy $a'\in\{a,s\}$, $b'\in\{b,\bar{b}\}$.
 By design, this map is invariant under $\Psi$. Thus, points that are in involution with respect to $\Psi$ belong to the same fiber.
 According to Remark~\ref{rmk:characterization-psi-nabla}, $\Psi(a,b) = (s, \bar{b})$. Therefore both $(a,b)$ and $(s,\bar{b})$ belong to the fiber.
 We now need to show that, generically, the points $(a, \bar{b})$ and $(s,b)$ do not belong to the fiber.

 Assume that $(a, \bar{b})$ belongs to the same fiber as $(a,b)$, namely, $A(a,\bar{b}) = A(a,b)$.
 Let us consider the matrix $\mathrm{M}_b$ that appears on Proposition~\ref{prop:app_psi} as an element of $\operatorname{PGL}(3, \C)$.
 Since $A(a,b) = [(a, \mathrm{M}_b(a))]$, equality $A(a,\bar{b}) = A(a,b)$ means that $\mathrm{M}_b(a) = \mathrm{M}_{\bar{b}}(a)$.
 Note that because~$\Psi$ is an involution, we must have that the composition $\mathrm{M}_{\bar{b}}\mathrm{M}_b$ is the identity (cf.\ Remark~\ref{T-inverse}). Thus, applying $\mathrm{M}_b$ on the left, we have that
 \begin{equation*}
 \mathrm{M}_b^2(a) = \mathrm{M}_b\mathrm{M}_{\bar{b}}(a) = a.
 \end{equation*}
 This implies that $a$, viewed as a line on $\C^3$, is an invariant linear subspace for the matrix $\mathrm{M}_b^2$. This imposes non-trivial polynomial conditions on the space $\P^2_a\times\P^2_b$, which are only satisfied in a proper subvariety of $\P^2_a\times\P^2_b$.
 The case $A(s,b) = A(a,b)$ is treated in the same way as above (it imposes the same conditions).

 We conclude that on a Zariski open subset of $\P^2_a\times\P^2_b$, the fibers of the map $A\times B$ consists of two points which are permuted by $\Psi$. This implies that $\operatorname{App}_C \times \operatorname{Bun}$ is generically injective.
\end{proof}

\begin{Remark}
 The Zariski closure of the image of $\operatorname{App}_C\times\operatorname{Bun}$ is a codimension 2 subvariety $X \subset \SP\times\P^1_z\times\P^1_w$. Unfortunately, we were unable to compute the polynomial equations that define the variety $X$.
\end{Remark}

\subsection*{Acknowledgements}

Most of the present work was carried out while the second author was a postdoc at IRMAR. He would like to thank the IRMAR and the Universit\'{e} de Rennes for hosting him during this period.
We would like to thank Thiago Fassarella and N\'{e}stor Fernandez Vargas for many valuable discussions on this topic.
We also thank Nicolas Tholozan who helped us to understand the action of $\Phi^{{\rm top}}$ on the symplectic $2$-form on the monodromy side.
We're also thankful to the anonymous referees for providing many suggestions to improve the content and clarity of the text.
F.L.~acknowledges the support of CNRS and the project Foliage ANR-16-CE40-0008.
V.R.~was supported by the grants PAPIIT IN-106217, CONACYT 219722, and the PRESTIGE post-doc program (coordinated by Campus France and co-financed under the Marie Curie Actions~- COFUND of the FP7). He also acknowledges the support of the Centre Henri Lebesgue ANR-11-LABX-0020-01.

\pdfbookmark[1]{References}{ref}
\LastPageEnding

\end{document}